\documentclass[11pt,a4paper]{article}
 \UseRawInputEncoding
\usepackage{epsfig}
\usepackage{graphics}
\usepackage{graphicx}
\usepackage{cite,indentfirst}
\usepackage{amsmath}    % Les symboles les plus fr�quents.
\usepackage{amssymb}    % Des symboles.
\usepackage{amsfonts}   % Des fontes, eg pour \mathbb.
\usepackage{verbatim}   % Pour les codes sources en informatique.
\usepackage{mathrsfs}   % Des lettres majuscules cursives (\mathscr).
\usepackage{dsfont}
\usepackage{vmargin}        % R�gler la taille de la feuille.

\makeatletter
\@addtoreset{equation}{section}
\makeatother

\renewcommand{\ln}{\log}

\def\sg{\sigma}

\def\eps{\varepsilon}

\def\tr{\operatorname{tr}}

\let\tend=\rightarrow

\newtheorem{theorem}{Theorem}[section]
\newtheorem{prop}{Proposition}[section]
\newtheorem{cor}{Corollary}[section]
\newtheorem{definition}{Definition}[section]
\newtheorem{conjecture}{Conjecture}[section]

   \newtheorem{rem}{Remark}[section]
  \newtheorem{lemme}{Lemma}
\def\Proof{\medskip\noindent {\it Proof --- \ }}

\let\qed=\cqfd

\newcommand\beq{\begin{equation}}
\newcommand\enq{\end{equation}}
\newcommand\bem{\begin{multline}}
\newcommand\enm{\end{multline}}

\def\beqa{\begin{eqnarray}}
\def\eeqa{\end{eqnarray}}
\def\ba{\begin{array}}
\def\ea{\end{array}}

\def\a{\alpha}

\def\eps{\epsilon}
\def\ga{\gamma}
\def\la{\lambda}
\def\sg{\sigma}

\newcommand{\f}[2]{{\ensuremath{%
    \mathchoice%
    {\dfrac{#1}{#2}}
    {\dfrac{#1}{#2}}
    {\frac{#1}{#2}}
    {\frac{#1}{#2}}
}}}
\newcommand{\tf}[2]{\ensuremath{#1/#2}}
\newcommand{\pa}[1]{\ensuremath{\left(#1\right)}}
\newcommand{\paa}[1]{\ensuremath{\left\{#1\right\}}}
\newcommand{\pac}[1]{\ensuremath{\left[#1\right]}}
\newcommand{\paf}[2]{\ensuremath{\left(\f{#1}{#2}\right)}}
\newcommand{\pab}[2]{\ensuremath{\pa{\ba{c} #1 \\ #2 \ea }}}

\def\eps{\epsilon}
\def\ga{\gamma}
\def\Ga{\Gamma}
\def\de{\delta}
\def\De{\Delta}
\def\la{\lambda}
\def\sg{\sigma}
\def\Sg{\Sigma}

\def\th{\theta}

\newcommand{\mc}[1]{\ensuremath{\mathcal{#1}}}

\newcommand{\ov}[1]{\ensuremath{\overline{#1}}}
\newcommand{\wt}[1]{\ensuremath{\widetilde{#1}}}
\newcommand{\msc}[1]{\ensuremath{\mathscr{#1}}}
\newcommand{\wh}[1]{\ensuremath{\widehat{#1}}}

\newcommand{\Int}[2]{\ensuremath{\int\limits_{#1}^{#2}}}

\newcommand{\sul}[2]{\ensuremath{\sum\limits_{#1}^{#2}}}
\newcommand{\pl}[2]{\ensuremath{\prod\limits_{#1}^{#2}}}

\newcommand{\R}{\ensuremath{\mathbb{R}}}
\newcommand{\Cx}{\ensuremath{\mathbb{C}}}

\newcommand{\Dp}[1]{\ensuremath{\partial_{#1}}}

\newcommand{\limit}[2]{\ensuremath{\underset{#1 \tend #2}{\longrightarrow} }}

\newcommand{\ex}[1]{\ensuremath{\e{e}^{#1}}}

\newcommand{\ddet}[2]{\ensuremath{\text{det}_{#1}\pac{#2}}}

\newcommand{\abs}[1]{\ensuremath{\left| #1 \right|}}
\newcommand{\norm}[1]{\ensuremath{\abs{\abs{#1}} }}

\newcommand{\dd}{\text{d}}
\newcommand{\e}[1]{\ensuremath{\mathrm{#1}}}

\newcommand{\intff}[2]{\ensuremath{\left [ \, #1 \,; #2 \, \right ] }}

\newcommand{\intoo}[2]{\ensuremath{\left ] \, #1 \,; #2 \, \right [ }}

\newcommand{\intn}[2]{\ensuremath{[\![ \, #1 \,;\, #2 \,]\!]}}

\def\ua{\uparrow}
\def\da{\downarrow}

\begin{document}

%%%%%%%%%%%%%%%%%%%%%%%%%%%%%%%%%%%%%%%%%%%%%%%%%%%%%%%%%%%%%%%%%%%%%%%%%%%%%%%%%

%\begin{flushright}
%LPENSL-TH-11/08\\
%\end{flushright}
%\par \vskip .1in \noindent

\vspace{24pt}

\begin{center}
\begin{LARGE}
{\bf  Truncated Wiener-Hopf operators with Fisher-Hartwig singularities.}
\end{LARGE}

\vspace{50pt}

\begin{large}

{\bf K.~K.~Kozlowski}\footnote[1]{Universit\'{e} de Bourgogne, Institut de Math\'{e}matiques de Bourgogne, UMR 5584 du CNRS, France,
karol.kozlowski@u-bourgogne.fr},~~
\par

\end{large}

\end{center}

\vspace{80pt}

\centerline{\bf Abstract} \vspace{1cm}
\parbox{12cm}{\small%
 }
We derive the asymptotic behavior of
determinants of truncated Wiener-Hopf operators generated by symbols having Fisher-Hartwig singularities.
This task is achieved thanks to an asymptotic resolution of the Riemann-Hilbert problem associated to some
generalized sine kernel. As a byproduct, we give yet another derivation of the asymptotic behavior of Toeplitz determinants having
Fisher-Hartwig singularities. The Riemann-Hilbert problem approach to these asymptotics yields a systematic although quickly
cumbersome way to compute their sub-leading asymptotics.

\vspace{5mm}


\vspace{5mm} 


\newpage

\section{Introduction}
There is a long history of understanding the asymptotic behavior of determinants of structured matrices.
It started in 1915 with the seminal work of Szeg\"{o} on the asymptotic behavior of large size Toeplitz matrices. His
result is known today as the strong Szeg\"{o} theorem. It states that for
non-vanishing and regular functions $b\in \mc{C}^{1}\pa{\mathbb{T, \R^+}}$
\beq
\ddet{m}{T\pac{b}}\underset{m \tend +\infty}{\sim} \pa{G\pac{b}}^m E\pac{b} \quad \e{with}\quad T_{jk}\pac{b}=c_{j-k}
\quad \e{and} \quad  c_k= \Int{0}{2\pi}  \f{\dd \th}{2\pi} \ex{i k \th}  b\pa{\th}
\enq
\noindent The constants $G\pac{b}$ and $E\pac{b}$ are expressed in terms of the Fourier coefficients of  $\ln b$:
\beq
G\pac{b}=\ex{ \pac{\ln b}_o} \qquad E\pac{b}=\ex{\sul{k=1}{+\infty} k \pac{\ln b}_k\pac{\ln b}_{-k}} \qquad
\pac{\ln b}_k=\Int{0}{2\pi} \f{\dd \th}{2\pi} \ex{ik \th} \ln b\pa{\th} \; .
\enq
The result of Szeg\"{o} underwent several refinements. In particular
Baxter \cite{BaxterToeplitzStrongSzego}, Hirschman \cite{HirschmanSzegolimittheoremwithweakerHypothesis} and, finally,
Ibragimov\cite{IbragimovFinalFormulationOfStrongSzego}
successively weakened Szeg\"{o}'s original assumptions on $b$.  Furthermore, Widom\cite{WidomSzegoLimitforBlockToeplitz}
considered determinants of block Toeplitz matrices and provided a clear interpretation of the constant $E\pac{b}$ in terms of an operator determinant. Despite the possibility to consider matrix valued functions $b$ in the Szeg\"{o}
theorem, there are limitations of its applicability, even in the scalar case.
Indeed, the theorem already breaks down in the case of symbols  having zeros,
power-law singularities or even jump discontinuities on the unit circle $\msc{C}$. Such symbols can be represented as
\beqa
\sg\pa{\th}&=&b\pa{\th} \pl{p=1}{n} \omega_{\de_p,\ga_p}\pa{\tf{\ex{i\th}}{a_p}}
\quad  \e{with} \; a_p \in \msc{C} \; ,\\
\e{and} \quad\omega_{\de_p,\ga_p}\pa{\ex{i\th}}&=& \f{\ex{i \pa{\th-\pi \e{sgn} \th }\de_p }}{\pa{2-2\cos \th}^{\ga_p}} \;\;
, \; \th \in \intoo{-\pi}{\pi} \; .
\eeqa
In such a representation one assumes that the function $b$ is regular enough, non vanishing on $\msc{C}$ and has a vanishing winding number.
The conjecture about the asymptotic behavior of Toeplitz matrices generated by such
symbols goes back to Fisher and Hartwig in 1968 \cite{FischerHartwigTheConjecture}. More precisely, they
claimed that
\beq
\ddet{m}{T\pa{\sg}} \underset{m \tend +\infty}{\sim}
\pa{G\pac{b}}^m m^{\sul{i=1}{n} \ga_i^2-\de_i^2 } C\pa{\pac{b},\paa{\de}_1^n,\paa{\ga}_1^n}
\enq
The Fisher-Hartwig conjecture was indorsed by a couple of examples where the authors were able to
compute the Toeplitz determinants explicitly. The value of the constant was later conjectured to be equal to
\bem
C\pa{\pac{b},\paa{\de}_1^n,\paa{\ga}_1^n} = E\pac{b} \pl{i=1}{n}
\f{G\pa{1-\ga_i-\de_i}G\pa{1-\ga_i+\de_i}}{G\pa{1-2\ga_i}}  \\
\times \pl{i=1}{n} b^{\ga_i-\de_i}_{-}\pa{a_i} b_{+}^{\ga_i+\de_i}\pa{a_i}
\pl{p\not=q}{} \pa{1-\tf{a_p}{a_q}}^{\pa{\de_p+\ga_p}\pa{\de_q-\ga_q}} \; .
\end{multline}
There $G$ is the Barnes' function and $b_{\pm}$ are the Wiener-Hopf factors of $b$, \textit{ie}
$b=b_+ G\pac{b} b_- $  with $b_+$, resp. $b_-$, a holomorphic function on the interior, resp.
exterior, of the unit disk and such that $b_{+}\pa{z=0}=1$, resp. $b_{-}\pa{z} \limit{z}{\infty}1$.

The above conjecture was first proved for some particular cases of the parameters $\nu_k$ and $\ov{\nu}_k$.
Basor \cite{BasorLocalizationThmForToeplitz} and, independently,  B\"{o}ttcher \cite{BottcherFHToeplitzPureJump}
treated the case of several jump discontinuities ($\forall \, p\; ,   \ga_p=0  $) under the restriction $\abs{\Re\pa{\de_p}}<\tf{1}{2}$.  In 1985, B\"{o}ttcher and Silbermann \cite{BottcherSilbProofOfFHConjInSomePartCases} proved
the conjecture in the case $\abs{\Re\pa{\de_p}}<\tf{1}{2}$ and $\abs{\Re\pa{\ga_p}}<\tf{1}{2}$. Finally, Ehrhardt
and Silbermann \cite{EhrhardtSilbermannOnePureFHSingularity} proved the conjecture in the case of
a single Fisher-Hartwig type singularity
for all ranges of parameters $\de_p$ and $\ga_p$ where it made sense. This allowed to
prove the conjecture
in most of the cases involving multiple Fisher-Hartwig singularities \cite{EhrhardtAsymptoticBehaviorOfFischerHartwigToeplitzGeneralCase}.
The proof was based on
the so-called separation technique developed by Basor \cite{BasorLocalizationThmForToeplitz}.
The Fisher--Hartwig conjecture breaks down in the case of the so-called ambiguous symbols. Basor and Tracy \cite{BasorTracyGeneralizedFischer-HartwigConjecture} raised a generalized Fisher-Hartwig conjecture for the behavior of Toeplitz determinants generated by such ambiguous symbols.
This conjecture was proven recently in the
framework of Riemann--Hilbert problems for orthogonal polynomials by Deift, Its and Krasovsky \cite{DeiftItsKrasovskyProofOfGeneralizedFHConjectureAndMoreAnnouncmentResults,DeiftItsKrasovskyAsymptoticsofToeplitsHankelWithFHSymbols}.

\vspace{3mm}

There exists a continuous analog of Toeplitz determinants, the Fredholm determinant of truncated Wiener-Hopf
operators. The underlying operators act on functions
$g\in L^{2}\pa{\R}$ according to the formula
\beq
\pa{1+K}.g\pa{t}= g\pa{t}+ \Int{0}{x} \dd t' \, K\pa{t-t'} g\pa{t'} \dd t' \;\; .
\enq
The kernel K is traditionally defined in terms of its Fourier transform $\sg-1$, \textit{ie}
$K\pa{t}=\mc{F}^{-1}\pac{\sg-1}\pa{t}$. In the following, we choose the below convention for the restriction of
the Fourier transform (and of its inverse) to $L^{2}\pa{\R}\cap  L^{1}\pa{\R}$:
\beq
\mc{F}^{-1}\pac{g}\pa{t}= \Int{\R}{} \f{\dd t}{2\pi} \, g\pa{ \xi  } \ex{-it \xi}  \quad \e{and} \quad
\mc{F}\pac{h}\pa{\xi} = \Int{\R}{} \dd t \, h\pa{t} \ex{it \xi} \;\;.
\enq
So that, for functions $\pa{\sg-1} \in L^{1}\pa{\R}$, we have an explicit integral representation
\beq
K\pa{t}= \f{1}{2\pi} \Int{\R}{}\!\! \dd \xi \pac{\sg\pa{\xi}-1} \ex{-it \xi} \;\; .
\enq
The question of the $x \tend +\infty$ asymptotics of $\ddet{}{I+K}$ were first addressed by Achiezer \cite{AchiezerKacFormulaforTruncatedWienerHopf}
and Kac \cite{KacAcheizerTruncatedWienerHopf}. More precisely, they showed
that for symbols $\sg$ regular enough
\beq
\ddet{}{I+K} \underset{x \tend +\infty }{ \sim} \exp\paa{x \Int{\R}{} \ln \pac{\sg\pa{\xi}} \, \dd \xi
+ \Int{0}{+\infty} \xi \, \ln \pac{\sg\pa{\xi}} \,  \ln \pac{\sg\pa{-\xi}} \;  \dd \xi  } \;\; .
\enq
There exist many generalizations of this formula. These either extend the result to less regular or matrix valued
 symbols $\sg$. However, just as in the Toeplitz case, the theorem  breaks down
when $\sg$ has some jump discontinuities or power-law behavior.
The continuous analogue $\sg$ of a symbol with Fisher-Hartwig type singularities reads
\beq
\sg\pa{\xi}= F\pa{\xi} \pl{k=1}{n} \sg_{\nu_k ,\ov{\nu}_k}\pa{\xi-a_k}  \quad, \quad
\sg_{\nu,\ov{\nu}}\pa{\xi} =
\paf{\xi+i}{\xi+i0^+}^{\nu }
\paf{\xi-i}{\xi-i0^+}^{\ov{\nu}} \;\; .
\label{definition de sigma}
\enq
Note that the exponents in the definition of $\sg_{\nu,\ov{\nu}}$ should be understood in the sense of the principal branch of the logarithm, \textit{ie} $\e{arg} \in \intoo{-\pi}{\pi}$.
In the above decomposition for $\sg$, the function $F$ is supposed regular enough and
$\sg_{\nu_k,\ov{\nu}_k}\pa{\xi}$ has a singularity at $0$:
\beqa
\sg_{\nu_k, \ov{\nu}_k}\pa{\xi} \underset{\xi \tend 0}{\sim} \f{\ex{i\pi \de_k \e{sgn}\pa{\xi}}}{ \abs{\xi}^{2\ga_k} } \quad
\e{with}\;\;\;
2 \ga_k =\nu_k + \ov{\nu}_k  \;\; \e{and} \;\;
2 \de_k=\ov{\nu}_k- \nu_k \; .
\label{behavior of a single sg nu nubar at zero}
\eeqa
There exists a continuous analogue of the Fisher-Hartwig conjecture for Toeplitz determinants and
it is due to B\"{o}ttcher \cite{BottcherFHConjectureWienerHopf}. It was
inspired by the study of truncated Wiener-Hopf operators that are generated by rational symbols; indeed, in the latter case, the author was able to estimate the Fredholm determinants explicitly. This conjecture was confirmed
in many particular cases: the case \cite{BottcherSilberWidomProofFHConjWienerHopfDiscontinuity} of pure jump type singularities,
the case \cite{BottcherSilbWienerHopfAsympAllNuorBarNuZero} where all $\nu_k=0$ or all $\ov{\nu}_k=0$  and finally in the case \cite{BasorWidomWienreHopfwithOneFHSingularity} of a single pure Fisher-Hartwig
singularity $\sg_{\nu,\ov{\nu}}$ under the restriction $\abs{\Re\pa{\nu \pm \ov{\nu}}}<1$.
All these results were established thanks to some identity relating determinants of truncated Wiener-Hopf
operators to determinants of Toeplitz matrices, and then the use of the formulae for the asymptotic behavior of
Toeplitz determinants.

Yet there exists an alternative
approach to the asymptotic analysis of large size determinants of structured matrices. In true, it is well
known that a Hankel matrix can be expressed as a product of the leading coefficients of the orthogonal
polynomials with respect to the measure defining its entries. It was observed by Fokas, Its
and Kitaev \cite{FokasItsKitaevIsomonodromyPlusRHPforOrthPly} that one can recast the problem
of computing orthogonal polynomials into a certain Riemann-Hilbert Problem (RHP). This problem can be solved
asymptotically  for a large case of weights \cite{DeiftKriechMcLaughVenakZhouOrthogonalPlyExponWeights,DeiftKriechMcLaughVenakZhouOrthogonalPlyVaryingExponWeights,
KuilajaarsMVVUniformAsymptoticsForModifiedJacobiOrthogonalPolynomials}.
As noticed by Krasovsky\cite{KrasovskyHankelDetAsymForPowerLike}, one can relate the
asymptotic solution of the RHP for orthogonal polynomials with respect to a weight having a finite number
of power-like singularities to the asymptotic behavior of the Hankel determinant defined in terms of this weight. Its and Krasovsky
used an analogous identity to establish the asymptotic behavior of a Hankel determinant
defined in terms of a gaussian weight having a jump discontinuity \cite{ItsKrasHankelDetAsymForJumpSing}.
 Moreover, still in the framework of RHP for orthogonal polynomials, Krasovsky estimated the asymptotic behavior of Toeplitz matrices on an arc and generated by symbols having jump type discontinuities \cite{KrasovskyToeplitzFHtypeOnArc} thanks to the relationship between polynomials orthogonal
on an arc and those on a line segment. His approach presented no obstruction for a generalization to the
case of root type singularities. The case of Toeplitz, Hankel and Hankel+Toeplitz determinants
with Fisher-Hartwig singularities has been treated recently in the framework of the Riemann--Hilbert problem for orthogonal polynomials by Deift, Its and Krasovsky \cite{DeiftItsKrasovskyAsymptoticsofToeplitsHankelWithFHSymbols,DeiftItsKrasovskyProofOfGeneralizedFHConjectureAndMoreAnnouncmentResults}

Truncated Wiener-Hopf operators being continuous analogs of Toeplitz matrices,
there arises a natural question concerning a RHP approach to study the large $x$ behavior of Fredholm determinants
for such operators. In this paper we will show how to tackle, in the framework of RHP, the $x\tend +\infty$ asymptotics of determinants of truncated Wiener-Hopf operators generated by Fisher-Hartwig symbols. This treatment is based on a relationship between truncated Wiener-Hopf operators and the so-called
generalized sine kernel acting on $\R$. The latter kernel is an integrable integral operator. Such operator can be analyzed
by a RHP as observed in \cite{ItsIzerginKorepinSlavnovDifferentialeqnsforCorrelationfunctions}. We asymptotically solve this RHP.  The construction of its asymptotic solution is an extension of the work \cite{KozKitMailSlaTerRHPapproachtoSuperSineKernel}, the latter being itself a generalization of an unpublished study on the pure sine kernel by Deift, Its and Zhou.
This approximate resolvent allows to compute the leading asymptotics of the Fredholm determinants of truncated Wiener-Hopf operators.
The latter constitutes the main result of this article:
\begin{theorem}
Let I+K be a truncated Wiener-Hopf operator acting on the segment $\intff{0}{x}$ and generated by the symbol
$\pa{\sg-1}\in L^{2}\pa{\R}$ with
\beq
\sg\pa{\xi}=F\pa{\xi} \pl{k=1}{n}\sg_{\nu_k,\ov{\nu}_k}\pa{\xi-a_p} \;  \qquad a_i \in \R\, , \quad a_1<\dots < a_n
\enq
\noindent where
\begin{itemize}

\item $F$ is non-vanishing and holomorphic in some open neighborhood U of the real axis ;

\item  $F-1\in L^{2}\pa{\R}$ and moreover $F\pa{\xi}-1=\e{O}\pa{\abs{\xi}^{-\f{1+\kappa}{2}}}$, for some $\kappa>0$ and $\xi \tend +\infty$ in $U$;

\item $\abs{\Re\pa{\de_k}}<\tf{1}{2}$, and $\Re\pa{\ga_k}<\tf{1}{4}$ for $k\in \intn{1}{n}$.

\end{itemize}
Then the leading asymptotics of the 2-regularized determinant of $I+K$ read:
\bem
\ddet{2}{I+K}=G^x_2\pac{\sg} \cdot  \paf{x}{2}^{\sul{p=1}{n} \ga^2_p-\de_p^2} E\pac{F} \pl{k=1}{n} \f{G\pa{1+\de_k-\ga_k}G\pa{1-\de_k-\ga_k} }{ G\pa{1-2\ga_k} }
\\
\pl{k=1}{n} \f{F^{\ov{\nu}_k}_{ + }\pa{a_k}}{F^{\ov{\nu}_k}_{+}\pa{a_k+i}}
\f{F^{\nu_k}_{-}\pa{a_k}}{F^{\nu_k}_{-}\pa{a_k-i}}  \pl{k\not=p}{n} \paf{\pa{a_k-a_p+i}^2 }{ \pa{a_{k}-a_p+2i}\pa{a_k-a_p}  }^{\ov{\nu}_k \nu_p}\; \pa{1+\e{O}\pa{x^{\rho-1}}} \; .
\end{multline}
There
\beqa
G^x_2\pac{\sg}&=&\exp\paa{x \Int{\R}{} \f{\dd \xi}{2\pi} \pac{\ln\pac{\sg\pa{\xi}}+1-\sg\pa{\xi}  }} \; , \\
E\pac{F} &=& \exp\paa{\Int{0}{+\infty} \dd \xi  \; \xi \,   \mc{F}^{-1}\pac{\ln F}\pa{\xi} \, \mc{F}^{-1}\pac{\ln F}\pa{-\xi}  } \; ,
\eeqa
and $F_{\pm}$ are the Wiener--Hopf factors of $F$:
\beq
\log F_{\pm}\pa{z} = \Int{\R}{} \f{\dd \xi }{2i\pi} \f{ \log F\pa{\xi} }{z-\xi}\;.
\enq
The estimates for the correction involve the constant $\rho=2\underset{k}{\max}{\abs{\Re\pa{\de_k}}}<1$. Finally G stands for the Barnes G function and we remind that
$\ddet{2}{I+K}=\ddet{}{\pa{I+K}\ex{-K}}$.
\end{theorem}

Note that a similar result can also be established for symbols $\sg-1 \in L^{1}\pa{\R}$.
The above theorem reproduces all of the aforementioned results that were obtained for particular cases of singularities.
However, it shows that the original conjecture doesn't hold in its whole generality. Indeed the conjecture \cite{BottcherFHConjectureWienerHopf,BottcherSilAnalysisOfToeplitzOperators} predicts
the presence of
\beq
\pl{k<p}{n}\paf{\pa{\pa{a_k-a_p}^2+1}^2 }{ \pa{\pa{a_{k}-a_p}^2+4}\pa{a_k-a_p}^2  }^{\ov{\nu}_k \nu_p} \;\; ,
\enq
whereas we find
\beq
\pl{k<p}{n} \paf{\pa{a_k-a_p+i}^2 }{ \pa{a_{k}-a_p+2i}\pa{a_k-a_p}  }^{\ov{\nu}_k \nu_p}
\pl{k<p}{n} \paf{\pa{a_k-a_p-i}^2 }{ \pa{a_{k}-a_p-2i}\pa{a_k-a_p}  }^{\ov{\nu}_p \nu_k}  \; .
\enq
 Of course both results coincide in all the cases previously investigated $\textit{ie}$ $\nu_k=\pm \ov{\nu}_k$,   $\forall k \; , \nu_k=0$ or $\forall k \; , \ov{\nu}_k=0 $.
  Moreover our approach opens a way, at least in principle, to asymptotically
inverting truncated Wiener-Hopf operators generated by Fisher--Hartwig symbols. Such an inversion could  be carried
out in the spirit of the inversion for holomorphic symbols $\sg$ proposed in \cite{KozKitMailSlaTerRHPapproachtoSuperSineKernel}.

As a byproduct, using a formula \cite{DeiftItsZhouSineKernelOnUnionOfIntervals} relating
certain  Fredholm determinants of a generalized sine kernel to Toeplitz matrices, we derive the asymptotic behavior of Toeplitz
determinants generated by symbols having Fisher-Hartwig singularities. This computation reproduces the result of
Ehrhardt \cite{EhrhardtAsymptoticBehaviorOfFischerHartwigToeplitzGeneralCase} for the leading asymptotics and of
Deift, Its and Krasovsky for the sub-leading ones \cite{DeiftItsKrasovskyProofOfGeneralizedFHConjectureAndMoreAnnouncmentResults}.
Hence, we see that Toeplitz determinants and those of truncated Wiener-Hopf operators are both related to a determinant
of a generalized sine kernel. The only difference being the interval on which the generalized sine kernel acts.

The article is organized as follows. In the first Section we establish the link between
Fredholm determinants of truncated Wiener-Hopf operators and those of a generalized sine kernel. The second part of Section \ref{Les acteurs se placent sur la scene}
is devoted to introducing some notations. In particular, we shall consider two
types of symbols for $K$. The first one belongs to $L^2\pa{\R}$ whereas the second to $L^{1}\pa{\R}$. Both have the
same Fisher-Hartwig singularities, but in the second case, part of the function $F$  \eqref{definition de sigma}
contains some prefactor depending on $\de$ in order to ensure that $\pa{\sg-1} \in L^{1}\pa{\R}$.

In Section  \ref{SectionRHP} and \ref{Section Parametrix} we asymptotically solve the Riemann--Hilbert problem
associated with the generalized sine kernel. In Section \ref{Section Estimates for the Resolvent} we derive the
asymptotic behavior of the resolvent of the generalized sine kernel. In Section \ref{Section asymptotiques du determinant} we use this asymptotic resolvent to obtain the asymptotic behavior of the Fredholm determinants of truncated Wiener-Hopf operators under investigation.
We compare our results with the already existing ones.
In the  last Section we adapt the RHP associated to the generalized sine kernel in order to study the asymptotic behavior
of Toeplitz matrices with Fisher-Hartwig symbols. We obtain the leading asymptotics of such Toeplitz determinants
and also compute their first sub-leading corrections. The structure of these corrections indicates that at least
part of the asymptotic series can be deduced by making shifts of certain parameters appearing in the leading asymptotics. This leads us to raise a generalization of the Basor--Tracy conjecture.

\section{Some preliminary results and definitions.}
\label{Les acteurs se placent sur la scene}
In this Section we establish the link between truncated Wiener-Hopf operators and generalized sine kernels
in the case of general  $L^1$ and $L^2$ symbols $\sg-1$. We also present some determinant identities that will be useful in our proofs.

\subsection{Truncated Wiener-Hopf and the modified sine kernel}
Let $I+K$ be a truncated Wiener-Hopf operator acting on $L^{2}\pa{\R}$, \textit{ie}. $I+K$ acts on functions
$g\in L^{2}\pa{\R}$ as follows
\beq
\pa{I+K}.g\pa{t}= g\pa{t}+ \Int{0}{x} \dd t' K\pa{t-t'} g\pa{t'} \dd t'
\label{operateur integral de Wiener-Hopf}
\enq
It is useful to define the kernel K in terms of its Fourier transform $\sg-1$, \textit{ie}
$K\pa{t}=\mc{F}^{-1}\pac{\sg-1}\pa{t}$.

 We shall focus on two cases of interest:  $\sg\pa{\xi}-1 \in L^{1}\pa{\R}$ and
$\sg\pa{\xi}-1 \in L^{2}\pa{\R}$. In the first case, $I+K$ is trace class and hence its determinant is well defined.
In the second one, $I+K$ is Hilbert-Schmidt; hence, one ought to consider the 2-regularized determinant \cite{SimonsInfiniteDimensionalDeterminants}
of $I+K$ \textit{ie} $\ddet{2}{I+K}=\ddet{}{\pa{I+K}\ex{-K}}$. For the purpose of this section only, we introduce  the index $\ell$ in $\sg_{\pa{\ell}}$. This  means that $\pa{\sg_{\pa{\ell}}-1} \in L^{\ell}\pa{\R}$, \textit{ie} $\pa{\sg_{\pa{1}}-1}\in L^{1}\pa{\R}$
 and $\pa{\sg_{\pa{2}}-1}\in L^{2}\pa{\R}$.
\begin{lemme}
\label{Lemme egalite Wiener Hopf Sinus modifie}
Let $K_{\pa{\ell}}=\mc{F}^{-1}\pac{\sg_{\pa{\ell}}-1}$ , then one has the identity
\beq
I+K_{\pa{1}}=\mc{F}^{-1} \circ M \circ D \circ \pa{I+V_{\pa{1}}}\circ D^{-1} \circ M^{-1} \circ \mc{F} \; .
\label{1+K et 1+V}
\enq
Similarly,
\beq
I+K_{\pa{2}}=\mc{F}^{-1} \circ D \circ \pa{I+ V_{\pa{2}} } \circ D^{-1} \circ \mc{F}  \; .\label{1+K et 1+Vtilda}
\enq
\noindent The operators $I+V_{\pa{1}}$, resp. $I+V_{\pa{2}}$, act on  $L^{2}\pa{\R}$ with kernels
\beqa
 V_{\pa{1}}\pa{\xi,\eta} &=&\f{\sqrt{\sg_{\pa{1}}\pa{\xi}-1}\sqrt{\sg_{\pa{1}}\pa{\eta}-1}} {2\pi
 i\pa{\xi-\eta}} \pac{\ex{ix\f{\xi-\eta}{2}}-\ex{-ix\f{\xi-\eta}{2}} } \quad;
\label{SineL1} \\
V_{\pa{2}}\pa{\xi,\eta}&=&\pac{\sg_{\pa{2}}\pa{\xi}-1} \f{\ex{ix\f{\xi-\eta}{2}}-\ex{-ix\f{\xi-\eta}{2}}  }{2i\pi \pa{\xi-\eta}}  \;. \label{SineL2}
\eeqa
\noindent  and $D,\, M$ are the multiplication operator on $L^{2}\pa{\R}$
\beq
\pa{M.g}\pa{\xi}= \sqrt{\sg_{\pa{1}}\pa{\xi}-1} \; g\pa{\xi} \;\;, \hspace{1.5cm}
\pa{D.g}\pa{\xi}=\ex{i\f{x \xi}{2}} \; g\pa{\xi} \; .
\enq
\end{lemme}

\Proof
In any of these two cases,  we have that
\bem
\forall g\in L^{2}\pa{\R} \quad \mc{F}\pac{\Int{0}{x} K_{\pa{p}}\pa{t-t'} g\pa{t'} \dd  t' }\pa{\xi} = \\
\pac{\sg_p\pa{\xi}-1} \Int{\R}{} \dd \eta \mc{F}\pac{g}\pa{\eta} \f{\ex{ix\pa{\xi-\eta}}-1}{2i\pi \pa{\xi-\eta}}
= \pac{O_p \circ \,  V_{\pa{p}} \, \circ\, O_p^{-1}\, \circ \,  \mc{F} } \pac{g} \pa{\xi}  \;\; .
\end{multline}
\noindent Where $O_1=M\circ D $ and $O_2=D$. Here, we precise that $V_{\pa{1}}\circ O^{-1}_1$ is indeed a well defined operator on $L^{2}\pa{\R}$. $\quad \Box$

\vspace{3mm}

These two identities relate the truncated Wiener-Hopf operator to a
generalized sine kernel acting on $\R$. If one is able to construct the resolvent for this operator, then one is
able to invert the corresponding Wiener-Hopf operator by using \eqref{1+K et 1+V} or \eqref{1+K et 1+Vtilda}.
This correspondence has already been used in \cite{KozKitMailSlaTerRHPapproachtoSuperSineKernel} to build the resolvent of
truncated Wiener-Hopf operators whose symbols
are holomorphic, non-vanishing functions on some strip around the real axis that are decaying fast enough
at infinity. In the case of symbols
$\sg_{\ell}$ having Fisher-Hartwig singularities as in \eqref{definition de sigma}, the expression for
the leading asymptotic resolvent of the underlying generalized sine kernel is much more involved than for operators considered in
\cite{KozKitMailSlaTerRHPapproachtoSuperSineKernel}, hence taking the Fourier transform and then obtaining some manageable result
might be complicated.

The two identities given in Lemma \ref{lemme egalite determinant Wiener-Hopf Sinus modifie} allow to
establish a connection between sufficiently regularized Fredholm determinants of $I+K_{\pa{\ell}}$ and those
of the corresponding  generalized sine kernel. Hence, to study the asymptotics of $\ddet{\ell}{I+K_{\pa{\ell}}}$
it is enough to focus on the ones of the associated generalized sine kernels.
\begin{lemme}
\label{lemme egalite determinant Wiener-Hopf Sinus modifie}
Let $K_{\pa{\ell}}=\mc{F}^{-1}\pac{\sg_{\pa{\ell}}-1}$ be the kernel of the integral operator given in \eqref{operateur integral de Wiener-Hopf},
then

\beq
\ddet{2}{I+K_{\pa{2}}}=\ddet{2}{I+V_{\pa{2}}} \; ,
\enq
\noindent and
\beq
\ddet{}{I+K_{\pa{1}}}=\ddet{}{I+V_{\pa{1}}} \; .
\enq
\end{lemme}

\Proof
%The second equality follows from the fact that $V$ is trace class as $\sg-1 \in L^{1}\pa{\R}$ and
%$\tf{\sin x t}{2t}$ is bounded on $\R$.
The second equality can be obtained thanks to the Fredholm's series representation
for the determinant of a trace class integral operator:
\beqa
\ddet{}{I+K_{\pa{1}}}&=&  \sul{n=0}{+\infty} \f{1}{n!} \Int{0}{x} \dd^n t \; \ddet{n}{K_{\pa{1}}\pa{t_i-t_j}} \nonumber \\
&=&  \sul{n=0}{+\infty} \f{1}{n!} \Int{0}{x} \dd^n t \Int{\R}{} \f{\dd^n \xi}{2\pi} \pl{p=1}{n} \pa{\sg_{\pa{1}}\pa{\xi_p}-1}
\ddet{n}{ \ex{-i\xi_p\pa{t_p-t_j}} } \nonumber \\
&=& \sul{n=0}{+\infty} \f{1}{n!} \Int{\R}{} \dd^n \xi \, \ddet{n}{V_{\pa{1}}\pa{\xi_p,\xi_j}}
\eeqa
\noindent Now, let us prove the first identity. Define  $R_2\pa{K_{\pa{2}}}=\pa{I+K_{\pa{2}}}\ex{-K_{\pa{2}}}-I$.
Since $K_2$ is Hilbert-Schmidt, we have that $R_{2}\pa{K_{\pa{2}}}$ is trace class \cite{SimonsInfiniteDimensionalDeterminants}.
Moreover,  equation $\eqref{1+K et 1+Vtilda}$ implies that
\beq
R_{2}\pa{K_{\pa{2}}}= \mc{F}^{-1} \circ D \circ R_{2}\pa{V_{\pa{2}}}\circ D \circ \mc{F} \;
\enq
\noindent and $R_{2}\pa{V_{\pa{2}}}$ is trace class as $V_{\pa{2}}$ is Hilbert-Schmidt. Moreover, the Fourier tansform
$\mc{F}$ its inverse $\mc{F}^{-1}$ as well as $D$ and $D^{-1}$ being continuous operators on $L^{2}\pa{\R}$ we have that $R\pa{V_{\pa{2}}}\circ D^{-1}\circ \mc{F}$ is trace class. We can thus change the order in the operator product
appearing in the determinant so that
\beqa
\ddet{2}{I+K_{\pa{2}}}&=&\ddet{}{I+\mc{F}^{-1}\circ D \circ R_{2}\pa{V_{\pa{2}}}\circ D^{-1}\circ \mc{F}} \nonumber\\
&=&\ddet{}{I+R_{2}\pa{V_{\pa{2}}} \circ D^{-1} \circ \mc{F} \circ \mc{F}^{-1} \circ D}=
\ddet{}{I+R_{2}\pa{V_{\pa{2}}}}  \quad \Box \; .
\nonumber
\eeqa
\subsection{General assumptions}
\label{soussection general assumptions}

Motivated by the latter results, we consider the two Fredholm operators $I+V_{\ell}$ acting on $L^2\pa{\R}$
and defined through equations $\eqref{SineL1}$ and $\eqref{SineL2}$. The operators are defined in terms of
the symbols $\sg_{\ell}$ below
\beq
\sg_{ \pa{\ell} }\pa{\xi} =  F_{ \pa{\ell} }\pa{\xi} \pl{k=1}{n} \sg_{\nu_k ,\ov{\nu}_k}\pa{\xi-a_k} \; .
\label{sigma pour deux valeurs de i}
\enq
\noindent The functions $F_{ \pa{\ell} }\pa{\xi}$ are chosen according to
\beqa
F_{\pa{1} }\pa{\xi}&=& b_{\pa{1}}\pa{\xi} \pl{k=1}{n} \pa{1+\f{2i\de_k}{i+\xi} }  \; \; ,\\
F_{\pa{2}} \pa{\xi}&=& b_{\pa{2} }\pa{\xi} \; \; .
\eeqa
and we assume that the functions $b_{\pa{\ell}}$  are such that
\begin{itemize}
\item $b_{\pa{\ell}}$ is  holomorphic on some open neighborhood $U$ of $\R$;
\item   $b_{\pa{\ell}}$  never vanishes on U ;
\item  $b_{\pa{\ell}}\pa{\xi}-1=\e{O}\pa{\abs{\xi}^{-\f{1+\kappa}{\ell}}} \; $ for some $\kappa>0$, when $\xi \tend \infty$ in $U$.
\end{itemize}
\noindent We also make some assumptions on the exponents $2\de_{\ell}=\ov{\nu}_{\ell}-\nu_{\ell}$ and
$2 \ga_{\ell}=\ov{\nu}_{\ell}+\nu_{\ell}$:
\begin{itemize}
\item $\forall \, k \,,  \; \Re\pa{\ga_k}<\tf{1}{2}$ in the $L^{1}\pa{\R}$ case;
\item $\forall \, k \, ,\;  \Re\pa{\ga_k}<\tf{1}{4}$ in the $L^{2}\pa{\R}$ case;
\item $\forall \, k \, , \;  \abs{\Re\pa{\de_k}}<\tf{1}{2}$  independently of the $L^{1}\pa{\R}$ or the $L^2\pa{\R}$ case.
\end{itemize}
The behavior of $\sg_{\nu, \ov{\nu}}\pa{\xi}$ around $\xi=0$
shows that the last restriction on $\de_k$ covers almost all the possible types of jump singularities $\sg_{\pa{\ell}}\pa{\xi}$ could have.
However, the case $\Re\pa{\de_k}=\pm\tf{1}{2}$ for some $k$'s should be treated separately. In particular, one expects additional corrections to the asymptotic formula \eqref{LA formule de l'article}. These should have the same structure  as those appearing in the Basor-Tracy conjecture for Toeplitz matrices \cite{BasorTracyGeneralizedFischer-HartwigConjecture}.
We have the
\begin{lemme}
Under the above assumptions, the symbols $\sg_{ \pa{\ell} }$ given in \eqref{sigma pour deux valeurs de i}
 are such that $\sg_{\pa{1}}-1 \in L^{1}\pa{\R}$  and $\sg_{\pa{2}}-1\in L^2\pa{\R}$.
\end{lemme}
\Proof
We give the proof in the $L^{1}\pa{\R}$ case only.
The assumptions on the parameters $\ga_k$ and the local behavior of $\sg_{\nu_k,\ov{\nu}_k}\pa{\xi-a_k}$
together with $a_1<\dots<a_n$ enure that $\sg_{\pa{1}}-1 \in L^{1}\pa{ \intff{-M}{M} }$ for any finite $M$.
It remains to check the integrability at infinity. $\sg_{\nu_k,\ov{\nu}_k}$ decreases at
infinity as:
\beq
\sg_{\nu_k,\ov{\nu}_k}\pa{\xi}=1-2i\de_k\xi^{-1} + \e{O}\pa{\xi^{-2}}
\enq
Therefore, we have
\bem
\sg_{\pa{1}}\pa{\xi}=b_{\pa{1}}\pa{\xi}\pa{1+ \f{2i}{\xi}  \sul{k=1}{n} \de_k +\e{O\pa{\xi^{-2}}}} \pl{k=1}{n}
\pa{1-2i\de_k\xi^{-1} + \e{O}\pa{\xi^{-2}}}\\
= b_{\pa{1}}\pa{\xi}-1 + \e{O}\pa{\xi^{-2}} =\e{O}\pa{\abs{\xi}^{1+\kappa}} \; ,
\end{multline}
and the claim follows.  $\quad \Box$

\subsection{The resolvent}
\noindent Let $f^{\pa{p}}_{\pm}$ be the solutions to the integral equations
\beq
f_{\pm}^{\pa{p}}\pa{\xi}+ \Int{\R}{}\pa{\sg_{\pa{p}}\pa{\eta}-1} \f{\sin x \tf{\pa{\xi-\eta}}{2} }{\pi \pa{\xi-\eta}}f^{\pa{p}}_{\pm}\pa{\eta} \dd \eta
= \ex{\pm i x \f{\xi}{2}} \; \; .
\label{eqn integral pour f plus moins}
\enq
\noindent It is well known \cite{ItsIzerginKorepinSlavnovDifferentialeqnsforCorrelationfunctions} that the resolvent  operator $I-R_{\pa{1}}$, resp. $I-R_{\pa{2}}$, of $I+V_{\pa{1}}$, resp. $I+V_{\pa{2}}$, has a simple
expression in terms of $f^{\pa{p}}_{\pm}$. Indeed
\beqa
R_{\pa{1}}\pa{\xi,\eta}&=& \f{\sqrt{\sg_{\pa{1}}\pa{\xi}-1}\sqrt{\sg_{\pa{1}}\pa{\eta}-1}}{2i\pi \pa{\xi-\eta}}
\pac{f^{\pa{1}}_{+}\pa{\xi}f^{\pa{1}}_{-}\pa{\eta}-f^{\pa{1}}_{+}\pa{\eta}f^{\pa{1}}_{-}\pa{\xi}} \; ,\\
R_{\pa{2}}\pa{\xi,\eta}&=& \f{ \sg_{\pa{2}}\pa{\xi}-1 }{2i\pi \pa{\xi-\eta}}
\pac{f^{\pa{2}}_{+}\pa{\xi}f^{\pa{2}}_{-}\pa{\eta}-f^{\pa{2}}_{+}\pa{\eta}f^{\pa{2}}_{-}\pa{\xi}} \; .
\label{Resolvent en termes de f_+ et f_-}
\eeqa
\begin{lemme}
Suppose that $\ddet{}{I+V_{\pa{1}}}\not=0$, resp. $\ddet{2}{I+V_{\pa{2}}}\not=0$, then the solutions $f^{\pa{p}}_{\pm}$ of
\eqref{eqn integral pour f plus moins} are entire functions.
\end{lemme}
\Proof
Suppose that $\ddet{}{I+V_{\pa{1}}}\not=0$. Then $I+V_{\pa{1}} :L^{2}\pa{\R} \tend L^{2}\pa{\R} $ is invertible and its inverse
$I-R_{\pa{1}}:L^{2}\pa{\R} \tend L^{2}\pa{\R}$ can be constructed in terms of a Fredholm series.
Thus, since $ \ex{\pm i x\f{\xi}{2}} \sqrt{\sg_{\pa{1}}\pa{\xi}-1} \in L^{2}\pa{\R}$, the unique solution
 $\wt{f}_{\pm}=\sqrt{\sg_{\pa{1}}\pa{\xi}-1} \; f^{\pa{1}}_{\pm }\pa{\xi}$ of
\beq
\wt{f}_{\pm}\pa{\xi}+ \Int{\R}{}V_{\pa{1}}\pa{\xi,\eta} \wt{f}_{\pm}\pa{\eta} \dd \eta
= \ex{\pm i x\f{\xi}{2}} \,  \sqrt{\sg_{\pa{1}}\pa{\xi}-1}
\label{eqn integral pour tilde f plus moins}
\enq
\noindent belongs to $L^{2}\pa{\R}$. We also have that
\beq
 \forall \xi \in \Cx \;\; , \quad \eta \mapsto \sqrt{\sg_{\pa{1}}\pa{\eta}-1} \,
 \f{\sin x \tf{\pa{\xi-\eta}}{2} }{\pi \pa{\xi-\eta}}  \in L^{2}\pa{\R} \;\; \; .
\enq
\noindent Therefore,
\begin{itemize}
\item $\eta \mapsto F\pa{\xi,\eta} \equiv f^{\pa{1}}_{\pm}\pa{\eta}\pa{\sg_{\pa{1}}\pa{\eta}-1} \f{\sin x \tf{\pa{\xi-\eta}}{2} }
{\pi \pa{\xi-\eta}} \in L^{1}\pa{\R}\, $, for all $\xi \in \Cx$;
\item  $\xi \mapsto F\pa{\xi,\eta}$ is entire for almost all $\eta\, $.
\end{itemize}
\beq
\e{Thus, } \hspace{2cm} \xi \mapsto \Int{\R}{} F\pa{\xi, \eta} \dd \eta
\enq
\noindent is an entire function. By \eqref{eqn integral pour f plus moins}, so is $f^{\pa{1}}_{\pm}$.
In the case $\ddet{2}{I+V_2}\not=0$, we have immediately that $\pa{1-\sg_{\pa{2}}\pa{\xi}} f^{\pa{2}}_{\pm} \in L^{1}\pa{\R}$, and the rest
of the proof goes the same. $\Box$

\subsection{Determinant identity}
We shall now derive an important determinant identity. This identity allows to obtain the leading
asymptotics of truncated Wiener-Hopf operators from thoses of its resolvent.
\begin{lemme}
\label{lemme derivee de determinants}
Let $V_{\pa{1}}$ and $V_{\pa{2}}$ be as in \eqref{SineL1} and $\eqref{SineL2}$ and defined in terms of the symbol $\sg_{\ell}$ \eqref{sigma pour deux valeurs de i} subject to the assumptions of subsection \ref{soussection general assumptions}.
We also assume that $\ddet{}{I+V_{\pa{1}}}\not=0$ and $\ddet{2}{I+V_{\pa{2}}}\not=0$. Suppose that
 $\beta_p$ equals $\de_p$ or $\ga_p$, then the following identities hold
\beqa
\Dp{\beta_p} \ln \ddet{}{ I+V_{\pa{1}} }&=&\Int{\R}{}  \f{R_{\pa{1}}\pa{\xi,\xi}}{\sg_{\pa{1}}\pa{\xi}-1} \Dp{\beta_p}\sg_{\pa{1}}\pa{\xi} \dd \xi \;\; , \label{derivee partielle trace class}\\
\Dp{\beta_p} \ln \ddet{2}{ I+V_{\pa{2}} }&=&\Int{\R}{}
\paa{ \f{R_{\pa{2}}\pa{\xi,\xi}}{\sg_{\pa{2}}\pa{\xi}-1}\Dp{\beta_p}\sg_{\pa{2}}\pa{\xi}-\Dp{\beta_p} V_{\pa{2}}\pa{\xi,\xi} } \dd \xi  \;\; . \label{derivee partielle HilbertSchmidt}
\eeqa
\end{lemme}

\Proof
We first treat the $L^{1}\pa{\R}$ case.
Let $\paa{\de^{0}_p,\ga^{0}_p}_{p=1}^{n}$ be a point in $\Cx^{2n}$ fulfilling the assumptions of subsection \eqref{soussection general assumptions} for the $L^1\pa{\R}$ case and such that $\ddet{}{I+V_{\pa{1}}}\not=0$. It then follows from the Fredholm series for $\ddet{}{I+V_{\pa{1}}}$, that the latter is a holomorphic and non-vanishing function on some open neighborhood of $\paa{\de^{0}_p,\ga^{0}_p}_{p=1}^{n}$. It is in particular differentiable and its derivatives can be expressed by using the resolvent operator $I-R_{\pa{1}}$.
Setting $\ex{2G_1\pa{\xi}}=  \sg_{\pa{1}}\pa{\xi}-1$, we get
\begin{multline}
\Dp{\beta_p}\ln \ddet{}{I+V_{\pa{1}}}= \e{tr}\paa{\pa{I-R_{\pa{1}}}\cdot \Dp{\beta_p} V_{\pa{1}} }
 \nonumber \\
= \Int{\R}{} \dd \xi \; \Dp{\beta_p} G_1\pa{\xi}\pac{V_{\pa{1}}\cdot \pa{I-R_{\pa{1}}}}\pa{\xi,\xi}+
\pac{\pa{I-R_{\pa{1}} }\cdot  V_{\pa{1}} }\pa{\xi,\xi} \Dp{\beta_p} G_1\pa{\xi} \nonumber\\
=2\,\Int{\R}{} \dd \xi \; R_{\pa{1}} \pa{\xi,\xi} \, \Dp{\beta_p} G_1\pa{\xi} =
\Int{\R}{} \dd \xi \; R_{\pa{1}} \pa{\xi,\xi} \, \f{\Dp{\beta_p} \sg_{\pa{1}}\pa{\xi}}{\sg_{\pa{1}}\pa{\xi}-1} \; .
\end{multline}
In the intermediary equalities we used
the symmetry of the kernels  as well as the cyclicity of the trace. The $L^2\pa{\R}$ case is proved by density.
Let $\chi_{\eps}$ be the characteristic function of $\intoo{-\eps}{\eps}$. Then $I+V_{2;\eps}$, with $V_{2;\eps}\pa{\xi,\eta} \equiv \chi_{\eps}\pa{\xi} V_{\pa{2}}\pa{\xi,\eta} \chi_{\eps}\pa{\eta}$, is trace class for all $\eps>0$ so that
\beq
\ddet{2}{I+V_{2;\eps}}=\ddet{}{I+V_{2;\eps}}\ex{-\tr{V_{2;\eps}}} \;\; .
\enq
As $\ddet{2}{I+V_{\pa{2}} }\not=0$, $\ddet{2}{I+V_{2;\eps}}\not=0$ for $\eps$ large enough, and hence
$\ddet{}{I+V_{2;\eps}}\not=0$ as well.
One can then apply the results for $L^{1}\pa{\R}$ kernels  for the $\beta_p$ derivative. We get,
\beq
\Dp{\beta_p} \ln \ddet{2}{I+V_{2;\eps}} = \Int{\R}{}\! \dd \la\,  \chi_{\eps}\pa{\la} \f{R_{2;\eps}\pa{\la,\la}}{\sg_2-1} \Dp{\beta_p}\sg_2\pa{\la}- \Int{\R}{} \dd \la \chi_{\eps} \pa{\la} \Dp{\beta_p}V_{\pa{2}}\pa{\la,\la} \; .
\enq
The $\eps\tend +\infty$ limit  in the $rhs$ becomes licit after merging the two integrals into one. $\qquad \Box$.

\vspace{3mm}

The asymptotic solution of the RHP presented in the upcoming sections will allow us to construct approximate in $x$ resolvents of $V_1$
and $V_2$ uniformly in respect to the parameters $\de_p$ and $\ga_p$. It will then remain to use these approximations to compute,
in the large $x$ limit, the integrals appearing in \eqref{derivee partielle HilbertSchmidt} or \eqref{derivee partielle trace class}. Once this is done, it is enough to integrate the result from $\beta_p=0$ to $\beta_p$. Such an integration is
 analogous to the separation technique
in the operator approach to asymptotics of Toeplitz determinants \cite{BasorLocalizationThmForToeplitz}. It is then not a problem to repeatedly apply the procedure
so as to obtain the asymptotics of the determinant. At this stage, it becomes clear why, in our approach, the $L^1\pa{\R}$ case doesn't follow from the $L^2\pa{\R}$ one. As a matter of fact, if we want to keep jump singularities and still have an $L^{1}\pa{\R}$ kernel, we ought to add an additional factor depending on the $\de_p$'s as it was
explicitly done for the function $F_{\pa{1}}$, \textit{cf} \eqref{sigma pour deux valeurs de i}. This modifies the $\de_p$ dependence of the integrand in \eqref{derivee partielle trace class} and hence the integration procedure. The result should also be, in principle modified, but eventually we see that the $L^1\pa{\R}$ case can be obtained from the $L^2\pa{\R}$ one by restricting
correctly the parameters and replacing $F_{\pa{2}}$ by $F_{\pa{2}}$. However, since these are only minor modifications, from now on we focus on the $L^2\pa{\R}$ case. The interested reader will find no problem in adapting the proofs to the $L^1\pa{\R}$ case. Accordingly, from  now on, we drop the $\ell$ subscript labeling $\sg_{\pa{\ell}}$, the kernels $V_{\pa{\ell}}$ and the resolvents $R_{\pa{\ell}}$. We will denote these quantities by $\sg$, $V$ and $R$ and assume that $\sg=\sg_{\pa{2}}$ as defined in \eqref{sigma pour deux valeurs de i}.

\section{The Riemann--Hilbert Problem}
\label{SectionRHP}
We start this Section by introducing a RHP for a matrix $\chi$. This type of RHP is adapted for constructing resolvents of integrable integral operators \cite{ItsIzerginKorepinSlavnovDifferentialeqnsforCorrelationfunctions} such as the generalized sine kernel. We then perform a few transformations of this original RHP, in order
to boil it down to one where the jump matrices will be $I_2+\e{o}\pa{1}$  uniformly away from the points $a_k$
and in respect to the $x \tend +\infty$ limit.
The first step will consist in finding a scalar valued function $\a$ such that $\chi \a^{\sg_3}$ has a jump matrix
with 1 in its lower diagonal entry. Then we deform the original cut.
The jump matrix on the new contour has the desired properties.

\subsection{The initial Riemann-Hilbert problem}
\label{SectionRHPsousSectionRHPchi}
As first observed in \cite{ItsIzerginKorepinSlavnovDifferentialeqnsforCorrelationfunctions },
the problem of finding the resolvent of any integrable integral operator is equivalent a Riemann-Hilbert problem (RHP).
Indeed, let $f_{\pm}$ be the solutions of \eqref{eqn integral pour f plus moins} for the corresponding $\sg$.
Then, there exists \cite{ItsIzerginKorepinSlavnovDifferentialeqnsforCorrelationfunctions} a matrix $\chi$ allowing to reconstruct the solutions $f_{\pm}$ of the integral equation
\eqref{eqn integral pour f plus moins}, and hence the resolvent:
\beq
\pa{\ba{c} f_+\pa{\xi} \\ f_{-}\pa{\xi} \ea }= \chi\pa{\xi}\pa{\ba{c} \ex{i x \f{\xi}{2}} \\ \ex{-i x \f{\xi}{2}} \ea }
 \quad \e{and}\quad \pa{-f_{-}\pa{\xi},f_+\pa{\xi}}=
\pa{-\ex{-i x \f{\xi}{2}}, \ex{i x \f{\xi}{2}} } \chi^{-1}\pa{\xi} \; .
\label{f_+ et f_- en termes de la solution du RHP}
\enq
This matrix $\chi$ solves the RHP:

\begin{itemize}
\item $\chi$ is analytic on   $\mathbb{C}\setminus\R$ \; ;
\item $\forall  k \in \intn{1}{n} \;,$ there exists $M_k\in \e{GL}_2\pa{\Cx}$ such that \vspace{2mm}\newline
 $ \chi = M_k \paa{I_2+g\pa{z} B_k   +\abs{z-a_k}\pa{g\pa{z}+1}  \e{O}\pa{\ba{cc}1 & 1 \\ 1 & 1 \ea } }
\;\; , \; z \tend a_k \;$ \; ;
\item $\chi \underset{z \tend \infty}{\tend} I_2\equiv \pa{ \ba{cc} 1 & 0 \\ 0& 1\ea }$ \; ;
\item $\chi_{+}\pa{z} G\pa{z}=\chi_-\pa{z}  \;\; ; \quad z \in \R \; \; .$
\end{itemize}

The matrix $B\pa{z}$ appearing in the estimates around $a_k$ is a rank one matrix that takes the precise form
\beq
B\pa{z} =  \pa{\ba{cc} -1 & \ex{ix z}\\ -\ex{-ix z} & 1 \ea } \;.
\label{definition matrice Bk}
\enq
The function $g$ reads
\beq
g\pa{z}= \Int{\R}{} \f{\dd s}{2i\pi}  \f{\sg\pa{s}-1}{z-s} \;
\label{definition fonction g}
\enq
and has the local behavior \cite{GakhovBoundaryValueProblems} at $z\tend a_k$:
\beq
 g\pa{z}=\left\{ \ba{cc c}  \e{O}\pa{1}+\e{O}\pa{\pa{z-a_k}^{-2 \ga_k}}  &  \e{for} & \ga_k\not=0 \\
                          \e{O}\pa{\ln \pa{z-a_k}} & \e{for} & \ga_k=0  \ea \right.
\label{comportement local de g}
\enq
The matrices $M_k$ are $a priori$ unknown and will be determined once the solution is known.
What only matters for the solvability is the invertibility of  $M_k$.
Lastly, we adopt the convention that the symbol $\e{O}\pa{M}$ for some matrix $M$ is to be understood entrywise \textit{ie} $\chi=\e{O}\pa{M}$ means that $\chi_{ij}=\e{O}\pa{M_{ij}}$.

\noindent The jump matrix $G$ appearing in the RHP reads
\beq
 G\pa{z}     = \left(\ba{cc}
                   2-\sg\pa{z} & \pa{\sg\pa{z}-1} \ex{ i x z} \\
                   \pa{1-\sg\pa{z}} \ex{-i x z} & \sg\pa{z}
                   \ea \right) \;\; .
\enq
\noindent Finally, $\chi_{+}\pa{t}$ (resp. $\chi_-\pa{t}$) stands for the non-tangential limits of $\chi\pa{z}$
as $z$ approaches a point $t$ of the contour from its $+$ (resp. $-$) side.
\begin{figure}[h]

\begin{center}
\includegraphics{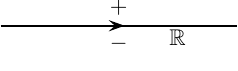}
\end{center}
\caption{Original contour for the RHP.\label{Contour du RHP pour chi}}
\end{figure}

\begin{prop}
\label{proposition unicite et existance RHP}
Whenever $\Re\pa{\ga_k}<\tf{1}{2} \; , \; \forall k \in \intn{1}{n}$ and $\ddet{2}{I+K}\not=0$, the solution to the RHP exists and is unique.
\end{prop}

Although most of the proof is standard, see for example  \cite{KuilajaarsMVVUniformAsymptoticsForModifiedJacobiOrthogonalPolynomials} (careful handling of the singularities) and \cite{DeiftOrthPlyAndRandomMatrixRHP} (general exposure), we include it for the sake of completeness and the reader's convenience. The last part of the proof dealing with a cancelation of the singularities due $\e{rank}\pa{B\pa{z}}=1$ is new.

\Proof

The RHP is equivalent to the singular integral equation for $\chi$:
\beq
\chi\pa{z}=I_2+ \Int{\R}{} \f{\dd s}{2i\pi} \f{\sg\pa{s}-1}{z-s} \chi_+\pa{s}
\pa{\ba{cc } -1 & \ex{i x s}\\ -\ex{- i x s}& 1 \ea}  \;\; , \quad z \in \mathbb{C}\setminus\R \; .
\label{equation integral singuliere pour Chi}
\enq
Since we assume that $\ddet{2}{I+K}\not= 0$, as already mentioned, the resolvent operator $I-R$ exists.
Hence, one can express a solution of equation \eqref{equation integral singuliere pour Chi} in terms of $R$ or,
equivalently, in terms of $f_{\pm}$ \eqref{eqn integral pour f plus moins}
(whose existence follows from the existence of the resolvent):
\beqa
%
%\chi\pa{z}&=&I_2+\Int{\R}{} \f{\dd s}{z-s}  \sqrt{\sg\pa{s}-1}\paa{\sqrt{\sg-1}\pa{\ba{c} e_+ \\ e_- \ea}.\pa{I-R} }\pa{s}
%\otimes \pa{-e_-\pa{s},e_+\pa{s}}
%\nonumber \\
%
%
\chi\pa{z}&=& I_2 + \Int{\R}{} \f{\dd s}{2i\pi} \f{\sg\pa{s}-1}{z-s}
\pa{\ba{cc } - f_+ \pa{s} \ex{- i x \f{s}{2}} & f_+ \pa{s} \ex{i x \f{s}{2}}
            \\ -f_- \pa{s} \ex{- i x \f{s}{2}}& f_-\pa{s} \ex{ i x \f{s}{2}} \ea}  \label{forme explicite de chi} \;\; .
\eeqa
\noindent This proves the existence of solutions provided we show that \eqref{forme explicite de chi} has the desired  behavior around each point $a_k$. The solution of the jump conditions given in \eqref{forme explicite de chi} can be written as
\begin{multline}
\chi\pa{z}=  M\pa{z} + \pa{\ba{cc } -f_+ \pa{z} \ex{- i x \f{z}{2}}
                & f_+ \pa{z} \ex{i x \f{z}{2}} \\
             -f_- \pa{z} \ex{- i x \f{z}{2}}
            & f_-\pa{z} \ex{ i x \f{z}{2}} \ea} \Int{\R}{} \f{\dd s}{2i\pi} \f{\sg\pa{s}-1}{z-s} \;
\end{multline}
with
\bem
M\pa{z}=I_2+ \Int{\R}{} \f{\dd s}{2i\pi} \f{\sg\pa{s}-1}{z-s} \left\{
\pa{\ba{c} f_+\pa{s} \\ f_-\pa{s} \ea  } \pa{-\ex{-ix\f{s}{2}}, \ex{ix\f{s}{2}}}  \right.\\
\left.- \pa{\ba{c} f_+\pa{z} \\ f_-\pa{z} \ea  } \pa{-\ex{-ix\f{z}{2}}, \ex{ix\f{z}{2}}}  \right\} \; .
\nonumber
\end{multline}

The matrix $M\pa{z}$ is holomorphic around  $z=a_k$, and it is easily seen from \eqref{f_+ et f_- en termes de la solution du RHP} that
\beq
\pa{\ba{c} f_+\pa{z} \\ f_-\pa{z} \ea} = M\pa{z} \pa{\ba{c} \ex{ix\f{z}{2}} \\\ex{-ix\f{z}{2}} \ea} \; .
\enq
Therefore,
\beq
\chi\pa{z}=M\pa{z}\paa{I_2+ g\pa{z}  \pa{\ba{cc}-1 & \ex{ix z} \\ \ex{-ix z} & 1  \ea} } \quad \e{with} \quad
g\pa{z}=\Int{\R}{} \f{\dd s}{2i\pi} \f{\sg\pa{s}-1}{z-s} \; .
\enq
The local behavior of $g$ around $z=a_k$ can be inferred from \cite{GakhovBoundaryValueProblems}.
Then, the claim for the local estimates follows by expanding the holomorphic matrices around $z=a_k$ and setting $M\pa{a_k}=M_k$. We shall now prove that $M_k\in \e{GL}_2\pa{\Cx}$

Given any solution  $\chi$ to the above RHP, $\ddet{}{\chi}$ is analytic on $ \mathbb{C}\setminus \R$. Since $\ddet{}{G}=1$
we have that $\ddet{}{\chi}$ is continuous across $\R \setminus\cup_{p=1}^{n}\paa{a_p}$ and can thus
 be extended to an analytic function on $\mathbb{C}\setminus\cup_{p=1}^{n}\paa{a_p}$.
Using the estimates for $\chi$ and the multilinearity of the determinant,  we get that for $z\tend a_k$
\beq
\ddet{}{\chi}=\ddet{}{M_k}\ddet{}{I_2+g\pa{z} B\pa{z}} + \e{O}\pa{g\pa{z}\pa{1+g\pa{z}} \abs{z-a_k}} \; .
\enq
As $B\pa{z}$ is a rank one matrix $\ddet{}{I_2+g\pa{z} B\pa{z}}= g\pa{z} \e{tr}\pac{B\pa{z}}=0$, and thus the first term is a $\e{O}\pa{1}$. The estimates for the local behavior of $g\pa{z}$ together with the hypothesis
on the parameters $\ga_k$ ensure that
$\e{O}\pa{g\pa{z}\pa{1+g\pa{z}} \abs{z-a_k}}=\e{o}\pa{\abs{z-a_k}^{-1}}$. $\ddet{}{\chi}$ has thus no pole at $z=a_k$. Its singularities at the $a_k$'s are thus of a removable type and hence $\ddet{}{\chi}$ is holomorphic around $a_k$.
It follows that $\ddet{}{\chi}$ an entire function that is bounded at
infinity in virtue of the normalization $\chi \limit{z}{\infty} I_2 $.  By Liouville's theorem, $\ddet{}{\chi}=1$.
In particular $\ddet{}{\chi}\pa{a_k}=1$, which can only happen if $\ddet{}{M_k}=1$.

We end the proof by showing the uniqueness of solutions.
Let $\chi_1$ and $\chi_2$ be two solution of the RHP for $\chi$. Then as $\ddet{}{\chi_2}=1$, $\chi_2$ is analytically invertible on $\mathbb{C}\setminus\R$. The matrix $\chi_1\chi^{-1}_2$
is holomorphic on $\mathbb{C}\setminus\R$ and continuous across $\R\setminus \bigcup\limits_{k=1}^{n} \paa{a_k}$.
Moreover, using the local behavior at $z=a_k$ we get
\beq
\chi_2^{-1}\pa{z} =  \pa{ I_2 + ^t\e{Comat}\pa{B\pa{z}} g\pa{z} } M_{k,2}^{-1} + \e{O}\pa{\pa{1+g\pa{z}}\abs{z-a_k} \pa{\ba{cc} 1&1\\1&1\ea}}\, .
\enq
Here, we used the fact that the inverse of $\chi_2$ is given by the transpose of its comatrix due to $\ddet{}{\chi_2}=1$. We have also introduced two, \textit{a priori} distinct, matrices $M_{k,\ell}$ associated with each of the solutions $\chi_{\ell}$.
Computing the matrix products and using that $B\pa{z}+^t\e{Comat}\pa{B\pa{z}}=0$,  $B\pa{z}\cdot ^t\e{Comat}\pa{B\pa{z}}=0$, we get
\beq
\chi_1\pa{z} \chi_2^{-1}\pa{z}=M_{k,1}M_{k,2}^{-1} + \e{O}\pa{g\pa{z}\pa{1+g\pa{z}}\abs{z-a_k} \pa{\ba{cc} 1&1\\1&1\ea}} \;\;\;,
\enq
The local estimates for $g$ imply that $g\pa{z}\pa{1+g\pa{z}}\abs{z-a_k} = \e{o}\pa{\abs{z-a_k}^{-1}}$. Hence $\chi_1\pa{z} \chi_2^{-1}\pa{z}$ has no poles at $z=a_k$. The singularities at the $a_k$'s are thus removable and,  because of the asymptotic condition,
 we have $\chi_1\chi^{-1}_2=I_2 \; $.
This guarantees the uniqueness of the solution to the RHP, at least for $\Re\pa{ \ga_k }<\tf{1}{2}$. $\quad \Box$

Note that one doesn't have to make such fine estimates for proving the uniqueness of solutions if one would assume that $\Re\pa{\ga_k}<\tf{1}{4}$. However, we presented here this more complex proof as it also holds in the
$L^{1}\pa{\R}$ case.

\subsection{A helpful scalar Riemann-Hilbert problem.}
\label{SectionRHPsousSectionRHPalpha}
 Let $\a$ be the solution of the following scalar Riemann-Hilbert problem
\beq
\a \; \e{is}\; \e{analytic} \; \e{on}\;  \Cx\setminus \R \; , \quad  \a_-\pa{\xi}=\a_+\pa{\xi} \sg\pa{\xi},\qquad
 \xi\in \R \setminus
 \bigcup\limits_{k=1}^{n} \paa{a_k}\; ,  \quad \a\pa{\xi} \underset{\xi \tend \infty}{\longrightarrow} 1\; .
\label{RHP for Alpha}
\enq
This Riemann-Hilbert problem can be solved almost explicitly for the particular form  \eqref{definition de sigma} of $\sg$ we study, namely, $\a\pa{z}=\a_{\ua}\pa{z}$ for $z\in \mc{H_+}$ and $\a\pa{z}=\a_{\da}\pa{z}$
for $z\in \mc{H_-}$, where
\beqa
\a_{\ua}\pa{z}&=F_+^{-1}\pa{z} \pl{k=1}{n} \paf{z-a_k}{z-a_k+i}^{\nu_k}
& \quad z \in \mc{H}_+\nonumber  \; \; , \\
\a_{\da}\pa{z}&= F_-\pa{z} \pl{k=1}{n} \paf{z-a_k-i}{z-a_k}^{\ov{\nu}_k}
 &  \quad z \in \mc{H}_-\nonumber \; \; .
\eeqa
\noindent  $\mc{H}_{\pm}$ is the upper/lower half-plane and $F_{\pm}$  are the Wiener-Hopf factors of $F$,
\textit{ie} $F=F_+F_-$  with $F_+$ (resp. $F_-$) analytic in the upper (resp. lower) half-plane and going to
1 at $z\tend \infty$ in $\mc{H_+}$ (resp. $\mc{H}_-$). There is no constant factor $F_0$ in the Wiener-Hopf decomposition of
$F$ as $F\limit{z}{\pm \infty}1$.

The Wiener-Hopf factors of $F$ have an integral representation either in terms of Cauchy or Fourier transforms of $\ln F$:
\beqa
\ln F_+\pa{z} &=  \textrm{{\Large \Int{\R}{}}} \f{ \dd \xi }{2i\pi}   \f{\ln F\pa{\xi}}{\xi-z}=
\mc{F}\pac{\Xi\pa{\xi}\mc{F}^{-1}\pac{\ln F}\pa{\xi} }\pa{z} & \;, \;  z \in \mc{H}_+  \; ,\nonumber \vspace{5mm}\\
\ln F_-\pa{z} &=  -\textrm{{\Large \Int{\R}{}}} \f{ \dd \xi }{2i\pi}   \f{\ln F\pa{\xi}}{\xi-z}=
\mc{F}\pac{\Xi\pa{-\xi}\mc{F}^{-1}\pac{\ln F}\pa{\xi} }\pa{z} & \;, \;   z \in \mc{H}_-   \nonumber \;\;,
\eeqa
and $\Xi$ is Heaviside's step function. As $F$ is analytic and non zero in $U$,  it is clear from these integral representations that $F_+$ and $F_-$ have an  analytic continuation to $U$. Moreover, as $F$ is non-vanishing in $U$
and the decomposition $F=F_+F_-$ is still valid on $U$,
$F_+$ and $F_-$ have no zeroes on $U$.

\noindent We introduce the auxiliary functions
\beqa
\sg_p\pa{z}&=&F\pa{z} \pl{k=1}{p} \paf{z-a_k+i}{z-a_k}^{\nu_k} \paf{z-a_k-i}{z-a_k}^{\ov{\nu}_k} \nonumber\\
&& \hspace{1.5cm} \times \pl{k=p+1}{n} \paa{\paf{z-a_k+i}{a_k-z}^{\nu_k} \paf{z-a_k-i}{a_k-z}^{\ov{\nu}_k} \ex{2i\pi \de_k} }\; ,
\label{sigma continuation sur ap-1 ap}\\
\wh{\sg}_p\pa{z} &=& \pa{z-a_p}^{2\ga_p} \, \sg_p\pa{z} \; .\label{sigma regularise en ap}
\eeqa
\noindent The function\footnote{We stress that the subscript $p$ appearing in $\sg_p$ has nothing to do with the notations of section 2. In the following $\sg_p$ will always refer to the definition \eqref{sigma continuation sur ap-1 ap}} $\sg_p$, resp. $\wh{\sg}_p$, is holomorphic on $\paa{z\; : \; a_p< \Re{z}< a_{p+1}}\cap U$,
resp. $D_{a_p,\eps}=\paa{z\in\Cx \; : \; \abs{z-a_p}<\eps}$. Here and in the following, $\eps$ is such that $D_{a_p,\eps} \subset U$ and $D_{a_p,\eps}\cap D_{a_q , \eps}=\emptyset$ for $p\not=q$.
The functions $\sg_p$ and $\wh{\sg}_p$ can be though of as the analytic parts of the local behavior of $\sg$
 on the real axis. Namely, one can continue $\sg$ by analyticity to the domains below
\beq
\sg\pa{z} = \left\{  \ba{ccl}
        \sg_p\pa{z} &\quad  &z \in \paa{z \, :\;  a_p < \Re\pa{z} < a_{p+1} }\cap U  \vspace{3mm}\\
        \wh{\sg}_{p}\pa{z} \f{\ex{ 2i\pi\ov{\nu}_p\Xi\pa{ \Re\pa{a_p-z}}  } }{\pa{z-a_p}^{2\ga_p}}  & &z \in \mc{H}_+ \cap D_{a_p,\eps}  \vspace{3mm}\\
        \wh{\sg}_{p}\pa{z} \f{\ex{-2i\pi\nu_p\Xi\pa{ \Re\pa{a_p-z}} } }{\pa{z-a_p}^{2\ga_p}}  & &z \in \mc{H}_- \cap D_{a_p,\eps} \ea
\right.
\label{continuation sg sur l'axe app+1}
\enq
In much the same way, we split the formula for $\a$ into a holomorphic and a singular part:
\beqa
\a_{\ua}^2\pa{z} &=& \wh{\sg}^{-1}_p\pa{z}  \pa{z-a_p}^{2\nu_p } x^{-2\de_p} \ex{i x a_p} K_{p}\pa{z}
 \quad , \;\; z \in \mc{H}_+\cap D_{a_p,\eps} \; ;\nonumber \\
\a_{\da}^2\pa{z} &=& \wh{\sg}_p\pa{z}  \pa{z-a_p}^{-2\ov{\nu}_p } x^{-2\de_p }  \ex{i x a_p} K_{p}\pa{z}
 \quad, \;\; z \in \mc{H}_-\cap D_{a_p,\eps} \; . \nonumber
\eeqa
There,
\beq
K_{p}\pa{z}= \f{ x^{2\de_p} }{ \ex{ixa_p} } \f{F_-\pa{z}}{F_+\pa{z}} \pl{k=1}{n}
\f{ \pa{z-a_k-i}^{ \ov{\nu}_k }  }{ \pa{z-a_k+i}^{ \nu_k } }
\pl{k=1}{p-1} \f{ \pa{z-a_k}^{ \nu_k } }{ \pa{z-a_k}^{ \ov{\nu}_k } }
 \pl{k=p+1}{n} \paa{\f{ \pa{a_k-z}^{ \nu_k } }{ \pa{a_k-z}^{ \ov{\nu}_k } } \ex{2i\pi\ga_k } } \label{definition de K_p}
\enq
\noindent The analyticity of $F\pa{z}$ on $U$ guarantees that $K_p\pa{z}$ is holomorphic on the disk $D_{a_p,\eps}$.
Lastly, we define $\wh{\a}_{\ua/\da}^{\pa{p}}$, the regularized version of $\a_{\ua/\da}$ around $a_p$, according to
\beq
\wh{\a}_{\ua}^{\pa{p}}\pa{\xi}=\a_{\ua}\pa{z} \; \pa{z-a_p}^{-\nu_p} \qquad , \;\;
\wh{\a}_{\da}^{\pa{p}}\pa{\xi}=\a_{\da}\pa{z} \; \pa{z-a_p}^{ \ov{\nu}_p} \; .
\label{alpha plusmoins regularise}
\enq
So that,
\beq
\f{\wh{\a}_{\da}^{\pa{p}}\pa{z}}{\wh{\a}_{\ua}^{\pa{p}}\pa{z}}=\wh{\sg}_p\pa{z} \; .
\enq
\subsection{The first step $\chi \tend \Phi$}
\label{SectionRHPsousSectionRHPchiversPhi}
\noindent  Let $\Phi$ be related to $\chi$ by
\beq
\Phi\pa{z} = \chi\pa{z} \pac{\a\pa{z}}^{\sg_3} \;\; , \qquad  \sg_3\equiv  \pa{\ba{c c} 1 & 0 \\ 0 & -1 \ea } \; .
\label{substitution 1 }
\enq
Then  $\Phi\pa{z}$ satisfies the following RHP:
\begin{itemize}
\item $\Phi$  is analytic in $\mathbb{C}\setminus\R\;$;
\item  $\forall  k \in \intn{1}{n} \;,$ there exists $M_k\in \e{GL}_2\pa{\Cx}$ such that:
\vspace{3mm}
\newline $  \Phi =  M_k \paa{I_2+g\pa{z} B\pa{z}
+\abs{z-a_k}\pa{g\pa{z}+1}  \e{O}\pa{\ba{cc}1 & 1 \\ 1 & 1 \ea } }
 \pac{\a\pa{z}}^{\sg_3}$
\item $\Phi \underset{z \tend \infty}{\tend} I_2 \; $ ;
\item $\Phi_{+}\pa{z} G_{\Phi}\pa{z}=\Phi_-\pa{z}  , \quad z \in \R \; ;$
\end{itemize}
The function $g$ and the rank one matrix $B\pa{z}$ are as given in \eqref{definition fonction g} and \eqref{definition matrice Bk}. In particular $g$ has a singular behavior at $z\tend a_k$ given by \eqref{comportement local de g}.
$\a$, the solution of the scalar RHP \eqref{RHP for Alpha}, introduces an additional singular behavior
at $z\tend a_k$:
\beq
 \a\pa{z}=\left\{ \ba{cc c}  \e{O}\pa{\abs{z-a_k}^{\Re\pa{\nu_k}}}  &  \e{for} & z\tend a_k \; ,\;  z\in \mc{H_+} \\
                          \e{O}\pa{\abs{z-a_k}^{\Re\pa{-\ov{\nu}_k}}}  &  \e{for} & z\tend a_k \; ,\;  z\in \mc{H_-}  \ea \right.
\enq
Finally, the jump matrix for $\Phi$ reads
 \beq
 G_\Phi\pa{z}=\left( \ba{cc}
                     1+P\pa{z}Q\pa{z}& P\pa{z}\ex{ix z}\\
                     Q\pa{z}\ex{-ixz}& 1
                    \ea\right),
\label{jump matrix for Phi}
\enq
and
 \beqa
 P\pa{z}&=&\pac{1-\sg^{-1}\pa{z}}\a_{+}^{-2}\pa{z}\; , \vspace{4mm} \\
 Q\pa{z}&=&\pac{\sg^{-1}\pa{z}-1}\a_{-}^{2}\pa{z}\; .
\label{Definition de P et Q}
\eeqa

Clearly the solution of the RHP for $\Phi$ exists as it can be built out of $\chi$. Its uniqueness can be seen along the same lines as in proposition 3.1.
Note that, because of the different possible analytic continuations of $\sg$ to the upper/lower half-planes
\eqref{continuation sg sur l'axe app+1}, $P$ and $Q$ will have different analytic continuations to the complex
plane depending on the value of $\Re\pa{z}$. In particular,
\beq
\ba{c c c c}
P\pa{z}&=&\a_{\ua}^{-2}\pa{z} - \f{\ex{-2i\pi \ov{\nu}_p \Xi\pa{\Re\pa{a_p-z}} } }{K_p\pa{z} \ex{ i x a_p} }
\pac{ x\pa{z-a_p} }^{2\de_p}   & z\in D_{a_p,\eps} \cap \mc{H}_+ \; , \vspace{3mm}\\
Q\pa{z}&=& K_p\pa{z} \f{\ex{2i\pi \nu_p \Xi\pa{\Re\pa{a_p-z}} +i x a_p } }{\pac{ x\pa{z-a_p} }^{2\de_p}  } - \a_{\da}^2\pa{z}
& z\in D_{a_p,\eps} \cap \mc{H}_- \; . \ea
\label{P et Q continues}\enq
\subsection{The second step $\Phi \tend \Upsilon $}
\label{SectionRHPsousSectionRHPPhiversUpsilon}
We now perform a transformation on $\Phi$. The resulting matrix $\Upsilon$ will have its jump matrices
exponentially close to the identity matrix, except in the vicinities of the singularities of $\sg$.
The jump matrix $G_{\Phi}$ can be  factorized into a product of an upper by a lower triangular matrix:
\beq
G_{\Phi}=M_{\ua}M_{\da} \;\;\; .
\label{factorization of GPhi}
\enq
\noindent The matrices $M_{\ua}$ (resp.  $M_{\da}$)
\beqa
 M_{\ua}\pa{z}&=& \left(\ba{cc}
                         1& P\pa{z}e^{ixz}\\
                         0& 1 \ea\right), \label{matrixM_+}\\
 M_{\da}\pa{z}&=& \left(\ba{cc}
                        1& 0\\
                         Q\pa{z}e^{-ixz}& 1 \ea \right),
\label{matrix M_-}
\eeqa
admit analytic continuations from the intervals $\intoo{-\infty}{a_1}$,   $\intoo{a_1}{a_{2}}$, \dots, $\intoo{a_n}{+\infty}$ to some interval depending domains in $U\cap \mc{H}_+$ (resp.$U\cap \mc{H}_+$). These analytic continuations are different if one starts from different intervals. In the following, so as to avoid any confusion, for $\Im z >0$, $M_{\ua}\pa{z}$ should be understood as the analytic continuation of $M_{\ua}\pa{\Re\pa{z}}$ from the interval containing $\Re\pa{z}$.
A similar statement holds for $M_{\da}\pa{z}$.

\begin{figure}[h]
\begin{center}

\includegraphics{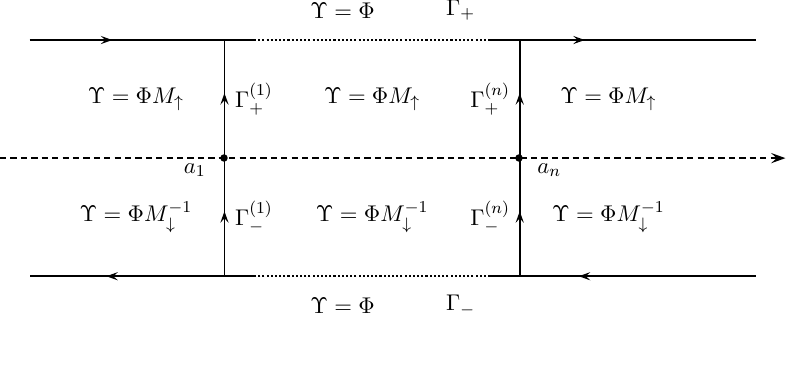}

\caption{Matrix $\Upsilon$ and its associated contour $\Ga=\Ga_+\cup\Ga_-\bigcup\limits_{k=1}^{n}\Ga_+^{\pa{k}}\cup \Ga_-^{\pa{k}}$.\label{contour pour le RHP de Y}}
\end{center}
\end{figure}

\noindent We draw a new contour
$\Ga=\Ga_+\cup\Ga_-\bigcup_{k=1}^{n} \Ga_+^{\pa{k}}\cup \Ga_-^{\pa{k}}$ in U. It allows to defines
a piecewise holomorphic matrix function $\Upsilon\pa{z}$  according to
Fig.\ref{contour pour le RHP de Y}.

As readily checked, $\Upsilon\pa{z}$ is continuous across $\R \setminus \bigcup_{i=1}^{n} \paa{a_i}$
and hence holomorphic in the interior of this new contour. By construction, $\Upsilon$ has cuts on the exterior
contour $\Ga_+\cup \Ga_-$. The additional cuts along  $\bigcup_{k=1}^{n}\Gamma^{\pa{k}}_{\pm} $
are due to the different analytic continuations for $P$ and $Q$ (and hence $M_{\ua}$ and $M_{\da}$)
to the strips $\paa{z: \; a_k<\Re\pa{z}<a_{k+1}}_{k=1}^{n-1}$ \textit{cf} \eqref{P et Q continues}.
\noindent The matrix $\Upsilon$ solves the following RHP:
\begin{itemize}
\item $\Upsilon$ is analytic in $\Cx\setminus \Ga \; ;$
\item  $\forall  k \in \intn{1}{n} \;,$ there exists $M_k\in \e{GL}_2\pa{\Cx}$ such that:
\vspace{3mm}
\newline
$\Upsilon\pa{z}  = M_k \paa{I_2+g\pa{z} B\pa{z}
+\abs{z-a_k} \pa{g\pa{z}+1}  \e{O}\pa{\ba{cc}1 & 1 \\ 1 & 1 \ea } }
M\pa{z} \;\, \quad z \tend a_k$
\item $\Upsilon \underset{z \tend \infty}{\tend} I_2 \; ;$
\item $\left\{ \ba{c} \Upsilon_{+}\pa{z} M_{\ua}\pa{z}=\Upsilon_-\pa{z}  ; \quad z \in \Gamma_+  \\
                      \Upsilon_{+}\pa{z} M_{\da}^{-1}\pa{z}=\Upsilon_-\pa{z}  ; \quad z \in \Gamma_-
                \ea \right. \; ;$
\item $\left\{ \ba{c} \Upsilon_{+}\pa{z} N^{\pa{l}}\pa{z}=\Upsilon_-\pa{z}  ; \quad z \in \Gamma^{(l)}_+  \\
                      \Upsilon_{+}\pa{z} \ov{N}^{\pa{l}}\pa{z}=\Upsilon_-\pa{z}  ; \quad z \in \Gamma^{(l)}_-
                \ea \right.  \; ,\qquad l \in \intn{1}{n} \; .$

\end{itemize}
There, the rank one matrix $B\pa{z}$ \eqref{definition matrice Bk} and the function $g$ \eqref{definition fonction g} are as in the RHP for $\chi$. We remind that $g$
has a singular behavior at $a_k$ given by \eqref{comportement local de g}.
The matrix $M$ is expressed in terms of $\a$ and $M_{\ua/\da}$ according to:
\beq
M\pa{z}= \left\{ \ba{cc}  \a^{\sg_3}_{\ua} M_{\ua}\pa{z}  & z \in \pa{ \mc{H}_+\setminus \cup_{k=1}^{n}\Ga_+^{\pa{k}} } \cap U
\vspace{2mm}\\
\a^{\sg_3}_{\da} M_{\da}\pa{z}  & z \in \pa{ \mc{H}_-\setminus \cup_{k=1}^{n}\Ga_-^{\pa{k}} } \cap U
\ea\right.
\enq
It is readily checked that the matrix $M\pa{z}$ has a singular behavior at $z\tend a_k$ given by
\beq
M\pa{z+a_k}=\left\{ \ba{cc c}
\e{O}\pa{ \ba{cc} \abs{z}^{\Re\pa{\nu_k}} & \abs{z}^{\e{min}\pa{-\Re\pa{\nu_k},\Re\pa{\ov{\nu_k}}}} \\ 0 & \abs{z}^{-\Re\pa{\nu_k}} \ea  }
                            &  \e{for} & z\tend 0 \; ,\;  z\in \mc{H_+} \vspace{2mm}\\
\e{O}\pa{ \ba{cc} \abs{z}^{-\Re\pa{\ov{\nu_k}}} & 0 \\ \abs{z}^{\e{min}\pa{-\Re\pa{\ov{\nu_k}},\Re\pa{\nu_k}}}
            & \abs{z}^{\Re\pa{\ov{\nu_k}}} \ea  }
                            &  \e{for} & z\tend 0 \; ,\;  z\in \mc{H_-}
\ea \right.  \; ;
\enq

\noindent Finally, the jump matrices $N^{\pa{l}}\pa{z}$, $\ov{N}^{\pa{l}}\pa{z}$ are defined by
\beqa
N^{\pa{l}}\pa{z}&=& \left( \ba{cc}
                                1 & n_l\pa{z}\ex{ixz} \\
                                0 & 1 \ea \right)
= \lim_{\substack{\eps \tend 0 \\ \Re\pa{\eps}>0 } } M_{+}^{-1}\pa{z-\eps}M_{+}\pa{z+\eps} \qquad z\in \Gamma_+^{\pa{l}} \; ,  \nonumber\\
\ov{N}^{\pa{l}}\pa{z}&=& \left( \ba{cc}
                                1 & 0 \\
                                \ov{n}_l\pa{z}\ex{-ixz} & 1 \ea \right)
= \lim_{\substack{\eps \tend 0 \\ \Re\pa{\eps}>0 } } M_{-}\pa{z-\eps}M_{-}^{-1}\pa{z+\eps} \qquad z\in \Gamma_-^{\pa{l}}  \;\;  \nonumber ;
\eeqa
\noindent and their entries read
\beqa
n_l\pa{z}&=& \f{ \pac{x\pa{z-a_l}}^{2\de_l} }{K_l\pa{z} \ex{ ix a_l}} \pa{\ex{-2i\pi \ov{\nu}_l}-1} \; ,\\
\ov{n}_l\pa{z}&=& \f{\ex{ix a_l} K_l\pa{z}}{\pac{x\pa{z-a_l}}^{2\de_l} } \pa{\ex{2i\pi\nu_l}-1} \;\; .
\eeqa

The solution of the RHP for $\Upsilon$ clearly exists and its uniqueness follows from a similar reasoning to proposition 3.1.
 The matrices $\Upsilon$ and $\chi$ are thus in a one-to-one correspondence.

Note that, apart from vicinities of the points $a_i$, $i \in \intn{1}{n}$,
the jump matrices for $\Upsilon$ are exponentially close to the identity. We have almost been
able to recast the original RHP into one suited for the Deift-Zhou steepest descent \cite{DeiftZhouSteepestDescentForOscillatoryRHP}; it only remains to build the parametrices around the $a_i$'s.

\section{Construction of the Parametrices, last transformation}
\label{Section Parametrix}
We first construct the parametrix for the model RHP on a small disc $D_{0,\eps}$ of radius $\eps>0$ and centered at 0.
This model parametrix will be the key ingredient of the parametrices around  the $a_i$'s.

\subsection{The model parametrix.}
\label{Section Parametrix sousSectModelPara}

\begin{figure}[h]
\begin{center}

\includegraphics{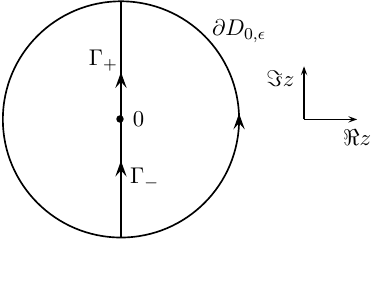}

\caption{Set of contours in the RHP for the model parametrix.\label{contours for the RHP for P}}
\end{center}
\end{figure}

\noindent The model parametrix $P$ is a solution to the following RHP:
\begin{itemize}
\item $P$ is analytic $D_{0,\eps}\setminus \paa{\Gamma_+\cup\Gamma_-} \; ;$
\item $P = I_2+ O\pa{\f{1}{\eps x}}\;\; , \;\;  z \in \partial D_{0,\eps} $ uniformly;
\item $\left\{ \ba{c} P_{+}\pa{z} N\pa{z}=P_-\pa{z}  \;\; , \quad z \in \Gamma_+ \cap D_{0,\eps} \\
                      P_{+}\pa{z} \ov{N}\pa{z}=P_-\pa{z} \;\; , \quad z \in \Gamma_- \cap D_{0,\eps}
                \ea \right. \; \; .$

\end{itemize}
The jump matrices of this model RHP read
\beqa
N\pa{z}&=& \left( \ba{cc}
                                1 & \f{  \pac{xz}^{2 \de }\ex{ixz}  }{K\pa{z}} \pa{\ex{-2i\pi\ov{\nu}}-1}  \\
                                0 & 1 \ea \right) \;\; ,\nonumber\\
\ov{N}\pa{z}&=& \left( \ba{cc}
                                1 & 0 \\
                                \f{K\pa{z}}{\pac{xz}^{2 \de}}\ex{-ixz}  \pa{\ex{2i\pi\nu}-1}  & 1 \ea \right)
\;\; ,  \nonumber
\eeqa
and the function $K\pa{z}$ is assumed to be holomorphic and non-vanishing on $D_{0,\eps}$. Note that the boundary
$\partial D_{0,\eps}$ of the disk $D_{0,\eps}$ is canonically oriented  just as
depicted in Fig. \ref{contours for the RHP for P}.
The RHP for P admits many solutions. For instance having one solution $P$, one can build another one by
multiplying $P$ on the left by a holomorphic matrix on $D_{0,\eps}$ that is
equal to $I_2$ up to corrections that are uniformly an  $O\pa{\tf{1}{\eps x}}$ on $\partial D_{0,\eps}$ and hence on the whole disk $D_{0,\eps}$.

A solution to the above RHP can be built thanks to the following procedure. We first assume that
$K$ is a constant. Then the Riemann-Hilbert problem for $P$ can be solved by the standard procedure. One
performs the transformation
\beq
P\pa{z} = \Theta\pa{\zeta} \pac{\ex{\f{i \zeta}{2}} \zeta^{\de} }^{-\sg_3} \qquad \e{with} \;\; \zeta=xz
\enq
\noindent so that the jump matrices for $\Theta$ are piecewise constant, and $\Theta$ is a solution of a RHP
on $D_{0,x\eps}$. This RHP is solved explicitly in the limit $x\eps \tend +\infty$ by the standard differential equation method \cite{ItsDifferentialMethodForParametrix}. It is then enough to go back to the original matrix $P$. Eventually,
we get that
\beq
P\!\pa{z} =\!\!  \pa{\!\! \ba{cc} \!\!
                        \Psi\pa{\ga-\de,1+2\ga;-ixz}
                                    &ib_{12} \Psi\pa{1+\ga+\de,1+2\ga;ixz} \\
                        -ib_{21} \Psi\pa{1+\ga-\de,1+2\ga;-ixz}
                                        & \Psi\pa{\ga+\de,1+2\ga;ixz} \ea \!\! }
\f{ L}{ \pa{xz}^{\de\sg_3-\ga}} \, .
\enq
Here $\Psi\pa{a,c;z}$ is Tricomi's confluent hypergeometric function (CHF) given in \eqref{Appendix Tricomi CHf definition}. We remind that it has a cut on $\R_-$. The piecewise constant matrix $L$ depends on $\de$ and $\ga$,
\beq
L=\left\{\ba{l c} \ex{i\f{\pi \de}{2}} \ex{-i\f{\pi \ga}{2}\sg_3}
                & -\tf{\pi}{2}<\e{arg}\pa{z}< \tf{\pi}{2} \vspace{4mm}\\
                 \ex{i\f{\pi \pa{\de-\ga}}{2}}\left( \ba{cc} 1  &0 \\
                                                0& \ex{-2i\pi \de-i\pi\ga}
                                                   \ea \right)
                & \tf{\pi}{2}<\e{arg}\pa{z}< \pi \vspace{4mm}\\
                 \ex{i\f{\pi \pa{\de+\ga}}{2}}\left( \ba{cc} \ex{-2i\pi \de+i\pi \ga}  &0 \\
                                     0& 1
                                    \ea \right)  & -\pi<\e{arg}\pa{z}< -\tf{\pi}{2}
        \ea \right. \;\; ,
\enq
\noindent whereas the coefficients $b_{12}$ and $b_{21}$ also exhibit an additional dependence on $K$
\beq
b_{12}= \f{i \ex{-i\pi\ga}}{K} \f{\Gamma\pa{1-\ga+\de}}{\Gamma\pa{-\ga-\de}}  \hspace{1cm},\hspace{1cm}
b_{21}=  -i K \ex{i\pi\ga} \f{\Gamma\pa{1-\ga-\de}}{\Gamma\pa{\de-\ga}} \; .
\enq
Lastly, one has $\ddet{}{P}=1$.

\vspace{2mm}
Using the asymptotic behavior of Tricomi's CHF \eqref{asy-Psi}, one readily checks that
$P$ has indeed the correct asymptotic behavior. The jump conditions can be checked thanks to the monodromy
properties of Tricomi's CHF. These are \eqref{cut-Psi-1} for the jump condition on $\Ga_+$ and $\eqref{cut-Psi-2}$ in the case of the jump condition on $\Ga_-$.
Moreover, P has  no jump across $D_{0,\eps}\setminus \intoo{-\eps}{0}\, :$ the discontinuity in $L$ is there to compensate
the one of $\pa{z}^{\ga-\de\sg_3}$. Hence $P$ is holomorphic on $\intff{-\eps}{0}$.
The fact that $\ddet{}{P}=1$ can be seen as follows. We first assume that $\abs{\Re\pa{\ga}}<\tf{1}{4}$. Then, using the local behavior at $z=0$ of Tricomi's CHF \eqref{Psi s'ecrit comme Phi}, we get that $\ddet{}{P}=\e{o}\pa{z^{-1}}$.  On the other hand, writing the jump conditions for $\ddet{}{P}$, one can easily convince oneself that the latter function is holomorphic on $\Cx\setminus \paa{0}$. Its singularity at $z=0$ is thus of a removable type. As $\ddet{}{P}\tend 1$ when $z\tend \infty$, we get that necessarily $\ddet{}{P}=1$. To reach the case of generic parameters $\ga$, we fix $z \not=0$ and invoque the fact that $\Psi\pa{a,c;z}$ is a holomorphic function of $a$ and $c$. It follows that $\ddet{}{P}\pa{z}$ is holomorphic in $\ga$.
As it is constant in the region $\abs{\Re\pa{\ga}}<\tf{1}{4}$, we get that it is constant on $\Cx$. Hence, $\ddet{}{P}\pa{z}=1$ for all $z\not=0$. Now we get that $\ddet{}{P}$ is bounded in every punctured neighborhood of $z=0$. It thus follows that it cannot have any power-law singularity at $z=0$, and $\ddet{}{P}=1$.

\vspace{3mm}

In this way we have built a solution of the RHP for $P$ in the case of constant functions $K$.
In order to extend this solution to functions $K$ that are holomorphic and non-vanishing in some open neighborhood of $\ov{D}_{0,\eps}$, it is
enough to notice that replacing the constant K appearing in the formulae above by a holomorphic non-vanishing function $K\pa{z}$ on $\ov{D}_{0,\eps}$ doesn't change the analyticity of P on $D_{0,\eps}\setminus \pa{\Ga_+\cup\Ga_-}$,
nor its asymptotic behavior on the boundary $\Dp{}D_{0,\eps}$ . As the jump conditions
hold pointwise, these are also satisfied.  Thence we get a solution to the general RHP for the parametrix.

\subsection{The parametrix around $a_k$}
\label{Section Parametrix sousSectLocalPara}

Let $P_{a_p}$ be defined as

\beq
P_{a_p}\!\!\pa{  z}\! = \!\! \pa{\!\! \ba{cc}\!\!
                        \Psi\pa{\ga_p-\de_p;-i\zeta_p}
                                    &ib^{\pa{p}}_{12}\pa{z} \Psi\pa{1+\ga_p+\de_p;i\zeta_p} \\
                        -ib^{\pa{p}}_{21}\pa{z} \Psi\pa{1+\ga_p-\de_p;-i\zeta_p}
                                        & \Psi\pa{\ga_p+\de_p;i\zeta_p} \ea \!\!}\!
\f{L_p}{\pa{\zeta_p}^{\de_p\sg_3-\ga_p} }\, .
\label{parametrix around ap}
\enq
There we have set $\zeta_p=x\pa{z-a_p}$ and the second argument of the CHF's is implicitly assumed to be
$1+2\ga_p$. The piecewise constant matrix $L_p$ reads
\beq
L_p=\left\{\ba{cc}  \ex{i\f{\pi \de_p}{2}} \ex{-i\f{\pi \ga_p}{2}\sg_3}
                & -\tf{\pi}{2}<\e{arg}\pa{z-a_p}< \tf{\pi}{2} \vspace{4mm}\\
                  \ex{i\f{\pi \pa{\de_p-\ga_p}}{2}}\left( \ba{cc} 1  &0 \\
                                                0& \ex{-2i\pi \de_p-i\pi\ga_p}
                                                   \ea \right)
                & \tf{\pi}{2}<\e{arg}\pa{z-a_p}< \pi \vspace{4mm}\\
                 \ex{i\f{\pi \pa{\de_p+\ga_p}}{2}}\left( \ba{cc} \ex{-2i\pi \de_p+i\pi \ga_p}  &0 \\
                                     0& 1
                                    \ea \right)  & -\pi<\e{arg}\pa{z-a_p}< -\tf{\pi}{2}
        \ea \right. \;\; ,
\enq
\noindent and the coefficients $b^{\pa{p}}_{12}\pa{z}$ and $b^{\pa{p}}_{21}\pa{z}$ are
\beq
b^{\pa{p}}_{12}\pa{z}= \f{i \ex{-i\pi\ga_p}}{K_p\pa{z}} \f{\Gamma\pa{1-\ga_p+\de_p}}{\Gamma\pa{-\ga_p-\de_p}}  \hspace{5mm} , \hspace{1.2cm}
b^{\pa{p}}_{21}\pa{z}=  -i K_p\pa{z} \ex{i\pi\ga_p} \f{\Gamma\pa{1-\ga_p-\de_p}}{\Gamma\pa{\de_p-\ga_p}} \;.
\enq

\begin{prop}
The matrix $P_{a_p}\pa{z}$ plays the role of
a parametrix around $a_p$ in the sense that:
\begin{itemize}
\item $\Upsilon\, P_{a_p}^{-1}$ is holomorphic inside of $D_{a_p,\eps}$ \; ,
\item $P_{a_p}^{-1}=I_2+\e{O}\pa{\tf{1}{x^{1-2\Re\pa{\de_p}}} }$  uniformly on
$\Dp{}D_{a_p,\eps}$.
\end{itemize}
\end{prop}

\Proof
It follows from $\ddet{}{P_{a_p}}=1$ that $P_{a_p}$ is invertible and that \newline $P_{a_p}^{-1}=  ^{t} \e{Comat}\pa{P_{a_p}}$. The fact that $ P^{-1}_{a_p}= I_2+\e{O}\pa{\tf{1}{x^{1-2\Re\pa{\de_p}}} }$ uniformly  on $\Dp{}D_{a_p,\eps}$  follows  from  the asymptotic behavior
on the boundary of $\Dp{}D_{0,\eps}$ of the matrix $P$ defined in the previous section together with the fact that the function  $K_p\pa{z}$ appearing in the definition of $b_{12}^{\pa{p}}\pa{z}$
and $b_{21}^{\pa{p}}\pa{z}$ depends on $x^{2\de_p}$,  \textit{cf} \eqref{definition de K_p}.

As $P_{a_p}$ and $\Upsilon$ have the same jump matrices on $\paa{\Ga_+^{\pa{p}}\cup \Ga_-^{\pa{p}}}\cap D_{a_p,\eps}$,
$\Upsilon P^{-1}_{a_p}$ can be analytically continued to the punctured disk $D_{a_p,\eps} \setminus \paa{a_p}$.
It remains to see that the singularity at $z=a_p$ is removable.

The local behavior of $\Upsilon$ at $z\tend a_p$ as well as the one of CHF at the origin (\textit{cf} appendix \eqref{Psi s'ecrit comme Phi})  imply that  $\Upsilon P^{-1}_{a_p}$  has at most a power-law singularity at $a_p$. In particular, it cannot have an essential singularity at $a_p$. The singularity can only be a pole of some finite order. To set aside this possibility, it is thus enough to check that $\Upsilon P^{-1}_{a_p}$ is bounded in the quadrant $0\leq \e{arg}\pa{z-a_p}\leq \tf{\pi}{2}$.

The addition formulae for the CHF \eqref{Phi s'ecrit comme Psi} allow to express the product $P_{a_p}M_{\ua}^{-1} \a^{-\sg_3}$ as
\bem
P_{a_p}\pa{z}M_{\ua}^{-1}\pa{z} \pac{\a_{\ua}\pa{z}}^{-\sg_3}= \f{\pa{-i\zeta_p}^{\ga_p-\de_p}}{\a_{\ua}\pa{z}}
%
%
%\pa{ \ba{c}  \Psi\pa{\ga_p-\de_p; -i\zeta_p}  & - \Psi\pa{\ga_p-\de_p; -i\zeta_p}  \\
%
%
%-ib_{21}^{\pa{p}}\pa{z}  \Psi\pa{1+\ga_p-\de_p; -i\zeta_p} &  i b_{21}^{\pa{p}}\pa{z}\Psi\pa{1+\ga_p-\de_p; -i\zeta_p}  \ea} \\
\pa{ \ba{c}  \Psi\pa{\ga_p-\de_p; -i\zeta_p}    \\
-ib_{21}^{\pa{p}}\pa{z}  \Psi\pa{1+\ga_p-\de_p; -i\zeta_p}   \ea}  \cdot \pa{1 \;  -\ex{ix z}} \\
+ \f{\ex{i\f{\pi}{2}\pa{\ga_p-\de_p}}}{\Ga\pa{1-2\ga_p}} \f{\a_{\ua}\pa{z}}{\zeta_p^{\ga_p-\de_p}}
\pa{\ba{cc}  0 & \Ga\pa{1+\de_p-\ga_p}\Phi\pa{1-\ga_p+\de_p, 1-2\ga_p;i\zeta_p} \\
   0 & \Ga\pa{1-\de_p-\ga_p} \ex{2i\pi \ga_p} \Phi\pa{\de_p-\ga_p, 1-2\ga_p;i\zeta_p}   \ea}
\label{formule explicite produit CHF parametrix et autre}
\end{multline}
Here, once again, the second argument of Tricomi's CHF is assumed to be $1+2\ga_p$.

Recall that the $\Phi$ functions are regular at $\zeta_p=0$ whereas the $\Psi$ functions have a power-law singularity of the type $\e{O}\pa{\zeta_p^{-2\ga_p}}$. Since $\tf{\a_{\ua}\pa{z}}{\zeta_p^{\ga_p-\de_p}=\e{O}\pa{1}}$ for $z\tend a_p$, we get that
\beq
P_{a_p}M_{\ua}^{-1} \pac{\a_{\ua}}^{-\sg_3}= \paa{\e{O}\pa{1} +\e{O}\pa{\abs{z-a_p}^{-2\ga_p}}+\e{O}\pa{\log\abs{z-a_p}} }
\e{O}\pa{\ba{cc} 1 & 1 \\ 1 & 1 \ea} \;
\enq
Moreover, using the representation \eqref{formule explicite produit CHF parametrix et autre}, one obtains that
\beq
P_{a_p}\pa{z}M_{\ua}^{-1}\pa{z} \a^{-\sg_3}\pa{z}  ^t \e{Comat}\pa{B\pa{z}}  = \e{O}\pa{\ba{cc} 1 & 1 \\ 1 & 1 \ea } \; .
\enq
Hence, using the local estimates for $\Upsilon$ around $a_p$, we get that
\beq
P_{a_p}\pa{z} \Upsilon^{-1}\pa{z} =\paa{ \e{O}\pa{1} +\e{O}\pa{\abs{z-a_p}^{-2\ga_p}}+\e{O}\pa{\log\abs{z-a_p}} }
\e{O}\pa{\ba{cc} 1 & 1\\ 1 & 1 \ea  } \;.
\enq
Here we have used that $\ddet{}{\Upsilon}=1$ what implies that, for $z\tend a_p$,
\beq
\Upsilon^{-1}\pa{z}=M_{\ua}^{-1} \a^{-\sg_3}_{\ua} \paa{I_2 + g\pa{z} ^t\e{Comat}\pa{B\pa{z}}
+ \abs{z-a_p} \pa{g\pa{z}+1} \e{O}\pa{\ba{cc} 1&1 \\1&1 \ea  }} M_{p}^{-1} \; ,
\enq
for some $M_p \in GL\pa{2,\Cx}$.
Therefore, as $\Re\pa{2\ga_p}<1$, we see that $ P_{a_p}\pa{z}\Upsilon^{-1}\pa{z}$ cannot have a pole at $z=a_p$. The singularity at $z=a_p$ is hence removable and $ P_{a_p}\pa{z}\Upsilon^{-1}\pa{z}$ is analytic on $D_{a_p,\eps}$.
As $\ddet{}{ P_{a_p}\pa{z}\Upsilon^{-1}\pa{z}}=1$, we have that
\beq
\Upsilon\pa{z}P_{a_p}^{-1}\pa{z}= ^t\e{Comat}\pa{ P_{a_p}\pa{z}\Upsilon^{-1}\pa{z}}%
\enq
is also analytic on $D_{a_p,\eps}$. \qed

\subsection{The last transformation $\Upsilon\tend \Omega$}
\label{SectionEstResRHPOmega}
\begin{figure}[h]
\begin{center}

\includegraphics{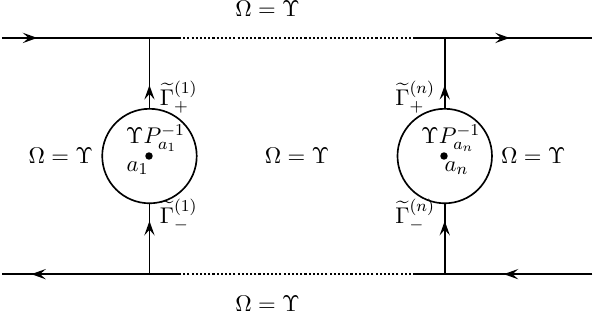}

\caption{Matrix $\Omega$ and new contours $\Sg_{\Omega}$.\label{contour pour le RHP de R}}
\end{center}
\end{figure}

The matrix
\beq
\Omega=\left\{  \ba{cc}
                    \Upsilon P_{a_k}^{-1} & z\in D_{a_k,\eps}\\
                    \Upsilon & z \in \mathbb{C}\setminus\bigcup\limits_{p=1}^{n}\ov{D}_{a_k,\eps}
            \ea \right.
\enq
satisfies the RHP
\begin{itemize}
\item $\Omega$ is analytic in $\mathbb{C}\setminus \Sg_{\Omega} \; ;$
\item $\Omega = I_2+ O\pa{\tf{1}{z}} \quad,\quad   z \tend \infty \; ;$
\item $ \left\{ \ba{rc c l c c c}  \Omega_{+}\pa{z} M_{\ua}\pa{z}&=\Omega_-\pa{z}  &,& \; z \in \Gamma_+  \; ;\\
                              \Omega_{+}\pa{z} M_{\da}^{-1}\pa{z}&=\Omega_-\pa{z}  &,& \;z \in \Gamma_- \; ; \\
                      \Omega_{+}\pa{z} N^{\pa{l}}\pa{z}&=\Omega_-\pa{z}  &, & \; z \in \widetilde{\Gamma}^{\pa{l}}_+ \; ; \\
                      \Omega_{+}\pa{z} \ov{N}^{\pa{l}}\pa{z}&=\Omega_-\pa{z}  &,& \; z \in \widetilde{\Gamma}^{\pa{l}}_- \; ; \\
                      \Omega_{+}\pa{z} P_{a_l}\pa{z}&=\Omega_-\pa{z}  &, &\; z \in \partial D_{a_l,\eps}  \;. \ea \right.$
\end{itemize}
The solution of this RHP for $\Omega$ exits and is unique as seen by already invoked arguments. The newly introduced  contours are all depicted in Fig.\ref{contour pour le RHP de R}.

The jump matrix $v_{\Omega}$ for $\Omega$ is uniformly $I_2+O\pa{\tf{1}{x^{1-\rho}}}$
in the $L^{2}\pa{\Sg_{\Omega}}$ and $L^{\infty}\pa{\Sg_{\Omega}}$ sense, \textit{ie} there exists a constant $c$
such that
\beq
\norm{v_{\Omega}-I_2}_{L^{2}\pa{\Sg_{\Omega}}}+\norm{v_{\Omega}-I_2}_{L^{\infty}\pa{\Sg_{\Omega}}}\leq c x^{\rho-1}
\; .
\enq
Here, we choose the following $L^{2}\pa{\R}$ and $L^{\infty}\pa{\R}$ matrix norms
\beq
\norm{A\pa{s}}_{L^{2}\pa{\Sg_{\Omega}}}^2= \Int{\Sg_{\Omega}}{} \e{tr}\pac{A^{\dagger}\pa{s} A\pa{s}} \abs{\dd s}
\; , \quad \norm{A\pa{s}}_{L^{\infty}\pa{\Sg_{\Omega}}}=\max_{i,j}  \norm{A_{ij}\pa{s}}_{L^{\infty}\pa{\Sg_{\Omega}}} \; .
\enq
We also remind that  $\rho=2\underset{k}{\max} \abs{\Re \pa{\de_k}}<1$. By using the matrix integral equation equivalent
to the RHP for $\Omega$ we see that $\Omega \limit{x}{+\infty} I_2$ uniformly. Moreover the first corrections to $\Omega$ are uniformly $\e{O}\pa{\tf{1}{x^{1-\rho}}}$. As a consequence, $I_2$ is the unique solution of the RHP for $\Omega$,
up to the uniformly $\e{O}\pa{\tf{1}{x^{1-\rho}}}$ corrections.

\section{Asymptotics of the resolvent}
\label{Section Estimates for the Resolvent}

Recall that the resolvent $I-R$ of $I+V$ can be expressed in terms of $f_+$ and $f_-$,
\textit{cf} \eqref{Resolvent en termes de f_+ et f_-}. Having solved the RHP for $\chi$ perturbatively in $x$,
we use this asymptotic solution in order to obtain the asymptotics to the leading order of the functions $f_{\pm}$.
The latter yield the leading asymptotic behavior of the resolvent.

\subsection{The zeroth order approximants to $f_{\pm}$.}
\begin{definition}
Let $\chi^{bk}$ and $\chi_p^{loc} $ be the matrices
\beqa
\chi^{bk}\pa{z}&=& M_{\ua}^{-1}\pa{z} \pac{\a_{\ua}\pa{z}}^{-\sg_3}  \;\;,\;   \nonumber\\
\chi_p^{loc}\pa{z}&=& P_{a_p}\pa{z}M_{\ua}^{-1}\pa{z} \pac{\a_{\ua}\pa{z}}^{-\sg_3} \;\;. \nonumber
\eeqa
\end{definition}

Let $\Omega_{\eps}$ be the solution of the RHP defined in subsection \ref{SectionEstResRHPOmega} and where all the circles in the contour $\Sg_{\Omega_{\eps}}$,
as depicted on fig.\ref{contour pour le RHP de R}, have a radius $\eps$. Then, if $\chi$ is the solution of the RHP defined in subsection \ref{SectionRHPsousSectionRHPchi}, we have
\beqa
\Omega^{-1}_{\f{\eps}{2}}\pa{z}\chi\pa{z}&=&\chi^{bk}\pa{z} \;\;,\;   z \in U\! \cap \paa{  z \in \mc{H}_+  : \Re\pa{z} \in \intoo{ a_p+\tf{\eps}{2} }{ a_{p+1}-\tf{\eps}{2} } } \;.
\label{definition ecart Chi bulk} \\
\Omega^{-1}_{2\eps}\pa{z}\chi\pa{z}&=&\chi^{loc}_p\pa{z} \;\; , \; z \in D_{a_p, 2\eps} \; .\label{definition chi local}
\eeqa
In a sense that will become clear in the following, $\chi^{bk}$ is the leading  solution of the RHP for $\chi$ when $z$ is
uniformly away from the singularities at the $a_p$'s and $\chi_p^{loc}$ is the leading  solution of the RHP for $\chi$ when $z$ belongs to the disk $D_{a_p,2 \eps}$.

The advantage of using the solution of Riemann-Hilbert problems with two sizes of disks around the $a_k$'s , $\Omega_{2\eps}$ and $\Omega_{\tf{\eps}{2}}$ is that the solution $\chi_+\pa{z}$ defined by \eqref{definition ecart Chi bulk} on $\R\setminus \cup_{k=1}^{n}\intff{a_k-\eps}{a_k+\eps}+i0^+$ and by \eqref{definition chi local}
 on $\cup_{k=1}^{n}\intff{a_k-\eps}{a_k+\eps}+i0^+$ has a smooth correction matrix $\Omega$ around the gluing point  $z=\a_p \pm \eps$. This will simplify our forthcoming analysis when integrating the solution $\chi$ around the points  $a_k=\pm\eps$.
If we would have used a single solution $\Omega_{\eps}$, then we should have had derived additional estimates for the behavior of this matrix around $a_k=\pm \eps$, as a priori it could exhibit a non-smooth behavior there.
Thence, we circumvent additional complications.

\begin{prop}
\label{proposition f plusmoins loc and bk}
Let
\beq
\pa{\ba{c} f_{+}^{bk}\pa{z} \\ f_-^{bk}\pa{z} \ea } \equiv \chi^{bk}\pa{z}\pa{\ba{c} e_{+}\pa{z} \\ e_-\pa{z} \ea } \;\; ,
\;\;
 \pa{\ba{c} f_{+;p}^{loc}\pa{z} \\ f_{-;p}^{loc}\pa{z} \ea } \equiv \chi_p^{loc}\pa{z}\pa{\ba{c} e_{+}\pa{z} \\ e_-\pa{z} \ea } \;.
\enq
Then
\beq
\left(\ba{c} f^{bk}_{+}\pa{z} \\
             f^{bk}_{-}\pa{z} \ea  \right) =\pa{\a_{\ua} e_-}^{-\sg_3}  \left(\ba{c} \sg_p^{-1}\pa{z} \\
             1 \ea  \right)  \;\;\; \e{for} \quad  z \in \paa{a_p< \Re\pa{z}<a_{p+1}} \cap U \; ,
\label{definition f plus-moins bulk}
\enq
with $\sg_p$ having already been defined in \eqref{sigma continuation sur ap-1 ap}. Also
\bem
\pa{ \ba{c} f^{loc}_{+;p}\pa{z} \\
             f^{loc}_{-;p}\pa{z} \ea  } =
  \f{\ex{i\f{\pi}{2}\pa{\ga_p\sg_3-\de_p }} x^{-\de_p\sg_3} }{\Gamma\pa{1-2\ga_p }x^{\ga_p}}
\pa{  \ba{c c}  e_{+}\pa{z}  \pac{\wh{\a}_{\da}^{\pa{p}}\pa{z} }^{-1} & 0 \\
                0 & \wh{\a}_{\ua}^{\pa{p}}\pa{z} e_-\pa{z} \ea  }
 \\
\qquad \times
\pa{ \ba{c} \Gamma\pa{1+\de_p-\ga_p } \Phi\pa{-\ga_p-\de_p,1-2\ga_p;-i\zeta_p}   \vspace{2mm}\\
          \Gamma\pa{1-\de_p-\ga_p} \Phi\pa{\de_p-\ga_p,1-2\ga_p;i\zeta_p }  \ea  } \;\;\; \e{for}\quad z \in D_{a_p;2\eps} \; .
\nonumber
\end{multline}
The functions $\wh{\a}_{\pm}^{\pa{p}}$ were defined in \eqref{alpha plusmoins regularise} whereas $\Phi$ is Humbert's CHF
\eqref{definition Humbert's CHF} and $\zeta_p=x\pa{z-a_p}$.
\end{prop}

\Proof

The proof of first formula \eqref{definition f plus-moins bulk} is straightforward. We explain how to derive an expression for $f^{loc}_{\pm;p}$ in a vicinity of $a_p$. In the intermediary calculations we suppose that $z$ belongs to the quadrant $\paa{z: a_p< \Re\pa{z}< a_{p} + \eps} \cap \mc{H}_+$.  The final result is however valid on the whole disk as can be seen through a direct computation carried out on the other quadrants. In order to lighten the notations, the second argument of Tricomi's CHF is undercurrented to be $1+2\ga_p$ and we have set $\zeta=x\pa{z-a_p}$.
\beqa
\left(\ba{c} f^{loc}_{+}\pa{z} \\
             f^{loc}_{-}\pa{z} \ea  \right) &=&
         \left( \ba{cc}
          \Psi\pa{\ga_p-\de_p;-i\zeta}   &ib^{\pa{p}}_{12}\pa{z} \Psi\pa{1+\ga_p+\de_p;i\zeta} \\
                        -ib^{\pa{p}}_{21}\pa{z} \Psi\pa{1+\ga_p-\de_p;-i\zeta}  & \Psi\pa{\ga_p+\de_p;i\zeta} \ea \right) \nonumber \\
&& \quad  \times \zeta^{\ga_p}   \ex{i\f{\pi}{2}\pa{\de_p-\ga_p \sg_3 } }
\pa{\zeta^{\de_p}\a_{\ua}\pa{z} e_-\pa{z}}^{-\sg_3}
\left(\ba{c} \sg_p^{-1}\pa{z} \\
             1 \ea  \right)
\nonumber\\
 &=&
\zeta^{\ga_p}   \ex{i\f{\pi}{2}\de_p}
\pa{\ex{i \f{\pi}{2} \ga_p}\zeta^{\de_p}\a_{\ua}\pa{z} e_-\pa{z}}^{-\sg_3}
                    \left( \ba{cc}  \sg_p^{-1}\pa{z} & 0   \\
                                     0              & 1 \ea \right) \\
&& \quad \times
         \left( \ba{cc}
          \Psi\pa{\ga_p-\de_p;-i\zeta}   &g^{\pa{p}}_{12} \Psi\pa{1+\ga_p+\de_p;i\zeta} \\
                        g^{\pa{p}}_{21} \Psi\pa{1+\ga_p-\de_p;-i\zeta}  & \Psi\pa{\ga_p+\de_p;i\zeta} \ea \right)
\left(\ba{c} 1 \\
             1 \ea  \right)\; .
\nonumber
\eeqa
\noindent With
\beq
g^{\pa{p}}_{12}= -\f{\Gamma\pa{1+\de_p-\ga_p}}{\Gamma\pa{-\de_p-\ga_p}}\ex{-i\zeta} \qquad
g^{\pa{p}}_{21}=-\f{\Gamma\pa{1-\ga_p-\de_p}}{\Gamma\pa{\de_p-\ga_p}}\ex{i\zeta}\; .
\enq
And finally using the recombination formulas for CHF \eqref{Phi s'ecrit comme Psi}, the claim follows.$\quad \Box$.

\begin{rem}
The function $f^{bk}_{\pm}$, resp. $f_{\pm;p}^{loc}$, are good approximates to $f_{\pm}$ in their
respective domains of validity. More precisely, for $z \in \R\setminus \cup_{p=1}^{n}\intff{a_p-\tf{\eps}{2}}{a_p+\tf{\eps}{2}}$
\beq
\pa{\ba{c} f_+\pa{z}-f_{+}^{bk}\pa{z} \\ f_-\pa{z}-f_-^{bk}\pa{z} \ea }= \pa{\chi\pa{z}-\chi^{bk}\pa{z}}
\pa{\ba{c} e_{+}\pa{z} \\ e_{-}\pa{z} \ea }=
 \pa{\Omega_{\f{\eps}{2}}\pa{z}-I_2}\pa{\ba{c} f_{+}^{bk}\pa{z} \\ f_{-}^{bk}\pa{z} \ea } \;\; ,
\enq
and for $z\in D_{a_p,\eps}$,
\beq
 \pa{\ba{c} f_+\pa{z}-f_{+;p}^{loc}\pa{z} \\ f_-\pa{z}-f_{-;p}^{loc}\pa{z} \ea } =
\pa{\chi\pa{z}-\chi^{loc}_p\pa{z}}\pa{\ba{c} e_{+}\pa{z} \\ e_-\pa{z} \ea }=
  \pa{\Omega_{2\eps}\pa{z}-I_2}\pa{\ba{c} f^{loc}_{+;p}\pa{z} \\ f^{loc}_{-;p}\pa{z} \ea } \;.
\enq

\end{rem}

\subsubsection{Uniform estimates for the resolvent}

\begin{definition}
Let $R_0\pa{\xi,\eta}$ be called the zeroth order resolvent. We define it in terms of $f_{\pm}^{bk}$ and
$f_{\pm;p}^{loc}$ as :
\beq
\ba{l}
R_{0}\pa{\xi,\eta}= \pa{\sg\pa{\xi}-1} \f{f_+^{bk}\pa{\xi} f_{-}^{bk}\pa{\eta}-f_{-}^{bk}\pa{\xi} f_{+}^{bk}\pa{\eta}  }{2i\pi\pa{\xi-\eta}}
\; ; \quad \!  \eta\; , \xi \in \R\setminus \bigcup\limits_{p=1}^{n} \intff{a_p\!-\!\eps}{a_p\!+\!\eps}  \\
R_{0}\pa{\xi,\eta}= \pa{\sg\pa{\xi}-1} \f{f_{+;p}^{loc}\pa{\xi} f_{-;p}^{loc}\pa{\eta}-f_{-;p}^{loc}\pa{\xi} f_{+;p}^{loc}\pa{\eta}  }{2i\pi\pa{\xi-\eta}}
 \; ; \quad \eta\; , \xi \in \intff{a_p-\eps}{a_p+\eps}   \; .\ea  \nonumber
\enq
\end{definition}

Using the explicit expression for $f_{\pm}^{bk}$ and $f_{\pm;p}^{loc}$ we get the local expressions for the diagonal zeroth order resolvent
\beq
2i\pi\sg\pa{\xi} \f{R_0\pa{\xi,\xi} }{\sg\pa{\xi}-1}
=ix-\Dp{\xi}\pa{\ln\a_{+} \a_{-}}\pa{\xi}  \; ,
\label{Resolvent sur interval non singulier}
\enq
where $\xi \in \R\setminus\bigcup\limits_{p=1}^{n}\intff{a_p-\eps}{a_p+\eps}$ and we remind that
$\a_+$ and $\a_-$ are the boundary values of $\a$ from the upper/lower half plane.
 Whereas for $\xi \in \intoo{a_p-\eps}{a_p+\eps}$
\beq
2i\pi\ex{i\pi \de_p} \wh{\sg}_p\pa{\xi}
\f{x^{2\ga_p}R_0\pa{\xi,\xi} }{\sg\pa{\xi}-1} = -ix \tau\pa{\ga_p,\de_p;\zeta_p }
-\Dp{\xi}\pa{\ln \wh{\a}^{\pa{p}}_{\ua,+} \wh{\a}^{\pa{p}}_{\da,-}}\pa{\xi} \varphi\pa{\ga_p,\de_p;\zeta_p} \; .
\label{resolvent autour des singlarites}
\enq
\noindent One should keep in mind that $\zeta_p=x\pa{\xi-a_p}$ and that $\wh{\a}^{\pa{p}}_{\ua,+}$
(resp. $\wh{\a}^{\pa{p}}_{\da,-}$) is the  boundary value of $\wh{\a}^{\pa{p}}_{\ua}$
(resp. $\wh{\a}^{\pa{p}}_{\da}$) from the upper (resp. lower) half-plane.
In the above formulae we have introduced two functions
\bem
\Ga\pab{1-2\ga,1-2\ga}{1+\de-\ga,1-\de-\ga}\tau\pa{\ga,\de;t}=
-\Phi\pa{-\ga-\de,1-2\ga;-it}\Phi\pa{\de-\ga,1-2\ga;it}  \\
+
\pa{\Dp{z}\Phi}\pa{-\ga-\de,1-2\ga;-it}\Phi\pa{\de-\ga,1-2\ga;it} \\
 + \Phi\pa{-\ga-\de,1-2\ga;-it}\pa{\Dp{z}\Phi}\pa{\de-\ga,1-2\ga;it} \; ,
\label{definition de tau}
\end{multline}
\noindent and
\beqa
\Ga\pab{1-2\ga,1-2\ga}{1+\de-\ga,1-\de-\ga}\varphi\pa{\ga,\de;t}=\Phi\pa{-\ga-\de,1-2\ga;-it}\Phi\pa{\de-\ga,1-2\ga;it}
 \;\; .
\label{definition de varphi}
\eeqa
\noindent Also, we have used the standard notation
\beq
\Ga\pab{a_1, \dots, a_n}{b_1, \dots, b_m}= \f{\pl{k=1}{n} \Ga\pa{a_k} }{ \pl{k=1}{m} \Ga\pa{b_k} } \;.
\enq

We now focus on the relationship between the exact resolvent $R$ and the zeroth order one $R_0$. Observe that one can write the exact resolvent $R$ as
\beq
R\pa{\xi,\eta}=\pa{-e_-\pa{\xi} \; e_+\pa{\xi}}
    \pac{\chi_{\eps'}^{\pa{0}}\pa{\xi}}^{-1}
    \f{\Omega_{\eps'}^{-1}\pa{\xi} \Omega_{\eps'}\pa{\eta} }{\xi-\eta}
    \chi_{\eps'}^{\pa{0}}\pa{\eta}
    \pa{\ba{c} e_+\pa{\eta} \\ e_-\pa{\eta} \ea  }  \;\; ,
\enq
Where $\chi_{\eps'}^{\pa{0}}\pa{\xi}\equiv \Omega_{\eps'}^{-1}  \pa{\xi}\chi\pa{\xi}$ is the leading order solution of the Riemann-Hilbert problem corresponding to the choice of the radius $\eps'$ for the disks $\Dp{}D_{a_p,\eps'}$
centered at the singularities of the symbol $\sg\pa{\xi}$. We remind that the subscript $\eps$ in $\Omega$ indicates that the matrix $\Omega$ is the solution of the RHP
corresponding to the contour $\Sg_{\Omega_{\eps}}$  of section \ref{SectionEstResRHPOmega} where all the circles centered at $a_p$ have a radius $\eps$.

Note that the exact resolvent $R\pa{\xi,\eta}$ \textit{does not} depend on the choice of $\eps'$. Hence, we can chose different values for $\eps'$, depending on the point $\pa{\xi,\eta}$ where we want to estimate the resolvent.
According to the above remark, we can present the exact resolvent as
\beq
R\pa{\xi,\eta}=R_{0}\pa{\xi,\eta}+R_c\pa{\xi,\eta} \; .
\enq
There we have introduced the correcting resolvent $R_c$. This resolvent is defined in terms of matrices $\Omega_{\eps'}$ where $\eps'$ takes different values depending on the point where we are placed. More precisely,
\beq
R_c\pa{\xi,\eta;\f{\eps}{2}}=\pa{-f_-^{bk}\pa{\xi} \; f_+^{bk}\pa{\xi}}\f{ \Omega_{\f{\eps}{2}}^{-1}\pa{\xi} \Omega_{\f{\eps}{2}}\pa{\eta}-I_2}{\xi-\eta} \pa{\ba{c} f_+^{bk}\pa{\eta} \\ f_-^{bk}\pa{\eta} \ea  }  \;\;
\label{correction resolvent bulk}
\enq
for $\eta\; , \xi \in \R\setminus \bigcup\limits_{p=1}^{n} \intff{a_p-\eps}{a_p+\eps}$, and
\beq
R_c\pa{\xi,\eta;2\eps}=\pa{-f_{-;p}^{loc}\pa{\xi} \; f_{+;p}^{loc}\pa{\xi}}\f{ \Omega_{2\eps}^{-1}\pa{\xi} \Omega_{2\eps}\pa{\eta}-I_2}{\xi-\eta} \pa{\ba{c} f_{+;p}^{loc}\pa{\eta} \\ f_{-;p}^{loc}\pa{\eta} \ea  }
\label{correcting resolvent at ap}
\enq
for $\eta\; , \xi \in \intoo{a_p-\eps}{a_p+\eps}$. The crucial point is that the correcting resolvent is indeed small in the sense of the integration:
\begin{prop}
Let $\beta_p\in\paa{\de_p,\ga_p}$. Then
\beq
\Int{\R}{}
\f{R\pa{\xi,\xi}-R_0\pa{\xi,\xi}}{\sg\pa{\xi}-1}\Dp{\beta_p}\sg\pa{\xi}  \dd \xi= \e{O}\pa{x^{\rho-1}}
\enq
\end{prop}
\Proof
It is a standard fact that a matrix $\Omega_{\eps}\pa{z}$ approaching the identity at $\infty$ and having a jump  $\Omega_{+;\eps} \pa{I_2+\De}=\Omega_{-;\eps}$
on $\Sg_{\Omega_{\eps}}$, can be expressed in terms of solutions of a singular integral equation:
\beq
\Omega_{\eps}\pa{\xi}=I_2+C_{\Sg_{\Omega_{\eps}}}\pac{\Omega_{+;\eps} \De}\pa{\xi} \quad \e{where} \quad \Omega_{+;\eps}\pa{\xi}=I_2+C_{\Sg_{\Omega_{\eps}}}^+\pac{\Omega_{+;\eps} \De}\pa{\xi} \; .
\label{Representation Cauchy pour solution SingIntEqn Omega}
\enq
It is readily seen from the asymptotics of CHF \eqref{Asymptotiques Humbert's CHF} that the jump matrix $I_2+\De$ for $\Omega_{\eps}$ is such that
$\norm{\De\pa{z}}_{L^{\infty}\pa{\Sg_{\Omega_{\eps}}}}=\e{O}\pa{x^{\rho-1}}$, but also $\norm{\De}_{L^{2}\pa{\Sg_{\Omega_{\eps}}}}=\e{O}\pa{x^{\rho-1}}$.  In the above equation, we have introduced the Cauchy operator
\beq
C_{\Sg_{\Omega_{\eps}}}\pac{f}\pa{\xi}=\Int{\Sg_{\Omega_{\eps}}}{} \f{\dd s}{2i\pi \pa{\xi-s}} f\pa{s} \quad ,\quad  C^+_{\Sg_{\Omega_{\eps}}}\pac{f}\pa{\eta}=\lim_{\xi\tend \eta^+} C_{\Sg_{\Omega_{\eps}}}\pac{f}\pa{\xi} \; .
\enq
$\xi\tend \eta^+$ means that $\xi$ approaches $\eta\in \Sg_{\Omega_{\eps}}$ non-tangentially from the left side of $\Sg_{\Omega_{\eps}}$.
It is a classical result that $C^+_{\Sg_{\Omega_{\eps}}}$ is a continuous operator on $L^{2}\pa{\Sg_{\Omega_{\eps}}}$. The fact that
$\norm{\De}_{L^2\pa{\Sg_{\Omega_{\eps}}}}=\e{O}\pa{x^{\rho-1}}$, leads to  $\norm{\Omega_{+;\eps}-I_2}_{L^2\pa{\Sg_{\Omega_{\eps}}}}=\e{O}\pa{x^{\rho-1}}$. In its turn, this means that the entries of $C$ such that
$\Omega_{\eps}\pa{z}=I_2+\tf{C}{z}+\e{O}\pa{z^{-2}}$ are a $\e{O}\pa{x^{\rho-1}}$.

We have that
\bem
\Int{\R}{}
\f{R\pa{\xi,\xi}-R_0\pa{\xi,\xi}}{\sg\pa{\xi}-1}\Dp{\beta_p}\sg\pa{\xi}  \dd \xi= \\
\Int{\R\setminus \bigcup\limits_{p=1}^{n}\intff{a_p-\eps}{a_p+\eps} }{} \hspace{-7mm} R_c\pa{\xi,\xi;\tf{\eps}{2}} \f{\Dp{\beta_p}\sg\pa{\xi}}{\sg\pa{\xi}-1}
+ \sul{p=1}{n}\Int{a_p-\eps}{a_p+\eps}R_c\pa{\xi,\xi;2\eps} \f{\Dp{\beta_p}\sg\pa{\xi}}{\sg\pa{\xi}-1} \; .
\label{calcul de la correction au remplacement du resolvent}
\end{multline}
There, we should substitute the values of the correcting resolvent $R_c$ given in \eqref{correcting resolvent at ap} and \eqref{correction resolvent bulk}.

The exact formula for the correcting resolvent is
\beq
R_c\pa{\xi,\xi;\eps'} \f{\Dp{\beta_p}\sg\pa{\xi}}{\sg\pa{\xi}-1}=\f{1}{2i\pi}\sul{i,j}{} f_{\mu_i}^{\pa{0}}\pa{\xi}f_{\mu_j}^{\pa{0}}\pa{\xi} \Dp{\beta_p}\sg\pa{\xi} \pa{\Omega_{\eps'}^{-1}\pa{\xi} \Dp{\xi}\Omega_{\eps'}\pa{\xi}}_{ij} \; .
\label{Resolvent correction sur l'intervall sans singularites}
\enq
There we agree upon $\mu_1=+$ and $\mu_2=-$, as well as on the fact that
\beq
f^{\pa{0}}_{\pm}\pa{\xi} =  \left\{\ba{cc}

f^{bk}_{\pm}\pa{\xi}  & \xi \in \mc{D}_{\e{bk}}= \R\setminus  \bigcup\limits_{k=1}^{n}\intff{a_k-\eps}{a_k+\eps} \vspace{1mm}\\
f^{loc}_{\pm;p}\pa{\xi}  & \xi \in \mc{D}_{\e{loc}} = \bigcup\limits_{k=1}^{n}\intff{a_k-\eps}{a_k+\eps} \ea \right. \;.
\enq
Also $\eps^{\prime}=\tf{\eps}{2}$ for $\xi \in \mc{D}_{\e{bk}}$ and  $\eps^{\prime}=2 \eps$ for $\xi \in \mc{D}_{\e{loc}}$. The choice of two possible values for $\eps'$  depending whether $\xi\in \mc{D}_{\e{bk}}$
or $\mc{D}_{\e{loc}}$ ensures that $\pa{\Omega_{\eps'}\pa{\xi}\Dp{\xi}\Omega_{\eps'}\pa{\xi}}_{ij}$ is smooth
on $\R$.

We first study the integral on $\R\setminus\cup_{p=1}^{n}\intff{a_p-\tf{\eps}{2}}{a_p+\tf{\eps}{2}}$. As already pointed out, the jump matrix $\De$ is uniformly $\e{O}\pa{x^{\rho-1}}$. This ensures that $\Omega\pa{\xi}-I_2$ is also an $\e{O}\pa{x^{\rho-1}}$, for $\xi$ bounded. Moreover, by Lemma \ref{lemme comportement Omega a l'infini},
\beq
\Omega_{\f{\eps}{2}}^{-1}\pa{\xi} \Dp{\xi} \Omega_{\f{\eps}{2}}\pa{\xi}=x^{\rho-1} \e{O}\pa{\abs{\xi}^{-\f{1+\min\pa{1,\kappa}}{2}} \pa{\ba{cc} 1&1\\ 1&1 \ea} } \; ,\quad
\e{for} \; \xi \tend \pm \infty \;.
\enq
As $f_{\pm}^{bk}\pa{\xi} \sim \ex{\pm i\f{x \xi}{2}}$
for $\xi \tend \infty$, and $\Dp{\beta_p}\sg\pa{\xi}=\e{O}\pa{\xi^{-1}}$, we get that
\eqref{Resolvent correction sur l'intervall sans singularites} is absolutely integrable on $\R\setminus\cup_{p=1}^{n}\intff{a_p-\tf{\eps}{2}}{a_p+\tf{\eps}{2}}$ and that the integral is a $\e{O}\pa{x^{\rho-1}}$.

It remains to study the integral of $\eqref{Resolvent correction sur l'intervall sans singularites}$
on $\intff{a_p-\eps}{a_p+\eps}$. $\Omega_{2\eps}\pa{\xi}$ is smooth on this interval and equal to  $I_2+\e{O}\pa{x^{\rho-1}}$. Moreover, for $\xi\tend a_p$,
\beq
\Dp{\beta_p}\sg\pa{\xi}=\e{O}\pa{\abs{\xi-a_p}^{-2\ga_p} \ln \abs{\xi-a_p}  }
\label{estimation grand 0 pour diff beta p de sigma}
\enq
for $\beta_p\in\paa{\de_p,\ga_p}$. The formula for $R_c$ involves four terms. They can all be treated similarly so we only explain here how to estimate the contribution of the integral around $a_p$ involving the
$\pac{f_{+;p}^{loc}}^2$. We denote $J_{+,+}$ this integral and have
\bem
\abs{J_{+,+}}\equiv\abs{\;\Int{a_p-\eps}{a_p+\eps} \Dp{\beta_p} \sg\pa{\xi} \pa{\Omega_{2\eps}^{-1}\pa{\xi}\Dp{\xi}\Omega_{2\eps}\pa{\xi}}_{21} \pac{ f_{+;p}^{loc} }^2\pa{\xi}
\dd \xi} \\
\leq x^{\rho-2-2\Re\pa{\de_p}} C'  \Int{-\eps x}{\eps x} \abs{z}^{-2\Re\pa{\ga_p}} \ln\pa{\tf{x}{\abs{z}}}
\abs{ \Phi\pa{-\de_p-\ga_p,1-2\ga_p;-iz} }^2 \, \dd z \;.
\label{estimation integrale J++}
\end{multline}
We have used the bounds on all the smooth functions appearing in the first line and used the $\e{O}$ estimates for $\Dp{\beta_p} \sg$ \eqref{estimation grand 0 pour diff beta p de sigma}. $C'$ is some constant that can depend on $\eps$ and we have explicitly extracted the factor $x^{\rho-1}$ from the bound on the product of $\Omega$ matrices.

The integral in the last line of \eqref{estimation integrale J++} is divergent for $x\tend +\infty$ so that its leading  $x \tend +\infty$ asymptotics
are obtained by substituting the $z \tend \infty$ behavior of the integrand. Since
\beq
\Phi\pa{-\de_p-\ga_p,1-2\ga_p;-iz}= c_{\pm} z^{\de_p+\ga_p}\pa{1+\e{o}\pa{1}} \;\; \e{for}\; z\tend \pm \infty \;\;, \;c_{\pm}\in \Cx \; ,
\enq
we have
\bem
\Int{-\eps x}{\eps x} \abs{z}^{-2\Re\pa{\ga_p} } \ln\pa{\tf{x}{\abs{z}}}
\abs{ \Phi\pa{-\de_p-\ga_p,1-2\ga_p;-iz} }^2 \, \dd z  =\e{O}\pa{\Int{-\eps x}{\eps x} \ln \paf{x}{\abs{t}} \;  \abs{t}^{2\Re\pa{\de_p}}  \dd t} \\
=\e{O}\pa{\f{\pa{\eps x}^{2\Re\pa{\de_p}+1}}{2\Re\pa{\de_p}+1}  \pa{-\ln \eps  +  \f{1}{2\Re\pa{\de_p}+1}  }  } \; .
\end{multline}
Hence we have that $\abs{J_{++}}\leq c x^{\rho-1}$ for some constant $c$, uniformly in $\de_p, \ga_p$.  $\qquad \Box$

\begin{lemme}
\label{lemme comportement Omega a l'infini}
$\forall \eps'>0$,
\beq
\Omega_{\eps^{\prime}}^{-1}\pa{\xi} \Dp{\xi} \Omega_{\eps^{\prime}}\pa{\xi} =
x^{\rho-1} \e{O}\pa{\abs{\xi}^{-\f{1+\e{min}\pa{1,\kappa}}{2}} \pa{\ba{cc} 1&1\\ 1&1 \ea} } \;\;\;,
\;\; \xi \tend \pm \infty \;\; .
\enq
Here, we remind that $\kappa$ is such that $F\pa{\xi}=1+\e{O}\pa{\abs{\xi}^{-\f{1+\kappa}{2}}}$ when $\xi \tend \pm \infty$.

\end{lemme}
\Proof
As already discussed, we have that $\forall \eps^{\prime}>0$, $\Omega_{\eps'}\pa{\xi}=I_2+\e{o}\pa{1}$ in the
$L^{2}\cap L^{\infty}\pa{\Sg_{\Omega_{\eps^{\prime}}}}$ sense. Hence, it remains to study the asymptotic behavior
of $\Dp{\xi}\Omega_{\eps^{\prime}}\pa{\xi}$.

We focus on the case $\xi \tend +\infty$ as the $\xi \tend - \infty$ case can be treated similarly.  We decompose the  contour
$\Sg_{\Omega_{\eps^{\prime}}}=\Sg^{L}_{\Omega_{\eps^{\prime}}}\cup\Sg^{R}_{\Omega_{\eps^{\prime}}}$, with
\beq
\Sg^{L}_{\Omega_{ \eps^{\prime} }}=\Sg_{\Omega_{ \eps^{\prime} }} \cap \paa{z \; : \; \Re\pa{z} < \f{\xi}{2} } \quad , \quad
\Sg^{R}_{\Omega_{\f{\eps}{2}}}=\Sg_{\Omega_{\f{\eps}{2}}} \cap \paa{z \; : \; \Re\pa{z} \geq \f{\xi}{2} } \; .
\enq
Then, using the Cauchy integral  representation for $\Omega_{\eps'}$ \eqref{Representation Cauchy pour solution SingIntEqn Omega}, we get
\beq
\abs{\pac{\Dp{\eps'}\Omega_{\eps'}\pa{\xi} }_{ij} } \leq \f{2}{\pi \xi^2} \norm{\Omega_{+;\eps'}}_{L^{2}\pa{\Sg_{\Omega_{\eps'}}}} \norm{\De}_{L^{2}\pa{\Sg_{\Omega_{\eps'}}}} + \abs{I^{R}_{ij}} \; .
\label{estimation derivee partielle omega}
\enq
Where $I+\De$ stands for the jump matrix for $\Omega_{\eps'}$ and
\beq
I^R=\Int{\Sg^R_{\Omega_{\eps'}}}{} \f{\dd s}{2i\pi \pa{\xi-s}^2} \Omega_{+;\eps'} \De\pa{s} \;\; .
\enq
Reminding that $\norm{A}^2_{L^{2}\pa{\Sg_{\eps^{\prime}}}}=\Int{ \Sg_{\eps^{\prime}} }{} \abs{\dd s} \e{tr}\pac{A^{\dagger}\pa{s}A\pa{s}} $, we get
\beq
\abs{I^R_{ij}}^2 \leq \f{1}{4\pi^2}\norm{\Omega_{+;\eps'}}^{2}_{L^{2}\pa{\Sg_{\Omega_{\eps'}}}}
\Int{\Sg_{\Omega_{\eps'}}^R}{} \f{\abs{\dd s}}{\abs{\xi-s}^4} \e{tr}\paa{ \De^{\dagger}\pa{s}\De\pa{s}} \; .
\enq
But taking the explicit formula for $\De\pa{s}$ and  using the estimates for the asymptotic behavior of  $\sg_{\nu_k,\ov{\nu}_k}\pa{\xi}$, as well as the assumptions on the asymptotic behavior of the function $F\pa{\xi}=1+\e{O}\pa{\abs{\xi}^{-\f{1+\kappa}{2}}}$, for some $\kappa>0$, we get that
\beq
\e{tr}\paa{ \De^{\dagger}\pa{s}\De\pa{s} }=  \f{ \ex{-2x \de} }{s^{1+\min\pa{1,\kappa}}}
\pa{1+\e{o}\pa{1}} \;\;\; \e{at} \quad s\tend \infty .
\enq
Where $\de=  \e{dist}\pa{\Sg^R_{\Omega_{\eps'}},\R } >0$ uniformly in $\abs{\xi}$ large enough. Therefore, for such $\xi$
\beqa
\Int{\Sg_{\Omega_{\eps'}}^R}{} \f{\abs{\dd s}}{\abs{\xi-s}^4} \e{tr}\paa{ \De^{\dagger}\pa{s}\De\pa{s}} &\leq&
2 \Int{\f{\xi}{2}}{+\infty} \dd s   \f{ \ex{-2x \de} }{ \abs{\xi-s-i\de}^4 \pa{\abs{s+i\de} }^{1+\min\pa{1,\kappa}}} \\
&=& \ex{-2x \de } \e{O}\paf{1}{ \xi^{1+\e{min}\pa{1,\kappa}}  } \; .
\eeqa
It follows that
\beq
\abs{I^R_{ij}} =  \ex{-2x \de } \e{O}\paf{1}{ \xi^{1+\e{min}\pa{1,\kappa}}  } \;,
\enq
for some constant $c$. The Lemma follows once upon replacing $\norm{\De}_{L^{2}\pa{\Sg_{\Omega_{\eps'}}  }}=\e{O}\pa{x^{\rho-1}}$ in \eqref{estimation derivee partielle omega}. $\qquad \Box$

\section{Asymptotics of the truncated Wiener-Hopf determinant}
\label{Section asymptotiques du determinant}
In this Section we shall compute the leading asymptotics of determinants of truncated Wiener-Hopf operators that are Hilbert-Schmidt.

\noindent Note that, for $\xi \in \R$, we can decompose $\sg_{\nu_k,\ov{\nu_k}}\pa{\xi-a_k}$ into
\beq
\sg_{\nu_k,\ov{\nu}_k}\pa{\xi-a_k} = \ex{g_{k}\pa{\xi}}  \ex{h_{k}\pa{\xi}} \; .
\enq
\noindent Here,
\beqa
g_{k}\pa{\xi}&=& \de_k \pac{\ln\pa{\xi-a_k-i}  - \ln \pa{\xi-a_k+i} + 2i\pi \Xi\pa{a_k-\xi}} \; , \\
h_{k}\pa{\xi}&=&\ga_k \pac{\ln\pa{\xi-a_k+i}+\ln\pa{\xi-a_k-i}-2\ln \abs{\xi-a_k}} \; ,
\eeqa
and we remind that $\Xi$ is Heaviside's step function. It follows,
\beq
\Dp{\beta_k}\sg \pa{z}=  \sg\pa{z} \Dp{\beta_k} \pa{g_k\pa{z} + h_k\pa{z}} \;\;\;
\e{for} \quad \beta_k\in \paa{\de_k, \ga_k} .
%
%\ln g_{\de_p}\pa{\xi-a_p}  h_{\ga_p}\pa{\xi-a_p}\f{\ex{-2 i \pi \de^{\pa{p}}H\pa{a_p-z}}}{\abs{z-a_p}^{2\ga^{\pa{p}}}} \qquad z \in \R
%
\enq

\subsection{Integration}
Before carrying out the integrals appearing in the formulas for $\Dp{\de_1}\ln\ddet{2}{I+V_2}$ and $\Dp{\ga_1}\ln\ddet{2}{I+V_2}$, we first establish a useful integration Lemma.

\begin{lemme}
\label{lemme preparatoire}
Let \mc{R} be a Riemann integrable function on $\R$, $g\in \msc{C}^{1}\pa{I,\Cx}$,
and $I$ an interval such that $\overset{\circ}{I}\supset \intff{-\eps_1}{\eps_2}$, where $\eps_1\geq 0$ , $\eps_2\geq 0$.
 Then
\beq
\Int{-\eps_1}{\eps_2}x \dd t g\pa{t} \mc{R}\pa{xt}= g\pa{0}
\Int{-\underset{x\tend + \infty}{\lim} x \eps_1}{ \underset{x\tend + \infty}{\lim} x \eps_2}
 \mc{R}\pa{t} + \e{o}\pa{1} \; .
\label{asymptotics of Riemann Intergrable integrals}
\enq
The $\e{o}\pa{1}$ corresponds to terms vanishing in the ordered limit $x \tend + \infty$ and then $\eps_i \tend 0$.
\end{lemme}
\Proof

One has
\beq
\Int{-\eps_1}{\eps_2} \pa{g\pa{t}-g\pa{0}} \mc{R}\pa{x t} x \, \dd t =
\Int{0}{\eps_2} \dd y g'\pa{y} \Int{xy}{x \eps_2} \dd t \mc{R}\pa{t} - \Int{-\eps_1}{0} \dd y g'\pa{y} \Int{-x \eps_1}{xy} \dd t \mc{R}\pa{t}
\label{decomposition dans integration lemma}
\enq
\noindent Since $\mc{R}$ is Riemann integrable, we have that  $\pa{a,b} \mapsto \Int{a}{b} \mc{R}\pa{t} \dd t$
is continuous and has a finite limit at $\infty$. Hence it is bounded, say by $M>0$. It follows from \eqref{decomposition dans integration lemma} that
\beq
\abs{\Int{-\eps_1}{\eps_2}x \dd t \pa{g\pa{t}-g\pa{0}} \mc{R}\pa{t} } \leq \pa{\eps_1+\eps_2} M \, \underset{\intff{-\eps_1}{\eps_2}}{\max} {\abs{g'}}  \;\;  \qquad  \qquad \Box .
\enq

We now establish the first separation identity.

\begin{prop}
\label{proposition derivee delta}
Let $I+V$ be the generalized sine kernel defined by a symbol $\sg-1$ satisfying the $L^{2}\pa{\R}$ assumptions. Then
the below identity holds
\bem
\ln \paf{\ddet{2}{I+V}} {\ddet{2}{I+V}_{\mid_{\de_1=0}}  } =
x\Int{\R}{} \f{\dd \xi}{2\pi} \paa{\ln \paf{\sg}{{\sg}_{\mid_{\de_1=0}}}  -\sg + {\sg}_{\mid_{\de_1=0}}  }  +\de_1^2 \ln \paf{2}{x}\\
 + \de_1 \ln \paf{\a_{\ua}^{\pa{2\dots n}}\pa{a_1+i}  \a_{\da}^{\pa{2\dots n}}\pa{a_1-i}  }
 {\a_{\ua}^{\pa{2\dots n}}\pa{a_1+i0^+}  \a_{\da}^{\pa{2\dots n}}\pa{a_1-i0^+}  }  \\
+ \ln\paf{ G\pa{1-\ga_1+\de_1} G\pa{1-\ga_1-\de_1}}{G\pa{1-\ga_1} G\pa{1-\ga_1} } +\e{o}\pa{1} \; .
\end{multline}
There, G stands for the Barnes' function and the $\e{o}\pa{1}$ vanishes in the $x\tend +\infty$ limit. The $\e{o}\pa{1}$ is uniform in $\de_p$ and $\ga_p$.

The functions $\a_{\ua/\da}^{\pa{2\dots n}}$ refer to the values in the upper/lower half planes of the solution to the scalar RHP for $\a$ with
vanishing exponents $\nu_1$ and $\ov{\nu}_1$, \textit{ie}
\beqa
\a_{\ua}^{\pa{2\dots n}}\pa{z} &=& F_+^{-1}\pa{z} \pl{k=2}{n} \paf{z-a_k}{z-a_k+i}^{\nu_k}
\quad  ,  \nonumber\\
\a_{\da}^{\pa{2\dots n}}\pa{z} &=&  F_-\pa{z} \pl{k=2}{n} \paf{z-a_k-i}{z-a_k}^{\ov{\nu}_k} \; .
\label{definition alpha updown site 2}
\eeqa
Finally,  the notation $\mid_{\de_1=0}$ means that we ought to set $\de_1=0$ without altering all the other parameters.
\end{prop}
\Proof
The asymptotic formulae for the resolvent \eqref{resolvent autour des singlarites}, \eqref{Resolvent sur interval non singulier} combined with \eqref{derivee partielle HilbertSchmidt} lead to
\bem
\Dp{\de_1}\ln \ddet{2}{I+V}= \Int{\R}{}
\paa{ \f{R_0\pa{\xi,\xi}   }{\sg\pa{\xi}-1}\Dp{\de_1}\sg\pa{\xi} -\Dp{\de_1}V\pa{\xi,\xi}} \dd \xi\\
+
\Int{\R}{}
\f{R\pa{\xi,\xi}-R_0\pa{\xi,\xi}}{\sg\pa{\xi}-1}\Dp{\de_1}\sg\pa{\xi}  \dd \xi \; .
\nonumber
\end{multline}
We have already shown in Proposition 5.2 that the second line is a $\e{O}\pa{x^{\rho-1}}$, or more precisely, a term going to $0$ in the ordered limit $x\tend +\infty$ and then $\eps \tend 0$.
We cannot send $\eps$ to zero first, as the constant in the estimate for the $\e{O}$ also depends on $\eps$, and goes to infinity when $\eps \tend 0$. Writing
the result of integration explicitly yields
\bem
\Dp{\de_1}\ln \ddet{2}{I+V}= x \Int{J_{\eps}}{} \f{\dd \xi}{2\pi} \paa{\Dp{\de_1} \ln \sg - \Dp{\de_1} \sg}
-\Int{ J_{\eps} }{} \f{ \dd \xi }{ 2i\pi } \pa{  \Dp{\xi}\ln \a_+ \a_-  }\pa{\xi} \Dp{\de_1}g_{1}\pa{\xi}  \\
%
%
%-x \sul{p=1}{n} \Int{a_p-\eps}{a_p+\eps} \f{\dd \xi}{2\pi} \Dp{\de_1} \sg \pa{\xi} +
%
 \sul{p=1}{n} \Int{a_p-\eps}{a_p +\eps}\!\!\!\! \paa{ -x \Dp{\de_1} \sg \pa{t}  + \Dp{\de_1} g_{1} \pa{t}
\pac{x-\f{2ix\de_p}{x\pa{t-a_p}+\e{sgn}\pa{t-a_p}} +i \Dp{t}\ln \wh{\a}^{\pa{p}}_{\ua,+} \wh{\a}^{\pa{p}}_{\da,-} \pa{t}  }} \f{\dd t}{2\pi} \\
+\sul{p=1}{n} \Int{-\eps}{\eps} \f{\dd t}{2i\pi }
\pac{\Dp{\de_1}\wt{g}_{1}\pa{t+a_p}} \pac{-ix \wt{\tau}\pa{\ga_p,\de_p;xt}
- \pa{\Dp{\xi}\ln \wh{\a}^{\pa{p}}_{\ua,+} \wh{\a}^{\pa{p}}_{\da,-}}\pa{t+a_p} \wt{\varphi}\pa{\ga_p,\de_p; xt} } \\
+\Int{-\eps}{0} \dd t
 \pac{-ix \wt{\tau}\pa{\ga_1,\de_1;xt}
- \pa{\Dp{\xi}\ln \wh{\a}^{\pa{1}}_{\ua,+} \wh{\a}^{\pa{1}}_{\da,-}}\pa{t+a_1} \wt{\varphi}\pa{\ga_1,\de_1; xt} } \;  +\e{O}\pa{x^{\rho-1}}
\label{asymptotique de Dplogdet en delta}
\end{multline}
There $J_{\eps}=\R \setminus \bigcup_{p=1}^{n} \intff{a_p-\eps}{a_p+\eps}$. Also, we have introduced $\wt{g}_{1}$,
the smooth part of  $g_{1}$, as well as $\wt{\tau}$ and $\wt{\varphi}$, the Riemann integrable regularizations of
the functions $\tau$ and $\varphi$ introduced in \eqref{definition de varphi}, \eqref{definition de tau}:
\beqa
\wt{g}_{1}\pa{\xi}&=&\de_1\pa{\ln \pa{\xi-a_1-i} - \ln \pa{\xi-a_1+i}   } \\
\wt{\tau}\pa{\ga_p,\de_p;xt}&=& \ex{-i\pi \de_p \e{sgn}\pa{t}} \abs{t}^{-2\ga_p}\tau\pa{\ga_p,\de_p;xt}+1-\f{2i \de_p}{xt + \e{sgn}\pa{t}}  \\
 \wt{\varphi}\pa{\ga_1,\de_1; xt} &=& \ex{-i\pi \de_p \e{sgn}\pa{t}} \abs{t}^{-2\ga_p}\varphi\pa{\ga_p,\de_p;xt}-1
\label{longue expression pour derivee delta}
\eeqa
The Riemann integrability of $\wt{\tau}$ and $\wt{\varphi}$ is part of the conclusions of Corollary \ref{Corolaire integrales CHF} given in Appendix \eqref{Appendix integration CHF}.

The integrals in the third and fourth line of \eqref{longue expression pour derivee delta} can be estimated using Lemma 7. One gets that
\bem
\Int{-\eps}{0} \dd t
 \pac{-ix \wt{\tau}\pa{\ga_1,\de_1;xt}
- \pa{\Dp{\xi}\ln \wh{\a}^{\pa{1}}_{\ua,+} \wh{\a}^{\pa{1}}_{\da,-}}\pa{t+a_1} \wt{\varphi}\pa{\ga_1,\de_1; xt} } \\
= -i\Int{-\infty}{0}\wt{\tau}\pa{\ga_1,\de_1;t}\dd t - x^{-1}\pa{\Dp{\xi}\ln \wh{\a}^{\pa{1}}_{\ua,+} \wh{\a}^{\pa{1}}_{\da,-}}\pa{a_1}
\Int{-\infty}{0}\wt{\varphi}\pa{\ga_1,\de_1;t}\dd t  \; +  \; \e{o}\pa{1}\\
=-2\de_1+\pa{\ga_1+\de_1}\psi\pa{-\ga_1-\de_1}+\pa{\de_1-\ga_1}\psi\pa{\de_1-\ga_1} + i\pi\ga_1 +\e{o}\pa{1}
\end{multline}
There we have used the value of the integrals of $\wt{\tau}\pa{\ga,\de;t}$ and $\wt{\varphi}\pa{\ga,\de;t}$
given in Corollary \ref{Corolaire integrales CHF}. The $\e{o}\pa{1}$ stands for terms that vanish
in the ordered limit $x\tend +\infty$ and then $\eps \tend 0$. Very similarly
\bem
\sul{p=1}{n} \Int{-\eps}{\eps} \f{\dd t}{2i\pi }
\pac{\Dp{\de_1}\wt{g}_{1}\pa{t+a_p}} \pac{-ix \wt{\tau}\pa{\ga_p,\de_p;xt}
- \pa{\Dp{\xi}\ln \wh{\a}^{\pa{p}}_{\ua,+} \wh{\a}^{\pa{p}}_{\da,-}}\pa{t+a_p} \wt{\varphi}\pa{\ga_p,\de_p; xt} }  \\
=\sul{p=1}{n} -\f{i}{2i\pi} \Dp{\de_1}\wt{g}\pa{a_p} \times -2\pi \ga_p
= \sul{p=1}{n} \ga_p  \Dp{\de_1}\wt{g}\pa{a_p}   \; .
\end{multline}
\noindent The equality
\beq
\sul{p=1}{n} \Int{a_p-\eps}{a_p +\eps} i\Dp{\de_1} g_{1} \pa{t}
\Dp{t}\ln \pa{ \wh{\a}^{\pa{p}}_{\ua,+} \wh{\a}^{\pa{p}}_{\da,-} \pa{t} }
\f{\dd t}{2\pi}=\e{O}\pa{\eps}
\enq
holds as we deal with a finite sum of integrals of piecewise smooth functions over intervals of length $2\eps$. In order to estimate
the  integral appearing in the second line of \eqref{longue expression pour derivee delta},
 we need to change its expression a little
\bem
 \sul{p=1}{n} \Int{a_p-\eps}{a_p +\eps} \Dp{\de_1} g_{1} \pa{t}
\f{2ix\de_p}{x\pa{t-a_p}+\e{sgn}\pa{t-a_p}}  \f{\dd t}{2\pi} \\
= 2i\pi\Int{-\eps}{0} \f{2\de_1 x }{ xt -1} \f{\dd t }{2\pi}
+ \Int{0}{\eps} \sul{p=1}{n} \Dp{\de_1} \paf{ \wt{g}_{1} \pa{t+a_p}- \wt{g}_{1} \pa{a_p-t} }{t} \f{2i\de_p xt}{xt+1}
\f{\dd t }{2\pi}
  \\
=-2\de_1\ln\pa{x\eps+1} +\e{O}\pa{\eps} = -2\de_1\ln{x\eps} +\e{o}\pa{1}
\label{somme qui donne le log epsilon}
\end{multline}
There we have used that $t^{-1}\pa{\wt{g}_{1} \pa{t+a_p}- \wt{g}_{1} \pa{a_p-t}}$ is smooth and $\tf{xt}{\pa{xt+1}}\leq 1$. This ensures that the corresponding integrals are a $\e{O}\pa{\eps}$. The $\e{o}\pa{1}$ term stands, once again, for terms vanishing in the ordered limit $x\tend+\infty$ and then $\eps\tend 0$.

We now explain how to treat the integral containing $\Dp{z}\ln\pa{\a_+ \a_-}\pa{z}$. There, the limit $\eps \tend 0$
cannot be taken directly as the integrand has non-integrable singularities at the $a_k$. We thus start
by integrating by parts:
\bem
-\Int{ J_{\eps} }{} \f{ \dd \xi }{ 2i\pi } \pa{  \Dp{\xi}\ln \a_+ \a_-  }\pa{\xi} \Dp{\de_1}g_{1}\pa{\xi}=
 \Int{J_{\eps}}{} \f{\dd \xi}{2i\pi} \ln \a_+\a_- \pa{\xi} \Dp{\de_1 , \xi}^2 \wt{g_1}\pa{\xi} \;
\\
 -\f{1}{2i\pi}\sul{p=1}{n} \Dp{\de_1} \wt{g}_{\de_1}\pa{a_p-\eps} \ln \pac{\a_+\a_-}\pa{a_p-\eps}
-  \Dp{\de_1}\wt{g}_{1}\pa{a_p+\eps} \ln \pac{\a_+\a_-}\pa{a_p+\eps} \\
- \, \ln\pa{\a_+\a_-}\pa{a_1-\eps} \; .
%
%
%+ \, \ln\pa{\a_+\a_-}\pa{a_1-\eps} + \sul{p=1}{n} \pac{ \pa{\Dp{\de_1} \wt{g}_{\de_1}} \ln \a_+\a_- }\pa{a_p-\eps} -  %\pac{ \pa{\Dp{\de_1}\wt{g}_{1}} \ln \a_+\a_- }\pa{a_p+\eps}
%
%
\label{calcul de l'integrale de alpha + et alpha -}
\end{multline}
There, we have explicitly separated the regular part $\Dp{\de_1}\wt{g}$ of $\Dp{\de_1}g$ from the one containing a jump. The sum appearing in the second line can be computed up to $\e{O}\pa{\eps\ln\eps}$ terms by using the local
behavior of $\a_+\a_-$ around $a_p$:
\beq
\a_+\a_-\pa{a_p\mp\eps}=\wh{\a}_{\ua}^{\pa{p}}\pa{a_p\mp\eps+i0^+}\wh{\a}_{\da}^{\pa{p}}\pa{a_p\mp\eps-i0^+} \ex{2i\pi \ga_p \Xi\pa{\pm 1}}\abs{\eps}^{-2\de_p}
\enq
Hence,  up to  $\e{O}\pa{\eps \ln\eps}$ terms, only the discontinuous part of $\a_+\a_-$ contributes to the sum:
\bem
 -\f{1}{2i\pi}\sul{p=1}{n} \paa{ \pac{\Dp{\de_1} \wt{g}_{\de_1} \ln \a_+\a_-}\pa{a_p-\eps}
-  \pac{\Dp{\de_1}\wt{g}_{1} \ln \a_+\a_-}\pa{a_p+\eps} } \\
= - \sul{p=1}{n}\ga_p \Dp{\de_1} \wt{g}_{\de_1}\pa{a_p}
+\e{O}\pa{\eps\ln \eps} \; .
\end{multline}
The singular behavior of $\a_+\a_-$ around $a_1$ leads to
\beq
-  \ln\pac{\a_+\a_-}\pa{a_1-\eps}  = 2\de_1 \ln \eps -i\pi \ga_1
- \ln \pac{\a^{\pa{2,\dots,n}}_{\ua}\pa{a_1+i0^+}\a^{\pa{2,\dots,n}}_{\da}\pa{a_1-i0^+}} \; .
\enq
 $\a^{\pa{2,\dots,n}}_{\ua/\da}$ have been defined in \eqref{definition alpha updown site 2}.

Up to $\e{o}\pa{1}$ corrections, it is now possible to replace the integral over $J_{\eps}$ appearing in the $rhs$ of \eqref{calcul de l'integrale de alpha + et alpha -} by one over $\R$ as the integrand has integrable singularities at the points $a_k$.  The resulting integral over $\R$ can then be evaluated by computing the residues at $\xi=a_1\pm i$
thanks to the fact that $\a_{\ua/\da}$ is analytic in the upper/lower half-plane and goes to 1
when $z\tend \infty$ in $\mc{H}_{\pm}$. One gets:
\beq
\Int{\R}{} \f{\dd \xi}{2i\pi} \ln\pa{\a_+\a_-}\pa{\xi} \Dp{\de_1\, , \xi}^2 \wt{g_1}\pa{\xi} =
\ln\pa{\a_+^{\pa{2\dots n}}\pa{a_1+i} \a_-^{\pa{2\dots n}}\pa{a_1-i}}+2\de_1\ln 2 \; .
\enq
At the end of the day we get,
\bem
-\Int{ J_{\eps} }{} \f{ \dd \xi }{ 2i\pi } \pa{  \Dp{\xi}\ln \a_+ \a_-  }\pa{\xi} \Dp{\de_1}g_{1}\pa{\xi}=
-\ln\paf{ \a_+^{\pa{2\dots n}}\pa{a_1} \a_-^{\pa{2\dots n}}\pa{a_1} }
{\a_+^{\pa{2\dots n}}\pa{a_1+i} \a_-^{\pa{2\dots n}}\pa{a_1-i}} \\
+2\de_1\ln 2\eps - i\pi \ga_1- \sul{p=1}{n}\Dp{\de_1} \wt{g}_1\pa{a_p} \ga_p  \; +\e{o}\pa{1}\; .
\label{integrale sur Jepsilon avec delta}
\end{multline}

 Putting all the different results together we get
\bem
\Dp{\de_1}\ln \ddet{2}{I+V_2}= x \Int{\R}{} \f{\dd \xi}{2\pi} \paa{\Dp{\de_1} \ln \sg - \Dp{\de_1} \sg}
+2\de_1\ln \paf{2}{x} -2\de_1\\
- \ln\paf{ \a_+^{\pa{2\dots n}}\pa{a_1} \a_-^{\pa{2\dots n}}\pa{a_1} }
{\a_+^{\pa{2\dots n}}\pa{a_1+i} \a_-^{\pa{2\dots n}}\pa{a_1-i}}
+\pa{\ga_1+\de_1}\psi\pa{-\ga_1-\de_1}  \\
+\pa{\de_1-\ga_1}\psi\pa{\ga_1-\de_1} + \e{o}\pa{1} \; .
\label{derivee partielle pour logdet en delta}
\end{multline}
\noindent In particular, the two  $\ln \eps$ contributions from \eqref{integrale sur Jepsilon avec delta} and
\eqref{somme qui donne le log epsilon} cancel each other. We stress that $\e{o}\pa{1}$ stands for vanishing terms in the ordered limit $x\tend +\infty$ and $\eps \tend 0$.
It remains to integrate this result from $0$  up to $\de_1$. This is licit as the remainders $\e{o}\pa{1}$ are uniform in $\abs{\Re\pa{\de_1}}<\tf{1}{2}$. The integral of the $\psi$ functions yields Barnes' functions due to the formula \eqref{fonction de Barnes}.
Finally, we get that
\bem
\ln \paf{\ddet{2}{I+V}} {\ddet{2}{I+V}_{\mid_{\de_1=0}}  }
x\Int{\R}{} \f{\dd \xi}{2\pi} \paa{\ln \paf{\sg}{{\sg}_{\mid_{\de_1=0}}}  -\sg + {\sg}_{\mid_{\de_1=0}}  }  +\de_1^2 \ln \paf{2}{x}\\
= \de_1 \ln \paf{\a_{\ua}^{\pa{2\dots n}}\pa{a_1+i}  \a_{\da}^{\pa{2\dots n}}\pa{a_1-i}  }
 {\a_{\ua}^{\pa{2\dots n}}\pa{a_1+i0^+}  \a_{\da}^{\pa{2\dots n}}\pa{a_1-i0^+}  }  \\
+ \ln\paf{ G\pa{1-\ga_1+\de_1} G\pa{1-\ga_1-\de_1}}{G\pa{1-\ga_1} G\pa{1-\ga_1} } +\e{o}\pa{1} \; .
\label{expression pour determinant ration at any delta}
\end{multline}
This means that, for all $\eps>0$, the limit
\beq
\!\! \underset{x \tend +\infty}{\lim}\!\!  \paa{\!\! \ln \paf{\ddet{2}{I+V}} {\ddet{2}{I+V}_{\mid_{\de_1=0}} } -
x\Int{\R}{} \f{\dd \xi}{2\pi} \paa{\ln \paf{\sg}{\sg_{\mid_{\de_1=0}}}  -\sg + \sg_{\mid_{\de_1=0}}  } + \de_1^2 \ln x}
\enq
exists. It is given by the $rhs$ of \eqref{expression pour determinant ration at any delta}, where $\e{o}\pa{1}$
are $\eps$ dependent terms that go to $0$ when $\eps\tend 0$. As the $lhs$ of \eqref{expression pour determinant ration at any delta} is $\eps$-independent the $x\tend +\infty$
limit cannot depend on $\eps$, therefore the value of the constant can be computed by sending $\eps\tend0$. The claim then follows.
$\Box$

\begin{prop}
\label{proposition derivee gamma}
Under the assumptions of the previous proposition, the following identity holds
\bem
\ln \paf{\ddet{2}{I+V}_{\mid_{\de_1=0}} } {\ddet{2}{I+V}_{\mid_{\de_1=\ga_1=0}}  } =
x\Int{\R}{} \f{\dd \xi}{2\pi} \paa{\ln \paf{\sg_{\mid_{\de_1=0}}}{\sg_{\mid_{\de_1=\ga_1=0}}}  -\sg_{\mid_{\de_1=0}} + \sg_{\mid_{\de_1=\ga_1=0}}  }  \\
+ \ga_1^2 \ln \paf{x}{2}
+ \ga_1 \ln \!\!\paf{\!\!\a_{\ua}^{\pa{2\dots n}} \pa{a_1+i}  \a_{\da}^{\pa{2\dots n}}\pa{a_1-i0^+}  }
{\a_{\ua}^{\pa{2\dots n}}\pa{a_1+i0^+} \a_{\da}^{\pa{2\dots n}} \pa{a_1-i}  } +
\ln\paf{ G\pa{1-\ga_1} G\pa{1-\ga_1}}{G\pa{1-2\ga_1} } +\e{o}\pa{1} \; .
\end{multline}
\end{prop}
\Proof
Following very analogous steps to the $\de_1$-derivative, one shows that
\bem
\Dp{\ga_1} \ln \ddet{2}{I+V}_{\mid_{\de_1=0}}= x \Int{\R}{} \f{\dd \xi}{2\pi} \paa{\Dp{\ga_1} \ln \sg_{\mid_{\de_1=0}} - \Dp{\ga_1} \sg_{\mid_{\de_1=0}}  }
\\- \Int{J_{\eps}}{} \f{\dd \xi}{2i\pi} \pac{ \Dp{\xi} \ln \pa{\a_+\a_-}}\pa{\xi} \Dp{\ga_1} h_{1}\pa{\xi }
  +\sul{p=2}{n} \ga_p \Dp{\ga_1} h_{1}\pa{a_p} \\
+\!\!\Int{-\eps}{\eps} \!\! \f{\dd t}{2i\pi} \Dp{\ga_1}h_{1}\pa{z+a_1}
 \pac{-ix \wt{\tau}\pa{\ga_1,0;xt}
- \pa{\Dp{\xi}\ln \wh{\a}^{\pa{1}}_{\ua,+} \wh{\a}^{\pa{1}}_{\da,-}}\pa{t+a_1} \wt{\varphi}\pa{\ga_1,0; xt} } \;
 + \e{o}\pa{1} .
\label{ecriture integrale pour gamma derivee}
\end{multline}
More precisely, we have replaced the integration of the exact resolvent by one involving the approximate one $R_0$
for the price of $\e{O}\pa{x^{\rho-1}}$ corrections. Then we have applied the integration Lemma 7 to estimate asymptotically the integrals around $\paa{a_k}_{k=2}^{n}$ that involve the CHF. These estimates produced the sum $\sul{p=2}{n} \ga_p \Dp{\ga_1} h_{1}\pa{a_p}$ and some $\e{o}\pa{1}$ corrections. These, as before, are vanishing in the ordered limit $x\tend +\infty$ and then $\eps\tend 0$.

The integral around $a_1$ should be considered separately as it is a little different in respect to the already studied integrals. We obtain
\bem
\Int{-\eps}{\eps} \f{\dd t}{2i\pi} \Dp{\ga_1}h_{1}\pa{t+a_1}
 \pac{-ix \wt{\tau}\pa{\ga_1,0;xt}
- \pa{\Dp{\xi}\ln \wh{\a}^{\pa{1}}_{\ua,+} \wh{\a}^{\pa{1}}_{\da,-}}\pa{t+a_1} \wt{\varphi}\pa{\ga_1,0; xt} } \\
=\e{o}\pa{1}+\e{o}\paf{\log x}{x\eps} + \Int{-x\eps}{x\eps} \f{\dd t}{2\pi} \log\paf{\abs{t}}{x}
\wt{\tau}\pa{\ga_1,0;t} \\
=  2 \ga_1 \pa{\psi\pa{1-\ga_1}  -2\psi\pa{1-2\ga_1}+1} - 2\ga_1 \ln x + \e{o}\pa{1} \; .
\end{multline}
During the estimation of the above integral, one finds that the contributions stemming from the regular part of $\Dp{\ga_1}h_1$ only produce $\e{o}\pa{1}$ corrections as it vanishes at $t=a_1$.
Also one gets that the integral of the irregular part (equal to $\log\abs{\xi}$) versus $\wt{\varphi}$
produces at most $\e{O}\pa{\tf{\log x}{x}}$ corrections. The remaining integral in the second line can be estimated thanks to corollary \ref{Corolaire integrales CHF}. One should however use that fact that as $\wt{\tau}\pa{\ga_1,0;t}$ decreases
at infinity as an oscillating factor dumped by $\tf{1}{t}$,
\beq
\ln x \Int{-\eps x}{\eps x} \f{\dd t }{2\pi} \wt{\tau}\pa{\ga_1,0;t}=\ln x \Int{\R}{} \f{\dd t }{2\pi} \wt{\tau}\pa{\ga_1,0;t} + \e{O}\pa{\f{\ln x}{x \eps}} \; .
\enq

It now remains to study the  limiting value of the integral in the second line of
  \eqref{ecriture integrale pour gamma derivee}.
 In the case of the $\ga_1$-derivative,  this integral should be handled with greater care. Indeed,
\bem
-\Int{J_{\eps}}{} \f{\dd \xi}{2i\pi} \pa{\Dp{\xi} \ln \a_+ \a_-}\pa{\xi} \Dp{\ga_1} h_1\pa{\xi}=  \\
- \ga_1 \Int{\R\setminus\intff{-\eps}{\eps}}{} \!\!\!\!\f{\dd \xi}{2i\pi}   \pa{\f{2i}{\xi^2+1} -
\f{1}{\xi-i0^+} +\f{1}{\xi+i0^+} } \pa{\ln\pa{\xi^2+1} -2\ln \abs{\xi}} \\
-\f{1}{2i\pi} \sul{p=2}{n} \pac{\ln{\pa{\a_{\ua,+}\a_{\da,-}}^{\pa{2\dots n}}} \Dp{\ga_1}h_1}\pa{a_p-\eps}
-\pac{ \ln{\pa{\a_{\ua,+}\a_{\da,-}}^{\pa{2\dots n}}} \Dp{\ga_1}h_1}\pa{a_p+\eps} \\
+\Int{\R}{} \ln \pac{\a_+ \a_-}^{\pa{2\dots n}} \pa{\xi+a_1}
\paa{ \f{i}{\pa{\xi-i}\pa{\xi-i0^+}}-\f{i}{\pa{\xi+i}\pa{\xi+i0^+}} } \f{\dd \xi}{2i\pi} \;\;.
\label{integration sur Jepsilon de a+a-}
\end{multline}
\noindent  We have integrated by parts and decomposed the result into the integration of the singular part
(line 2), and the regular part involving $\a_{\ua}^{\pa{2,\dots,n}}\a_{\da}^{\pa{2,\dots,n}}$ (lines 3 and 4).
The latter functions have already been defined in \eqref{definition alpha updown site 2}. We have
\beq
 \f{1}{\xi+i0^+}-\f{1}{\xi-i0^+}=2i\pi\de\pa{\xi}
\enq
in the sense of distributions ($\de\pa{\xi}$ stands for the Dirac mass at zero). We can thus drop this functions  from the
first integral appearing on the $rhs$ of the equation above. In particular, this integral is finite even if, a priori, it involves terms $\tf{\ln\abs{\xi}}{\pa{\xi\pm i0^+}}$ that are integrated at distance $\eps$ from zero.  The sum appearing in the second line
of\eqref{integration sur Jepsilon de a+a-} is handled similarly to the case of the $\de_1$-derivative, \textit{ie} by separating the
smooth/singular parts of $\a_{\pm}^{\pa{2,\dots,n}}$ around $a_p$ and then neglecting all the $\e{O}\pa{\eps\ln \eps}$ contributions. Finally, one can send $\eps$ to zero in the integral appearing in the
last line of \eqref{integration sur Jepsilon de a+a-}. This produces some corrections that go to zero with $\eps$
due to the integrability of the integrand. The resulting integral over $\R$ can then be computed by the residues at $\xi=\pm i 0^{+}$ and $\xi=\pm i$, exactly as it was done in the proof of proposition \ref{proposition derivee delta}.
At the end of the day,
\bem
-\Int{J_{\eps}}{} \f{\dd \xi}{2i\pi} \pa{\Dp{\xi} \ln \a_+ \a_-}\pa{\xi} \Dp{\ga_1} h_1\pa{\xi}= -2\ga_1 \ln 2
- \sul{p=2}{n} \ga_p \Dp{\ga_1}h_1\pa{a_p} \\
+ \ln\paf{  \a_{\ua}^{\pa{2\dots n}} \pa{a_1+i} \a_{\da}^{\pa{2\dots n}} \pa{a_1-i0^+}    }
{ \a_{\ua}^{\pa{2\dots n}} \pa{a_1+i0^+} \a_{\da}^{\pa{2\dots n}} \pa{a_1-i} }+ \e{o}\pa{1} \;\; .
\end{multline}
Hence,
\bem
\Dp{\ga_1}{} \ln \ddet{2}{I+V}_{\mid_{\de_1=0}}=x \Int{\R}{} \f{\dd \xi}{2\pi} \paa{\Dp{\ga_1} \ln \sg_{\mid_{\de_1=0}} - \Dp{\ga_1} \sg_{\mid_{\de_1=0}}  }
 + 2 \ga_1 \ln \paf{x}{2} \\
+ \ln\paf{  \a_{\ua}^{\pa{2\dots n}} \pa{a_1+i} \a_{\da}^{\pa{2\dots n}} \pa{a_1-i0^+}    }
{ \a_{\ua}^{\pa{2\dots n}} \pa{a_1+i0^+} \a_{\da}^{\pa{2\dots n}} \pa{a_1-i} }
+2\ga_1\pa{\psi\pa{1-\ga_1}-2\psi\pa{1-2\ga_1}+1} + \e{o}\pa{1}\;\; .
\label{formule finale avant integration gamma derivee}
\end{multline}
\noindent We now integrate \eqref{formule finale avant integration gamma derivee} with respect to $\ga_1$.
The operation preserves the $\e{o}\pa{1}$ symbols as they are uniform in $\ga_1$. The $\psi$ functions are integrated thanks to $\eqref{fonction de Barnes}$. Once upon integration,
sending first $x$ to infinity and then $\eps$ to zero settles the value of the constant term.
 $\quad \Box$

\subsection{Asymptotics of the Fredholm determinant}
\begin{theorem}
Let I+K be a truncated Wiener-Hopf operator acting on the segment $\intff{0}{x}$ and generated by the symbol
$\sg-1$ with
\beq
\sg\pa{\xi}=F\pa{\xi} \pl{k=1}{n}\sg_{\nu_k,\ov{\nu}_k}\pa{\xi-a_p} \;\; , \quad a_i \in \R \quad a_1<\dots < a_n\; ,
\enq
\noindent where
\begin{itemize}
\item $F$ is holomorphic and non-vanishing in some open neighborhood of the real axis
        such that $F-1\in L^{2}\pa{\R}$, and even $F\pa{\xi}-1=\e{O}\pa{\abs{\xi}^{-\f{1+\kappa}{2}}}$, for some $\kappa>0$;
\item $\Re\pa{\ga_k}< \tf{1}{4}$ and $\abs{\Re\pa{\de_k}}<\tf{1}{2} \;,\; \forall k \in \intn{1}{n}$.
\end{itemize}
Then the leading asymptotics of $\ddet{2}{I+K}$ read:
\bem
\ddet{2}{I+K}=G^x_2\pac{\sg} \paf{x}{2}^{\sul{p=1}{n} \ga^2_p-\de_p^2} E\pac{F} \pl{k=1}{n} \f{G\pa{1+\de_k-\ga_k}G\pa{1-\de_k-\ga_k} }{ G\pa{1-2\ga_k} }
\\
\pl{k=1}{n} \f{F^{\ov{\nu}_k}_{+}\pa{a_k}}{F^{\ov{\nu}_k}_+\pa{a_k+i}}
\f{F^{\nu_k}_{-}\pa{a_k}}{F^{\nu_k}_-\pa{a_k-i}}
\pl{k\not=p}{} \paf{\pa{a_k-a_p+i}^2 }{ \pa{a_{k}-a_p+2i}\pa{a_k-a_p}  }^{\ov{\nu}_k \nu_p}
\pa{1+\e{O}\pa{x^{\rho-1}}} \; .
\label{LA formule de l'article}
\end{multline}
We have defined
\beqa
G^x_2\pac{\sg}&=&\exp\paa{x\Int{\R}{} \f{\dd \xi}{2\pi} \pac{\ln\pa{\sg}\pa{\xi}+1-\sg\pa{\xi}  }} \; , \\
E\pac{F} &=& \exp\paa{\Int{0}{+\infty} \dd \xi  \; \xi \,   \mc{F}^{-1}\pac{\ln F}\pa{\xi} \, \mc{F}^{-1}\pac{\ln F}\pa{-\xi}  } \; .
\eeqa
We also remind that $\rho=2\max\abs{\Re\pa{\de_k}}$.
\end{theorem}
\Proof
The result follows by a recursive applications of propositions \ref{proposition derivee delta} and \ref{proposition derivee gamma}. At the end of the
recursion, one also needs to invoke the Aheizer-Kac formula for the $2$-determinant of the truncated Wiener-Hopf operator generated by the
$L^{2}\pa{\R}$ symbol $F-1$ so as to fix the constant $E\pac{F}$ and the $F$ dependent part of $G_2^x\pac{\sg}$.   $\Box$

\vspace{3mm}
The leading asymptotics of Wiener-Hopf operators generated by general Fisher-Hartwig symbols \eqref{LA formule de l'article} reproduces all the
previously know results:  $ \forall \; k\;  \nu_k=0  \;$, or  $ \;\forall \; k\; \ov{\nu}_k=0 \; $ proven in
\cite{BottcherSilbWienerHopfAsympAllNuorBarNuZero},  $\forall \; k \; \ga_k=0\;$  proven in \cite{BottcherSilberWidomProofFHConjWienerHopfDiscontinuity}
 and also the case of a pure Fisher-Hartwig singularity $n=1$ and $F=1$\cite{BasorWidomWienreHopfwithOneFHSingularity}.
We refer the reader to \cite{BottcherSilAnalysisOfToeplitzOperators} for a restatement of all the know results in a language very close to the one used
in this article.
One should only pay attention to the different definition of $\sg_{\nu_k,\ov{\nu}_k}$ between this article and the book
\cite{BottcherSilAnalysisOfToeplitzOperators}. Indeed $\nu_k$, resp. $\ov{\nu}_k$, differ by an overall minus sign
with respect to the conventions of the latter book.

However, our result disproves the continuous analog of the Fisher-Hartwig conjecture \cite{BottcherFHConjectureWienerHopf} in its broad generality.
Although most of the factors between the formula and the conjecture coincide, the latter predicts the presence of
\beq
\pl{k<p}{} \paf{\pa{a_k-a_p+1}^2 }{ \pa{\pa{a_k-a_p}^2+4}\pa{a_k-a_p}^2  }^{\ov{\nu}_k \nu_p} \;\; ,
\enq
\noindent whereas we find
\beq
\pl{k<p}{} \paf{\pa{a_k-a_p+i}^2 }{ \pa{a_k-a_p+2i}\pa{a_k-a_p}  }^{\ov{\nu}_k \nu_p}
\paf{\pa{a_k-a_p-i}^2 }{ \pa{a_k-a_p-2i}\pa{a_k-a_p}  }^{\ov{\nu}_p \nu_k} \;\; .
\enq
\noindent The difference comes from the presence of  $\ov{\nu}_p\nu_k$ in the second exponent instead of $\ov{\nu}_k\nu_p$. Of course in all the cases
previously investigated, the difference  between the conjecture and the
present result was not appearing as either the factor was not present or $\nu_p$ and $\ov{\nu}_p$ were related by a sign.

To end this Section we would like to stress that it is not a problem to obtain the sub-leading asymptotics of truncated Wiener-Hopf with
 Fisher-Hartwig symbols by the so-called $x$-derivative method:
\beq
\Dp{x} \ln \ddet{}{I+V}= -\f{i}{2} \e{tr}\paa{\chi_1 \sg_3 } \quad \e{with} \quad \chi\pa{z}=I_2+\f{\chi_1}{z}+\e{o}\pa{z^{-1}} \;\;\;
\e{for} \; z\tend \infty \; .
\nonumber
\enq
The above is a straightforward generalization of the identity for the pure sine kernel given in \cite{DeiftItsZhouSineKernelOnUnionOfIntervals}.
 The $x\tend+\infty$ asymptotics for $\chi_1$ can be obtained by solving perturbatively the singular integral equation satisfied by $\Omega$. We
  do not present the calculations here as we are going to derive the sub-leading asymptotics for the Toeplitz determinant case investigated in the next Section.

\section{Toeplitz matrices with Fisher-Hartwig type symbols}
In this Section, we adapt the previous analysis of the generalized sine kernel. In this way we obtain, in the framework of Riemann-Hilbert problems, the asymptotic behavior of Toeplitz determinants generated by symbols $\sg$ having Fisher-Hartwig singularities.
This approach is based on an observation made  by Deift, Its and Zhou in
\cite{DeiftItsZhouSineKernelOnUnionOfIntervals} concerning the relationship between
the Fredholm determinant of a sine kernel on a circle and a Toeplitz determinant.
Our results reproduce those obtained by T.Ehrhardt in his thesis \cite{EhrhardtAsymptoticBehaviorOfFischerHartwigToeplitzGeneralCase}.
Moreover, the Riemann-Hilbert approach allows to compute sub-leading asymptotics to any order
in a quite systematic, although quickly cumbersome way. At the end of this Section we shall
establish the first sub-leading asymptotics of Toeplitz matrices with Fisher-Hartwig singularities. We observe that these sub-leading asymptotics of
Toeplitz determinants partly restore the
independence on the choice of a Fisher-Hartwig type representant for the symbol $\sg$. Indeed the jumps of $\sg$ are
characterized by parameters $\de_k$. A shift of any $\de_k$ by an integer describes the same jump. The freedom of choice of a Fisher-Hartwig
representant for $\sg$ is broken if one considers the leading asymptotics only. These correspond to the choice $\abs{ \Re\pa{\de} }<\tf{1}{2}$. However,
part of the sub-leading asymptotics (the so-called oscillating ones) we obtain correspond to shifts $\de_i\mapsto \de_i+1\, , \, \de_j\mapsto \de_j-1$ in the parameters appearing in the
leading asymptotics.
These sub-leading asymptotics shed a light on the mechanism appearing in the asymptotics for ambiguous case type symbols  that has been conjectured by Basor and Tracy
\cite{BasorTracyGeneralizedFischer-HartwigConjecture} and proven recently by Deift, Its and Krasovsky \cite{DeiftItsKrasovskyAsymptoticsofToeplitsHankelWithFHSymbols}. Indeed, terms that were subdominant in the asymptotic series for a generic set of parameters become of the same order of magnitude as the leading asymptotics
when some of the parameters $\de_p$ and $\ga_p$ are set to these specific ambiguous values.
The global structure of the sub-leading asymptotics seems to follow the scheme already pointed out in \cite{KozKitMailSlaTerRHPapproachtoSuperSineKernel}. We
formulate a conjecture on this global structure at the end of this Section. Our conjecture can be seen as a generalization of the Basor-Tracy conjecture: we believe that the full asymptotic series for $\ddet{m}{T\pac{\sg}}$  results of a $1-\de$ periodization of
only a small part of the asymptotic series.

\subsection{The Riemann-Hilbert Problem}
Let us consider an integral operator acting on the unit circle $\msc{C}$ with the kernel
\beq
V\pa{z,z'}=\sqrt{\sg\pa{z}-1}\sqrt{\sg\pa{z'}-1}
\f{z^{\f{m}{2}}{z'}^{-\f{m}{2} } -z^{-\f{m}{2}}{z'}^{\f{m}{2} }  }{2i\pi \pa{z-z'}} \; ,
\label{noyau sinus sur le cercle}
\enq
where
\beq
\sg\pa{z}= b\pa{z} \pl{k=1}{n}\pa{1-\f{z}{a_k}}^{-\nu_k} \pa{1-\f{a_k}{z}}^{-\ov{\nu}_k}
\; ,\;\;  a_k \in \msc{C} \; .
\label{symbole sg FH pour Toeplitz}
\enq
\noindent There we assume that $b$ is holomorphic and non-vanishing in a vicinity of $\msc{C}$ and has zero winding number.
One can actually characterize the singular behavior of $\sg$ on the contour
of integration more explicitly. Namely,
\beq
\pa{1-\f{z}{a_k}}^{-\nu_k} \!\!\!\pa{1-\f{a_k}{z}}^{-\ov{\nu}_k}=
\f{\ex{i\de_k\pa{\th-\th_k-\pi\e{sgn}\pa{\th-\th_k} } } }{\pa{2-2\cos\pa{\th-\th_k}}^{\ga_k} }\; ,
\quad \f{z}{a_k}=\ex{i\pa{\th-\th_k}} \; ,\;  \th-\th_k \in \intoo{-\pi}{\pi} .
\enq
\noindent We have set, just as in the preceding sections,
\beq
 2 \de_k=\ov{\nu}_k-\nu_k \;\;\; 2 \ga_k=\ov{\nu}_k+\nu_k \; .
\enq
Here, we assume that $\abs{\Re\pa{\de_k}}<\tf{1}{2}$ and $\Re\pa{\ga_k}<\tf{1}{2}$.

The authors of \cite{DeiftItsZhouSineKernelOnUnionOfIntervals} observed that for $m\in \mathbb{N}$ the pure sine kernel on the unit circle is of finite rank. This
property persists in the case of the integral operator under investigation as:
\beq
V\pa{z,z'}=\sqrt{\sg\pa{z}-1} \sqrt{\sg\pa{z'}-1} \sul{p=1}{m} \f{z^{p-1-\f{m}{2}} {z'}^{\f{m}{2}-p}}{2i\pi} \; .
\enq
Hence we have that
\beqa
\ddet{\msc{C}}{I+V}&=&\ddet{m}{\de_{jk}+ \Int{0}{2\pi} \f{\dd \th}{2\pi} \pa{\sg\pa{\ex{i\th}}-1} \ex{i\th \pa{k-j}} } \\
&=&\ddet{m}{\Int{0}{2\pi} \f{\dd \th}{2\pi} \sg\pa{\ex{i\th}} \ex{i\th \pa{k-j}} } \; .
\eeqa
\noindent We used the subscript $\msc{C}$ in order to insist that the $lhs$ is the Fredholm determinant of an integral operator acting on the unit circle $\msc{C}$ whereas the $rhs$ is the determinant of an $m \times m$ matrix.

Hence, the
asymptotics of Toeplitz matrices with symbols having jump and power-law singularities
will follow from those of the Fredholm determinant of the integral operator defined in \eqref{noyau sinus sur le cercle}. The only significant difference between the kernel \eqref{noyau sinus sur le cercle} and the one considered
in the preceding Sections is the interval on which they act. Most of the steps in the derivation of the asymptotics are very similar. We only insist on the most striking differences.

\subsection{Asymptoic solution of the Riemann-Hilbert Problem}
We consider the RHP for a piecewise analytic matrix $\chi$ having a jump on the unit circle $\msc{C}$:
\begin{itemize}
\item $\chi$ is analytic on   $\mathbb{C}\setminus\msc{C} \; $;
\item $\forall  k \in \intn{1}{n} \;,$ there exists $M_k\in \e{GL}_2\pa{\Cx}$ such that
\vspace{2mm} \newline$ \chi = M_k \paa{I_2+g\pa{z} B\pa{z}   +\abs{z-a_k}\pa{g\pa{z}+1}  \e{O}\pa{\ba{cc}1 & 1 \\ 1 & 1 \ea }}
\;\; , \; z \tend a_k \;$ \; ;
\item $\chi \underset{z \tend \infty}{\tend} I_2 \; $;
\item $\chi_{+}\pa{z} G\pa{z}=\chi_-\pa{z}  \;\; ; \quad z \in \msc{C} \, .$
\end{itemize}
There, just as for the Wiener-Hopf case, the rank one matrices $B_k$ read
\beq
B_k = \pa{\ba{cc} -1 & z^m \\ - z^{-m} & 1 \ea } \;.
\enq
The function $g$ is also defined similarly
\beq
g\pa{z} = \Int{\msc{C}}{} \f{\dd s }{2i\pi} \f{\sg\pa{s}-1}{z-s} \; .
\enq
It has a singular behavior at $z=a_k$ of the type
\beq
 g\pa{z}=\left\{ \ba{cc c}  \e{O}\pa{1}+\e{O}\pa{\pa{z-a_k}^{-2 \ga_k}}  &  \e{for} & \ga_k\not=0 \\
                          \e{O}\pa{\ln \pa{z-a_k}} & \e{for} & \ga_k=0  \ea \right.
\;  \e{when} \; z\tend a_k \; .
\enq
We finally precise that the unit circle $\msc{C}$ is oriented canonically (\textit{ie} the $+$ side of the
contour corresponds to the interior of the circle) and that the  jump matrix $G$ reads
\beq
 G\pa{z}
            = \left(\ba{cc}
                   2-\sg\pa{z} & \pa{\sg\pa{z}-1} z^m \\
                   \pa{1-\sg\pa{z}} z^{-m} & \sg\pa{z}
                   \ea \right) \;\; .
\enq
We now define a new matrix $\Upsilon$ according to
\begin{itemize}
\item $\Upsilon=\chi \, \a^{\sg_3}\, , $ for $z$ being in the exterior of $\Ga_-$ and the interior of $\Ga_+\; ;$
\item $\Upsilon=\chi \, \a^{\sg_3}_{\da} M_{\da}^{-1}\, , $ for $z$  between $\Ga_-$ and $\msc{C}\; ;$
\item $\Upsilon=\chi \, \a^{\sg_3}_{\ua} M_{\ua}\, , $ for $z$ between $\Ga_+$ and $\msc{C}$.
\end{itemize}
Here $\a$ is the solution of the scalar RHP
\beq
\a\;  \e{analytic}\; \e{on} \; \Cx\setminus \msc{C} \quad \a_-=\sg \a_+  \; ,\qquad z \in \msc{C}\setminus\cup_{k=1}^{n}\paa{a_k}
\quad \a\tend 1 \; \e{when}\; z \tend \infty \;.
\enq
This scalar RHP can be solved explicitly in terms of the canonical Wiener-Hopf factors of $b$:
$b=b_+G\pac{b}b_-$. One has $\a=\a_{\ua}$ on $D_{0,1}$ and $\a=\a_{\da}$ on $\Cx\setminus\ov{D}_{0,1}$, with
\beqa
\a_{\ua}\pa{z}&=& b_+^{-1}\pa{z}G\pac{b}^{-1} \pl{k=1}{n} \pa{1-\f{z}{a_k}}^{\nu_k} \\
\a_{\da}\pa{z}&=& b_-\pa{z}\pl{k=1}{n} \pa{1-\f{a_k}{z}}^{-\ov{\nu}_k} \;.
\eeqa
\noindent   The matrices
 $M_{\ua/\da}$ defining $\Upsilon$ read
\beq
 M_{\ua}\pa{z}= \left(\ba{cc}
                         1& \pa{1-\sg^{-1}}\a_{\ua}^{-2}\pa{z} z^m \\
                         0& 1 \ea\!\! \right), \quad
 M_{\da}\pa{z}= \left(\ba{cc}
                        1   &    0\\
                         \pa{\sg^{-1}-1}\a_{\da}^{2}\pa{z}z^{-m}& 1 \ea \!\!\right)\; .
\enq
We stress that, just as for the Wiener-Hopf case, the matrix $M_{\ua/\da}\pa{r\ex{i\th}}$ should
be understood as the analytic continuation of $M_{\ua/\da}\pa{\ex{i\th}}$ from a small neighborhood of $\ex{i\th}$
to the ray $\intff{\ex{i\th}}{r\ex{i\th}}$.
\begin{figure}[h]
\centering
\includegraphics{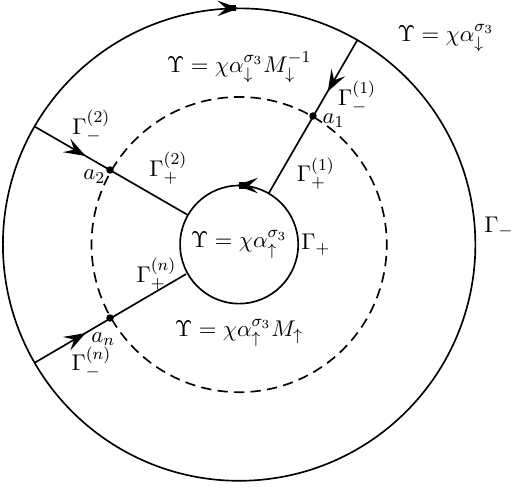}
\caption{ Contour for the RHP $\Upsilon$ and the associated contour $\Ga$. \label{Contour du RHP pour Y Toeplitz}}
\end{figure}
One readily sees that $\Upsilon$ satisfies the RHP
\begin{itemize}
\item $\Upsilon$ is analytic in $\Cx\setminus \Ga\; ;$
\item  $\forall  k \in \intn{1}{n} \;,$ there exists $M_k\in \e{GL}_2\pa{\Cx}$, such that:
\vspace{2mm} \newline
$\hspace{-3mm}\Upsilon\pa{z}  = M_k \paa{I_2+g\pa{z} B\pa{z}
+\abs{z-a_k}\pa{g\pa{z}+1}  \e{O}\pa{\ba{cc}1 & 1 \\ 1 & 1 \ea } }
M\pa{z} \;\, \;z \tend a_k$
\item $\Upsilon \limit{z}{\infty} I_2\; ;$
\item $\left\{ \ba{c} \Upsilon_{+}\pa{z} M_{\ua}\pa{z}=\Upsilon_-\pa{z}  ; \quad z \in \Gamma_+  \\
                      \Upsilon_{+}\pa{z} M_{\da}^{-1}\pa{z}=\Upsilon_-\pa{z}  ; \quad z \in \Gamma_-
                \ea \right. \; , $
\item $\left\{ \ba{c} \Upsilon_{+}\pa{z} N^{\pa{l}}\pa{z}=\Upsilon_-\pa{z}  ; \quad z \in \Gamma^{(l)}_+  \\
                      \Upsilon_{+}\pa{z} \ov{N}^{\pa{l}}\pa{z}=\Upsilon_-\pa{z}  ; \quad z \in \Gamma^{(l)}_-
                \ea \right.  \quad, \;\;  l \in \intn{1}{n}\;\;.$
\end{itemize}
The function $g$ and the rank one matrices $B\pa{z}$ are as defined above. Moreover, the matrix $M$ reads
\beq
M\pa{z} = \left\{ \ba{cc} \a^{\sg_3}_{\ua}M_{\ua} & z \in \paa{D_{0,1}\setminus \cup_{k=1}^{n}\Ga_+^{\pa{k}}}\cap U \vspace{2mm}\\
\a^{\sg_3}_{\da}M_{\da}^{-1} & z \in \paa{\Cx\setminus \pa{D_{0,1} \cup_{k=1}^{n}\Ga_+^{\pa{k}}}}\cap U  \ea \right. \; .
\enq
The local behavior at $z\tend a_k$ of $M$ can be inferred from the one of $\a_{\ua/\da}$ and the explicit formulae for $M_{\ua/\da}$.
\noindent The jump matrices $N^{\pa{l}}\pa{z}$, $\ov{N}^{\pa{l}}\pa{z}$ are defined as
\beqa
N^{\pa{l}}\pa{z}&=& \left( \ba{cc}
                                1 & n_l\pa{z}\paf{z}{a_l}^m \\
                                0 & 1 \ea \right)
 \qquad  z\in \Gamma_+^{\pa{l}}  \;\; ,
 \nonumber\\
\ov{N}^{\pa{l}}\pa{z}&=& \left( \ba{cc}
                                1 & 0 \\
                                \ov{n}_l\pa{z}\paf{a_l}{z}^{m} & 1 \ea \right)
\qquad z\in \Gamma_-^{\pa{l}}  \;\;  \nonumber ;
\eeqa
\noindent and their entries read
\beqa
n_l\pa{z} &=& \f{ \pac{-im\ln \pa{\tf{z}{a_l}} }^{2\de_l } }{K_l\pa{z} } \pa{\ex{-2i\pi \ov{\nu}_l}-1} \;\; ,\\
\ov{n}_l\pa{z}&=& \f{ K_l\pa{z}}{\pac{-im\ln \pa{\tf{z}{a_l}} }^{2\de_l } } \pa{\ex{2i\pi\nu_l}-1} \;\; .
\eeqa
We have defined, analogously to the Wiener-Hopf case,
\beq
K_l\pa{z}= \f{b_-\pa{z}}{b_+\pa{z}}\ex{-i\pi\ga_l}\paa{\f{1-\tf{z}{a_l}}{-\ln \pa{\tf{z}{a_l}}}}^{-2\de_l}\f{m^{2\de_l}}{a_l^m} \paf{z}{a_l}^{\ov{\nu}_l} \pl{\substack{r=1 \\ \not= l}}{n} \pa{1-\f{z}{a_r} }^{\nu_r} \pa{1-\f{a_r}{z}}^{-\ov{\nu}_r} \; .
\enq
It is not a problem to see that the parametrix around  $a_p$ can be chosen as
\beq
P_{a_p}\!\pa{ z} = \!\pa{\!\! \ba{cc}
                        \Psi\pa{\ga_p-\de_p;-i\zeta_p}
                                    &ib^{\pa{p}}_{12}\pa{z} \Psi\pa{1+\ga_p+\de_p;i\zeta_p} \\
                        -ib^{\pa{p}}_{21}\pa{z} \Psi\pa{1+\ga_p-\de_p;-i\zeta_p}
                                        & \Psi\pa{\ga_p+\de_p;i\zeta_p} \ea \!\!}
\f{L_p }{\pa{\zeta_p}^{\de_p\sg_3-\ga_p} } \; .
\enq
Where  we have set $\zeta_p=-im\ln\pa{\tf{z}{a_p}}$ and the second parameter of the CHF's is implicitly assumed to be
$1+2\ga_p$. The piecewise constant matrix $L_p$ reads
\beq
L_p=\left\{\ba{cc}  \ex{i\f{\pi \de_p}{2}} \ex{-i\f{\pi \ga_p}{2}\sg_3}
                & -\tf{\pi}{2}<\e{arg}\pa{\zeta_p} <\tf{\pi}{2} \; \vspace{4mm}\\
                 \ex{i\f{\pi \pa{\de_p-\ga_p}}{2}}\left( \ba{cc} 1  &0 \\
                                                0& \ex{-2i\pi \de_p-i\pi\ga_p}
                                                   \ea \right)
                 &\tf{\pi}{2}<\e{arg}\pa{\zeta_p}<\pi \vspace{4mm}\\
                  \ex{i\f{\pi \pa{\de_p+\ga_p}}{2}}\left( \ba{cc} \ex{-2i\pi \de_p+i\pi \ga_p}  &0 \\
                                     0& 1
                                    \ea \right)  & -\pi<\e{arg}\pa{\zeta_p}< -\tf{\pi}{2}
        \ea \right. \;\; ,
\enq
\noindent and the coefficients $b^{\pa{p}}_{12}\pa{z}$ and $b^{\pa{p}}_{21}\pa{z}$ are given by
\beq
b^{\pa{p}}_{12}\pa{z}= \f{i \ex{-i\pi\ga_p}}{K_p\pa{z}} \f{\Gamma\pa{1-\ga_p+\de_p}}{\Gamma\pa{-\ga_p-\de_p}}  \qquad ,\qquad
b^{\pa{p}}_{21}\pa{z}=  -i K_p\pa{z} \ex{i\pi\ga_p} \f{\Gamma\pa{1-\ga_p-\de_p}}{\Gamma\pa{\de_p-\ga_p}} \; .
\label{fonction b parametrix RHP cercle}
\enq

The matrix
\beq
\Omega=\left\{  \ba{cc}
                    \Upsilon P_{a_k}^{-1} & z\in D_{a_k,\eps}\\
                    \Upsilon & z \in \mathbb{C}\setminus\bigcup\limits_{p=1}^{n}\ov{D}_{a_p,\eps}
            \ea \right.
\label{matrice Omega pour RHP cercle}
\enq
\noindent  satisfies the RHP
\begin{itemize}
\item $\Omega$ is analytic in $\mathbb{C}\setminus \Sg_{\Omega} \;\; ,\quad \Sg_{\Omega}=\pa{\Ga\setminus\ga}\bigcup\limits_{p=1}^{n} \Dp{}D_{a_p,\eps} \;\; ,\;  \e{with}\; \ga=\Ga\cap\pa{\bigcup\limits_{p=1}^{n}D_{a_p,\eps} } \; ;$
\item $\Omega = I_2+ O\pa{\tf{1}{z}} \quad;\quad   z \tend \infty \; ,$
\item $\left\{ \ba{rc c l c c c} \Omega_{+}\pa{z} M_{\ua}\pa{z}&=\Omega_-\pa{z}  &,& \; z \in \Gamma_+  \; ;\\
                      \Omega_{+}\pa{z} M_{\da}^{-1}\pa{z}&=\Omega_-\pa{z}  &,& \;z \in \Gamma_- \; ; \\
                      \Omega_{+}\pa{z} N^{\pa{l}}\pa{z}&=\Omega_-\pa{z}  &, & \; z \in \widetilde{\Gamma}^{\pa{l}}_+
                       = \Ga_{+}^{\pa{l}}\setminus\paa{\Ga_+^{\pa{l}}\cap D_{a_l,\eps}}\; ; \\
                      \Omega_{+}\pa{z} \ov{N}^{\pa{l}}\pa{z}&=\Omega_-\pa{z}  &,& \; z \in
                      \widetilde{\Gamma}^{\pa{l}}_-
                        = \Ga_{-}^{\pa{l}}\setminus\paa{\Ga_-^{\pa{l}}\cap D_{a_l,\eps}}\; ;\\
                      \Omega_{+}\pa{z} P_{a_l}\pa{z}&=\Omega_-\pa{z}  &, &\; z \in \partial D_{a_l,\eps}\; .  \\
                \ea \right.$
\end{itemize}
The solution of the RHP for $\Omega$, clearly exits and is unique. Moreover
it is uniformly $I_2+\e{O}\pa{\tf{1}{x^{1-\rho}}}$ with $\rho=2\, \underset{k}{\max}{\abs{\Re \pa{\de_k} }}$.
\subsection{The resolvent and asymptotics of $\ln\ddet{}{I+V}$}

Just as in the case of kernels acting on $\R$, one can reconstruct the approximate resolvent $R_0$ in terms of Humbert's CHFs:
\beq
R_0\pa{z,z}=\f{1-\sg^{-1}\pa{z}}{2i\pi} \paa{\f{m}{z} - \Dp{z} \ln \pa{\a_+\a_-} \pa{z} } \; ,
\enq
for $z \in \msc{C}\setminus \pa{\bigcup\limits_{k=1}^{n}\msc{C}\cap D_{a_k,\eps}  }$ and
\bem
R_0\pa{z,z}= \ex{-i\pi \de_p} \f{\pa{\sg\pa{z}-1} \wh{\a}_{+}^{\pa{p}}\pa{z}}{2i\pi m^{2\ga_p} \wh{\a}_-^{\pa{p}}\pa{z} }
\left\{-\f{m}{z} \tau\pa{\ga_p,\de_p;-im \ln \tf{z}{a_p}}  \right.  \nonumber  \\
 \quad \left. - \Dp{z} \ln \pa{\wh{\a}^{\pa{p}}_+\wh{\a}^{\pa{p}}_-} \pa{z} \varphi\pa{\ga_p,\de_p;-im \ln \tf{z}{a_p}} \right\} +\e{o}\pa{1} \;\;,
\end{multline}
\noindent whenever $z \in \bigcup\limits_{k=1}^{n}\msc{C}\cap D_{a_k,\eps}$.
The functions $\tau\pa{\ga,\de;z}$ and $\varphi\pa{\ga,\de;z}$ have been defined in \eqref{definition de tau} and \eqref{definition de varphi}.

One should observe that the proof of the differential identities \eqref{derivee partielle trace class} is contour independent, so that we can also use them for a generalized sine kernel acting on $\msc{C}$. Moreover one can prove that, just as for the Wiener-Hopf kernel, if
$R$ denotes the exact resolvent of the kernel \eqref{noyau sinus sur le cercle},
\beq
\Int{\msc{C}}{} \dd z \f{R\pa{z,z}-R_0\pa{z,z}}{\sg-1}\Dp{\beta_p}\sg =\e{o}\pa{1} \;\; .
\enq
The $\e{o}\pa{1}$ terms vanish in the ordered limit $x \tend + \infty $ and then $\eps \tend 0^+$.   Harping on the steps for the integration in the case
of a kernel acting on $\R$ we get
\begin{prop}
\label{proposition derivee delta et gamma Toeplitz}
Up to $\e{o}\pa{1}$ terms in the $x\tend +\infty$ limit, one has
\bem
\ln\paf{\ddet{}{1+V}}{\ddet{}{I+V}_{\mid_{\de_1=0}}}=m \Int{\msc{C}}{} \f{\dd z}{ 2i\pi z} \ln \paf{\sg}{\sg_{\mid_{\de_1=0}} }
-\de_1^2 \ln m  \\
+ \de_1 \ln \paf{\a_{\ua}^{\pa{2\dots n}}\pa{0}}{ \a_{\ua,+}^{\pa{2\dots n}}\pa{a_1}\a_{\da,-}^{\pa{2\dots n}}\pa{a_1} } + \ln \paf{G\pa{1-\ga_1+\de_1}G\pa{1-\ga_1-\de_1}}{G\pa{1-\ga_1}G\pa{1-\ga_1}} + \e{o}\pa{1} \;\; .
\end{multline}
\bem
\ln\paf{\ddet{}{1+V}_{\mid_{\de_1=0}}}{\ddet{}{I+V}_{\mid_{\de_1=\ga_1=0}} }=m \Int{\msc{C}}{} \f{\dd z}{2i\pi z}
        \ln \paf{\sg_{\mid_{\de_1=0}} }{\sg_{\mid_{\de_1=\ga_1=0}} }
+\ga_1^2 \ln m \\
- \ga_1 \ln \paf{\a_{\ua,+}^{\pa{2\dots n}}\pa{a_1}}{ \a_{\ua,+}^{\pa{2\dots n}}\pa{0}\a_{\da,-}^{\pa{2\dots n}}\pa{a_1} }   + \ln \paf{G\pa{1-\ga_1}G\pa{1-\ga_1}}{G\pa{1-2\ga_1}} + \e{o}\pa{1} \;\; .
\end{multline}
There,
\beq
\a_{\ua}^{\pa{2\dots n}}\pa{z}= b_+^{-1}\pa{z}G\pac{b}^{-1} \pl{k=2}{n} \pa{1-\f{z}{a_k}}^{\nu_k} \qquad
\a_{\da}^{\pa{2\dots n}}\pa{z}= b_-\pa{z}\pl{k=2}{n} \pa{1-\f{a_k}{z}}^{-\ov{\nu}_k} \;.
\enq
Also, $\a_{\ua/\da, +/-}$ stand for their boundary values from the $\pm$ sides of $\msc{C}$.
\end{prop}
\Proof
The proof goes exactly the same as in the case of Wiener-Hopf operators. The only notable difference is that one needs an additional version of the integration lemma.

Let $\eps>0$, and $I$ be a sub-interval of $\R$ such that $\forall \overset{\circ}{I}\supset \intff{-\eps}{\eps}$.
Then $\forall \, h \in \msc{C}^1\pa{I,\Cx}$, and for any
Riemann integrable on $\R$ function $ \mc{R}\pa{t}$  that it is decreasing as a power-law when $\abs{t}\tend \infty$, one has
\beq
\Int{-\eps}{\eps} m \dd t \mc{R}\pa{m t }\, h\pa{t}\, \ln\pac{2 \pa{1-\cos t }}  =  h\pa{0}  \Int{\R}{}\ln \pac{t^2} \mc{R}\pa{t}   \dd t  \; + \; \e{o}\pa{1} \; .
\enq
\noindent The $\e{o}\pa{1}$ term is vanishing in the ordered limit  $m\tend +\infty$  and then  $\eps \tend 0$. The proof is straightforward by applying the original integration Lemma \ref{lemme preparatoire} to the function
$h\pa{t}\ln\pac{2 \pa{1-\cos t }}-h\pa{0}\ln \pac{t^2}$ that is $\msc{C}^{1}\pa{I,\Cx}$. Moreover the power-law decrease of $\mc{R}\pa{t}$ at infinity ensures that the contributions of the boundary which are of the type $\ln m \Int{ m \eps}{+\infty} \dd t \mc{R}\pa{t}$
will indeed be subdominant with respect to $\e{o}\pa{1}. \quad \Box$

We point out that
\beq
\Int{\msc{C}}{} \f{\dd z}{2i\pi z} \ln \paf{\sg_{\mid_{\de_1=0}} }{\sg_{\mid_{\de_1=\ga_1=0}} }=0=
\Int{\msc{C}}{} \f{\dd z}{2i\pi z} \ln \paf{\sg}{\sg_{\mid_{\de_1=0}} } \; .
\enq
We chose to include these factors in the proposition above so as to make the parallel with the Wiener-Hopf case more obvious. By repeatedly applying Proposition \ref{proposition derivee delta et gamma Toeplitz} and then invoking the
strong Szeg\#{o} limit theorem we get
\begin{theorem}
Let $T\pac{\sg}$ be an $m\times m$ Toeplitz matrix generated by the symbol
\beq
\sg\pa{\theta}=b\pa{z} \pl{k=1}{n} \pa{1-\f{z}{a_k}}^{-\nu_k}\pa{1-\f{a_k}{z}}^{-\ov{\nu}_k} \quad a_k \in \msc{C} \, , \; z=\ex{i\th}
\enq
where $b$ is analytic and non-vanishing in some open neighborhood of the unit circle $\msc{C}$ and $2\de_k=\ov{\nu}_k-\nu_k$, $2\ga_k=\ov{\nu}_k+\nu_k$ are such that $\abs{\Re\pa{\de_k}}<\tf{1}{2}$, $\Re\pa{\ga_k}<\tf{1}{2}$.

Then the leading asymptotics of $\ddet{m}{T\pac{\sg}}$ are given by
\bem
\ddet{m}{T\pac{\sg}}=\pa{G\pac{b}}^m \,E\pac{b} \, m^{\sul{p=1}{n}\ga_p^2-\de_p^2}
\pl{p=1}{n} \f{G\pa{1-\ga_p+\de_p}G\pa{1-\ga_p-\de_p}}{G\pa{1-2\ga_p}} \\
\pl{p=1}{n} b_{+}^{\ga_p+\de_p}\pa{a_p} b_{-}^{\ga_p-\de_p}\pa{a_p}
\pl{p\not=q}{n} \pa{1-\f{a_p}{a_q}}^{\pa{\ga_p+\de_p}\pa{\de_q-\ga_q}} \pa{1+\e{o}\pa{1}}
\label{asymptotiqueToeplitzFH}
\end{multline}
where,
\beq
G\pac{b}=\ex{ \pac{\ln b}_o} \qquad E\pac{b}=\ex{\sul{k}{+\infty} k \pac{\ln b}_k\pac{\ln b}_{-k}} \qquad
\pac{\ln b}_k=\Int{0}{2\pi} \f{\dd \th}{2\pi} \ex{i k \th} \ln b\pa{\th} \;\; .
\enq
\end{theorem}

This reproduces the result of T.Ehrhardt \cite{EhrhardtAsymptoticBehaviorOfFischerHartwigToeplitzGeneralCase}. One should pay attention that the
exponents $\de_i$ and $\ga_i$ correspond to the exponents $\beta_i$ and $-\a_i$ in Ehrhardt's notations.

\subsection{The sub-leading asymptotics}

The authors of \cite{DeiftItsZhouSineKernelOnUnionOfIntervals} found a way to express the discrete derivative of $\ln T_m\pac{\sg}$
in terms of the RHP data. Their result reads
\beq
\f{T_{m+1}\pac{\sg}}{T_{m}\pac{\sg} }= \chi_{11}\pa{z=0} \;\; .
\enq
Where $\chi_{11}\pa{0}$ stands for the upper diagonal entry of the solution to the RHP for $\chi$ given at the beginning of this Section. Their proof also works, word for word, in the case of the more complicated kernel we consider, so we omit it here.  It is now enough to determine the sub-leading asymptotics of the matrix $\Omega$ defined by \eqref{matrice Omega pour RHP cercle} thanks to the singular integral equation:
\beq
\Omega\pa{z}=I_2+\Int{\Sg_{\Omega}}{} \f{\dd s}{2i\pi\pa{z-s}} \Omega_+\pa{s} \De\pa{s}
\label{eqn singuliere integrale pour Omega RHP cercle}
\enq
with $I_2+\De$ being the jump matrix for $\Omega$. The method for computing the corrections is standard. We send the reader to \cite{KozKitMailSlaTerRHPapproachtoSuperSineKernel} for more details. We stress that we did not chose this reference because of its originality with respect to the perturbation theory of such integral equations. We rather chose it as there, the perturbation theory is applied in notation quite close to the ones we use.

It is easy to see that the jump matrix $\De$ is exponentially vanishing with respect to $x$ away of the disks $\Dp{}D_{a_p,\eps}$. However for
 $s \in \Dp{}D_{a_p;\eps}$ it admits the asymptotic expansion
\beq
\De\pa{s}=\sul{\ell=1}{M}\f{1}{\ell! m^{\ell} \pac{\ln \pa{\tf{s}{a_p}} }^{\ell}} \De^{\pa{p}}_{\ell}\pa{s} + \e{O}\pa{m^{-M+\rho-1}}\;\;. \;\;
\enq
We have set, using $\pa{a}_n=\tf{\Ga\pa{a+n}}{\Ga\pa{a}}$,
\bem
\De^{\pa{p}}_{\ell}\pa{s}=\pa{\ba{cc} 1 & \f{i \ell b_{12}^{\pa{p}}\pa{s} }{\de_p^2-\ga_p^2}  \\
                        -\f{i \ell b_{21}^{\pa{p}}\pa{s} }{\de_p^2-\ga_p^2} & 1 \ea }  \\
\times   \pa{ \ba{cc} \pa{\ga_p-\de_p}_{\ell}\pa{-\ga_p-\de_p}_{\ell} & 0 \\
                    0 &  \pa{-1}^{\ell}\pa{\ga_p+\de_p}_{\ell}\pa{\de_p-\ga_p}_{\ell} \ea } \; .
\end{multline}
and the functions $b_{12}^{\pa{p}}\pa{s}$ and $b_{21}^{\pa{p}}\pa{s}$ defined in \eqref{fonction b parametrix RHP cercle} exhibit a slight dependence on $m$ that is a $\e{O}\pa{m^{\rho}}$. The standard considerations of a perturbative approach to
\eqref{eqn singuliere integrale pour Omega RHP cercle} lead to
\beq
\Omega\pa{0}=I_2+\f{\Omega_1\pa{0}}{m} + \f{\Omega_2\pa{0}}{m^2}+ \e{O}\pa{\f{1}{m^3}, \f{ {\pa{\tf{a_i}{a_j}}}^m  m^{3\rho}}{m^3}} \;\; .
\enq
Where, for $z \in \Cx \setminus \bigcup\limits_{p=1}^{n} \partial D_{a_p,\eps}$,
\beq
\Omega_1\pa{z}=-\sul{p=1}{n} \f{ 1 }{1-\tf{z}{a_p}} \Big\{ \De_{1}^{\pa{p}}\pa{a_p}  \, + \, \f{ 1-\tf{z}{a_p} }{ \ln \big( \tf{z}{a_p}\big) }  \boldsymbol{1}_{D_{a_p,\eps}}(z)  \Big\} \; ,
\enq
with $\boldsymbol{1}_{A}$ being the indicator function of $A$, 
and the expression for $\Omega_2\pa{0}$ is already more involved:
\bem
\Omega_2\pa{0} = - \sul{p=1}{n} \f{\dd }{\dd s} \paa{ \f{\De_2^{\pa{p}}(s)}{2s} \pac{ \f{s-a_p}{\ln \pa{ \tf{s}{a_p}}} }^2 }_{s=a_p} + \sul{p\not=\ell}{n} \f{\De_1^{\pa{p}}\pa{a_p}\De_1^{\pa{\ell}}\pa{a_{\ell}} }{1-\tf{a_{\ell}}{a_p}} \\
 + \sul{p=1}{n}  \f{\dd }{\dd s} \paa{ \De_1^{\pa{p}}(s) \pac{ \f{s-a_p}{  \ln \pa{ \tf{s}{a_p}}} } }_{s=a_p} \cdot  \De_{1}^{\pa{p}}\pa{a_p}  \;\; .
\end{multline}
We get
\beq
\f{T_m\pac{\sg}}{T_{m-1}\pac{\sg}}=G\pac{\sg}\paa{1+\f{\pac{\Omega_1\pa{0}}_{11}}{m} + \f{\pac{\Omega_2\pa{0}}_{11}}{m^2}+ \e{O}\pa{\f{1}{m^3},  \f{ \pa{{\tf{a_i}{a_j}}}^m  m^{3\rho}}{m^3} }  } \; .
\enq
This leads to
\beq
\ln T_m\pac{\sg}=m \ln G\pac{\sg} + \ln m  \sul{p=1}{n}\pa{\ga^2_p-\de_p^2}+K+ \f{Osc}{m^2}\pa{1+\e{o}\pa{1}}+ \f{Nosc}{m}\pa{1+\e{o}\pa{1}} \;\; .
\enq
We recover the first two terms of the asymptotics of $\ln T_m\pac{\sg}$ given in \eqref{asymptotiqueToeplitzFH}. The discrete $m$ derivative method
does not allow to determine the constant $K$, but this is irrelevant in what concerns the structure of corrections.
$Osc$ stands for what we call the oscillating corrections, whereas $Nosc$ for the non-oscillating ones. More explicitly,
\bem
Osc=\sul{p\not=\ell}{} \Ga\pab{1-\ga_{\ell}-\de_{\ell}, 1-\ga_{p}+\de_{p}}{-\ga_{\ell}+\de_{\ell}, -\ga_{p}-\de_{p}}
\f{b_-\pa{a_{\ell}}b_+\pa{a_p}}{b_+\pa{a_{\ell}} b_-\pa{a_p}}  \f{m^{2\de_{\ell}-2\de_{p}}\pa{\tf{a_p}{a_{\ell}}}^{m} }{\pa{\tf{a_p}{a_{\ell}} -1}\pa{\tf{a_{\ell}}{a_p} -1}} \\
\times \f{\pl{\substack{r=1\\ \not=\ell}}{n}\pa{1-\tf{a_r}{a_{\ell}} }^{-\pa{\ga_r+\de_r}}\pa{1 - \tf{a_{\ell}}{a_r} }^{\ga_r-\de_r}    }
{\pl{\substack{r=1\\ \not=p}}{n}\pa{1-\tf{a_r}{a_{p}} }^{-\pa{\ga_r+\de_r}}\pa{1-\tf{a_{p}}{a_r} }^{\ga_r-\de_r}   }
\end{multline}
whereas
\beq
Nosc=\sul{p=1}{n}a_p \pa{\de_p^2-\ga_p^2}\pac{\Dp{z}\ln K_p\pa{z}}_{z=a_p} 
+\f{1}{2}\sul{p\not=\ell}{n} \pa{\de_{\ell}^2-\ga_{\ell}^2}\pa{\de_{p}^2-\ga_{p}^2}\f{\tf{a_{\ell}}{a_p}+1}{\tf{a_{\ell}}{a_p}-1}
\label{ecriture explicite terme sous dominant non oscillant}
\enq

The oscillating corrections have a very nice relationship with the leading term $T^{\pa{0}}\pa{\paa{\ga_i},\paa{\de_i};m}\pac{b}$ of the
asymptotics given in \eqref{asymptotiqueToeplitzFH}. Indeed, one readily checks that,
\beq
Osc=\sul{p \not= q}{n} T^{\pa{0}}\pa{\paa{\ga_i},
\paa{\paa{\de_i}_{\substack{i\not=\\p,q}}, \de_p+1,\de_q-1};m}\pac{b'_{pq}}
\enq
with $b_{pq}'\pa{z}=a_p \tf{b\pa{z}}{a_q}$. Such simultaneous changes $\pa{\de_p,\de_q,b}\mapsto\pa{\de_p+1,\de_q-1,b'_{pq}}$  leave the value of
the symbol $\sg$ given in \eqref{symbole sg FH pour Toeplitz} unchanged. Hence,  this gives strong arguments supporting the Basor-Tracy conjecture as
 already a small part of the different Fisher-Hartwig representations for $\sg$ is present in the sub-leading asymptotics for $T_m\pac{\sg}$.
However, our computations do not allow to give the proof of the Basor-Tracy conjecture in some particular cases where some terms in $Osc$ become of
the same order than $T^{\pa{0}}\pa{\paa{\ga_i},\paa{\de_i};m}\pac{b}$. Although we formally reproduce some particular results of the conjecture
(for instance $\Re\pa{\de_j}=\tf{1}{2}=\Re\pa{\de_k}$ for some j and k) we cannot consider this limiting case as we do not have a control of the remainder.
However, we raise the following generalization of the Basor-Tracy conjecture:

\begin{conjecture}
The sub-leading asymptotics of a Toeplitz matrix with Fisher-Hartwig symbols slowly restore the broken by $T^{0}$ independence with respect to the choice of a Fisher-Hartwig representation for $\sg$. More precisely the asymptotics have the structure
\beq
T_m\pac{\sg} \sim \sul{\substack{n_i \in \mathbb{Z} \\ \Sg n_i=0}}{}\tilde{T}\pa{\paa{\ga_i},\paa{\de_i+n_i};m}\pac{b_{n_i}}
\;\; ,\quad  b_{n_i}\pa{z}=\pl{p=1}{n}a_p^{n_p} \, b\pa{z} \; .
\enq
With
\beq
\tilde{T}\pa{\paa{\ga_i},\paa{\de_i};m}\pac{b}\sim T^{\pa{0}}\pa{\paa{\ga_i},\paa{\de_i};m}\pac{b}\pa{1+\sul{k=1}{+\infty}\f{C_k\pa{\paa{\ga_i},\paa{\de_i},b}}{m^k}}
\enq
having no oscillating terms with $m$.
Note that the notation $\sim$ stands for the equality in the sense of asymptotic series.

\end{conjecture}
This conjecture is a natural extension of the $\nu$ periodicity conjecture raised in \cite{KozKitMailSlaTerRHPapproachtoSuperSineKernel}, and, of course, of the Basor-Tracy conjecture. We stress that one could raise
a similar type of conjecture in what concerns the Wiener-Hopf case.

\section{Conclusion}
In this article we have proven the formula for the leading asymptotic of Fredholm determinants of truncated Wiener-Hopf operators generated by symbols having Fisher-Hartwig singularities. As a byproduct we reproduced, in the framework of Riemann-Hilbert problems, the leading asymptotics of Toeplitz matrices having Fisher-Hartwig singularities. We were also able to compute the first sub-leading asymptotics of Toeplitz matrices having Fisher-Hartwig singularities. These give  support to the Basor-Tracy conjecture. We proposed an extension of the latter conjecture.
Our results were based on a connection between Toeplitz determinants and those of Wiener-Hopf operators: both are related to the so-called generalized sine kernel acting either on the unit circle $\msc{C}$ or the real axis $\R$. In the case of Fisher-Hartwig singularities, this generalized kernel has some jump discontinuities and power-law singularities on the contour.

An open question is the construction of an explicit asymptotic resolvent for truncated Wiener-Hopf operators $I+K$ generated by symbols having Fisher-Hartwig type singularities.
Indeed the resolvent is known in the Fourier space, so it would be enough to take the inverse Fourier transform so as to have the resolvent of $I+K$. Also, it would be interesting to find a proof for the conjecture we have raised above.

%%%%%%%%%%%%%%%%%%%%%%%%%%%%%%%%%%%%%%%%%%%%%%%%%%%%%%%%%%%%%%%%%%%%%%%%%%%%%%%%%%
%%%%%%%%%%%%%%%%%%%%%%%%%%%%%%%%%%%%%%%%%%%%%%%%%%%%%%%%%%%%%%%%%%%%%%%%%%%%%%%%%%

\section*{Acknowledgments}

I would like to thank N. A. Slavnov for introducing me to the subject of Riemann-Hilbert problems. I am also indebted to T. Ehrhardt
who informed me on the actual status of the Fisher-Hartwig conjecture. I am grateful to J.-M. Maillet and N.A. Slavnov for their
numerous and valuable comments. I thank Nick Jones for pointing an error in \eqref{ecriture explicite terme sous dominant non oscillant} present in an 
earlier version of the manuscript. 
This work has been supported by the ANR programm GIMP ANR-05-BLAN-0029-01 and by the French ministry of education and research.

\appendix

\section{Some properties of confluent hypergeometric function}
\label{Appendix CHF formulas}
Tricomi's confluent hypergeometric function $\Psi\pa{a,c;z}$ is one of
the  solutions to the differential equation \cite{BatemanHigherTranscendentalFunctions}:
\beq
z y^{\prime \prime} + \pa{c-z} y^{\prime} -a y=0
\label{CHF eqn}
\enq
For generic $a$ and $c$, $\Psi\pa{a,c;z}$ has a power-law singularity at the origin,
and a cut on $\R^-$. It can be defined, for instance, by its Mellin-Barnes type integral
representation in terms of Euler's $\Ga$ function
\beq
\Psi\pa{a,c;z}= \Int{\ga-i\infty}{\ga+i\infty} \Ga\pab{a+s,-s,1-c-s}{a, a-c+1} z^{s} \f{\dd s}{2i\pi}
\quad ,
\enq
that is valid for $-\Re\pa{a}<\ga<\e{min}\pa{0,1-\Re\pa{c}}$ and $-\tf{3\pi}{2} < \e{arg}\pa{z} < \tf{3\pi}{2}$. The latter integral representation
 is then supplemented by an analytic continuation.  In the above formula we have used  the standard hypergeometric type notation
\beq
\Ga\pab{a_1,\dots , a_n}{b_1, \dots, b_m}= \f{\pl{k=1}{n} \Ga\pa{a_k}}{\pl{k=1}{m} \Ga\pa{b_k}} \; .
\label{Appendix Tricomi CHf definition}
\enq
Tricomi's CHF satisfies the monodromy properties
\beq
 \Psi\pa{a,c;z\ex{2i\pi}}= \ex{-2i\pi a} \Psi\pa{a,c;z}+
\f{2i \pi \ex{-i\pi a+z}}{\Ga\pa{a , 1+a-c} }
 \Psi\pa{c-a,c;-z}\;
\label{cut-Psi-1}
\enq
for $ \Im z<0\;$ and
\beq
 \Psi\pa{a,c;z\ex{-2i\pi}}= \ex{2i\pi a} \Psi\pa{a,c;z}-
 \f{2i \pi \ex{i\pi a+z}}{\Gamma\pa{a, 1+a-c}}
 \Psi\pa{c-a,c;-z}
\label{cut-Psi-2}
\enq
for $ \Im z>0\;$. $\Psi\pa{a,c;z}$ has an asymptotic expansion at $z \tend \infty$ given by
\beq
 \Psi\pa{a,c;z}=\sum_{n=0}^{M} \pa{-1}^n \f{ \pa{a}_n \pa{a-c+1}_n }{ n! }z^{-a-n} \,+\e{O}\pa{z^{-M-a}}\;  ,\qquad -\frac{3\pi}2<\arg(z)<\frac{3\pi}2.
\label{asy-Psi}
\enq
\noindent  Humbert's CHF $\Phi\pa{a,b;z}$ is another solution of \eqref{CHF eqn}.  $\Phi\pa{a,c;z}$ is an entire
function that is defined in terms of its series expansion around $z=0$
\beq
\Phi\pa{a,c;z}=\sul{n=0}{+\infty} \f{\pa{a}_n }{n! \pa{c}_n} z^n ,  \qquad c \not \in \mathbb{Z}^-\;.
\label{definition Humbert's CHF}
\enq
It has the asymptotic expansion around $\infty$
\begin{multline}
\Phi\pa{a,c;z} = \f{\Ga\pa{c}}{\Ga\pa{c-a}} \paf{\ex{i\pi \eps}}{z}^{a}
\sul{n=0}{M} \f{\pa{a}_n\pa{a-c+1}_n}{n!\pa{-z}^n}  + \e{O}\pa{\abs{z}^{-a-M-1}}\\
 +\f{\Ga\pa{c}}{\Ga\pa{a}} \ex{z} z^{a-c} \sul{n=0}{N} \f{\pa{c-a}_n \pa{1-a}_n}{n! z^n}
+ \e{O}\pa{\abs{\ex{z}z^{a-1-c-N}}}.
\label{Asymptotiques Humbert's CHF}
\end{multline}

There are many relations between these two different CHF. In particular
\bem
 \ex{\eps\pa{c-a+1}i\pi} z^{1-c}\Phi\pa{a-c+1,2-c;z}= \\
 \paa{\Ga\pab{2-c}{1-a} \Psi\pa{a,c;z} -
                    \Ga\pab{2-c}{1+a-c}\ex{z} \Psi\pa{c-a,c;-z} }
\label{Phi s'ecrit comme Psi}
\end{multline}
\noindent where $\eps= \e{sgn}\pa{\Im z}$ and it is assumed that $\arg\pa{z}\in \intoo{-\tf{\pi}{2}}{\tf{\pi}{2}}$.
One can also express $\Psi\pa{a,c;z}$ in terms of $\Phi\pa{a,c;z}$ thus allowing to access to the singularity structure of Humbert's CHF at the origin:
\beq
\Psi\pa{a,c;z}= \Ga\pab{1-c}{a-c+1} \Phi\pa{a,c;z}  + \Ga\pab{c-1}{a} z^{1-c} \Phi\pa{a-c+1,2-c;z}  \; .
\label{Psi s'ecrit comme Phi}
\enq

\noindent Actually a CHF is some limiting case of the Gauss hypergeometric function. This function is one of the
solutions of the hypergeometric equation. We recall its series expansion around $z=0$:
\beq
_2 F_1\pa{ \ba{cc} a,b \\ c \ea ; z } = \sul{n=0}{+\infty} \f{\pa{a}_n \pa{b}_n}{n! \pa{c}_n} z^n \;\; ,  \qquad c \not \in \mathbb{Z}^- \;\; .
\enq
\noindent The above solution is regular around $z=0$ and can be continued to large value of $z$ thanks to the identity
\bem
_2 F_1\pa{ \ba{cc} a,b \\ c \ea ; z }= \Ga \pab{c, b-a}{b,c-a} \pa{-z}^{-a}\, _2 F_1\pa{ \ba{cc} a,1+a-c \\ 1+a-b \ea ; z^{-1} }  \\
+ \Ga \pab{c, a-b}{a,c-b} \pa{-z}^{-b} \, _2 F_1\pa{ \ba{cc} b,1+b-c \\ 1+b-a \ea ; z^{-1} } \;\;  .
\label{Fonction de Gauss et continuation Analytique}
\end{multline}
 One can also consider multi-variable generalizations of hypergeometric functions. We give here the definition
of the Appell function of the second kind in terms of a double series that is convergent provided that
$\abs{y}+\abs{z}<1$:
\beq
F_2                 \pa{a; \ba{c} b \, , c \\
                                    d \, , e \ea ; y, z}=  \sul{n,m \geq 0}{}
\f{\pa{a}_{m+n} \, \pa{b}_n \, \pa{c}_m }{n! \, m!\, \pa{d}_n \pa{e}_m } y^n \, z^m \; .
\label{fonction de Appell}
\enq

We finally point out that the Barnes G-function admits an integral representation in terms of $\psi$, the
logarithmic derivative of Euler's Gamma function.
\beq
G\pa{z+1}= \pa{2\pi}^{\f{z}{2}}\exp\paa{ - \f{z\pa{z-1}}{2}+  \Int{0}{z} t \psi\pa{t} \dd t }  \;\; \Re\pa{z}>-1 \;\; .
\label{fonction de Barnes}
\enq

\section{Integrals of CHF}
\label{Appendix integration CHF}
\subsection{Series Expansion of Appell function of large arguments.}

Using the Mellin-Barnes type integral representation for CHFs \cite{BatemanHigherTranscendentalFunctions}, Erdelyi was able to
evaluate the Laplace transform of products of CHFs in terms of Lauricella's function \cite{ErdelyiTablesOfIntTransf}. The Lauricella function associated
to an integral involving a product of two Humbert's CHF is better know as the Appell function of the second kind \eqref{fonction de Appell}.
In terms of this function, Erdelyi's result reads:
\bem
f\pa{p,\a_1,\a_2,\eps;s}\equiv\Int{0}{+\infty} \f{\dd t \, \ex{-st}}{ t^{2\ga+1-\eps-p}}
 \Phi\pa{-\de-\ga,1-2\ga;-i\a_1t}\Phi\pa{\de-\ga,1-2\ga;i\a_2t} \\
=s^{2\ga-p - \eps} \, \Ga\pa{\eps+p-2\ga}
\,
 F_2\pa{p+\eps-2\ga; \ba{ccc} -\ga-\de &,& \de-\ga \\
                                                1-2\ga &,& 1-2\ga \ea ; -i\f{\a_1}{s},\, i\f{\a_2}{s}} \; .
\label{integrals of CHF and Appell functions}
\end{multline}
Such integrals have been considered in \cite{KarlueIntegralsOf2CHF}, in the case where $\eps$ and $p$ are integers
and the integral is absolutely convergent. Here, we study the behavior of such integrals when $\eps$ is close to zero.  In that case, one cannot apply the integration procedure presented in \cite{KarlueIntegralsOf2CHF} as it only applies to integer $\eps$.  Moreover, we
consider products of CHF that aren't decaying sufficiently fast at infinity. Hence  one should
regularize the integrals before taking the $s\tend 0^+$ limit. Once a regularization is performed, this
limit can be computed thanks to a series expansion of the second Appell function around $\infty$

\begin{figure}[h]
\begin{center}

\includegraphics{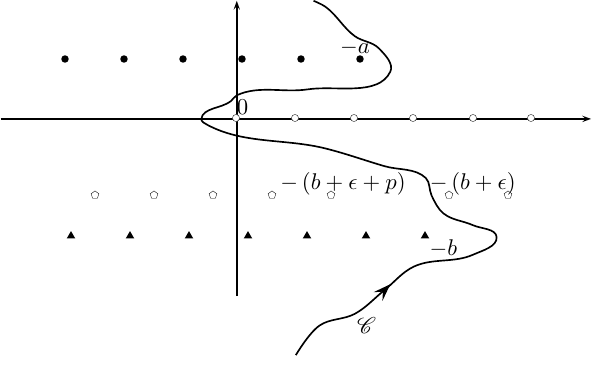}

\caption{Contour of integration for the Appel function.\label{contour pour la fonction d'Appel}}
\end{center}
\end{figure}

We study  $f\pa{p,\a_1,\a_2,\eps;s}$ as it is the generating function for all the integrals that appear in the evaluation
of the trace of the resolvent. The precise procedure for computing these integrals will be explained in this Appendix.
The idea is to derive a series expansion for $F_2$ at $\infty$ and then use it to compute,  after a proper regularization, the $s\tend 0^+$ limit of $f\pa{p,\a_1,\a_2,\eps;s}$.
\begin{lemme}
\label{developpement asymp fonction Appell}
Let $f\pa{p,\a_1,\a_2,\eps;s}$ be defined in terms of the Appell function of the second kind as in \eqref{integrals of CHF and Appell functions}.
Then $f$ admits a series expansion around $s=0$:
\bem
\Ga\pab{-\ga-\de,\de-\ga}{1-2\ga,1-2\ga} f\pa{p,\a_1,\a_2,\eps;s}=S_1\pa{p,\a_1,\a_2,\eps;s}\\
+ S_2\pa{p,\a_1,\a_2,\eps;s}+S_3\pa{p,\a_1,\a_2,\eps;s}
\end{multline}
Where the $S_i$ are series involving Gauss' functions
\bem
S_1\pa{p,\a_1,\a_2,\eps;s}=
\f{\ex{i\f{\pi}{2}\pa{p+\eps-2\ga}}}{\a_2^{n+p+\eps-2\ga}} \sul{n\geq0}{} \paf{\a_1}{\a_2}^{n}
 \!\! \!\!\!\; _2F_1 \pa{\ba{c} p+n -2\ga+\eps \, , p+n+\eps \\ 1-\ga-\de+p+n+\eps \ea ; -i\f{s}{\a_2}  }  \\
\times \Ga\pab{n+p+\eps-2\ga\, , \de+\ga -p-\eps-n \, , -\de - \ga +n}{1+n-2\ga \, , 1-p-n-\eps \, , 1+n}
\label{definition S1}
\end{multline}
\bem
S_2\pa{p,\a_1,\a_2,\eps;s} = \f{\ex{i\pi\de}}{ \ex{i\f{\pi}{2}\pa{p+\eps}}} \sul{n\geq1}{}
\, _2F_1 \pa{\ba{c} n+ p+\de-\ga \, , \ga+\de+p+n \\ 1+p+n \ea ; -i\f{s}{\a_2}  } \\
\times \f{\a_1^{n+\de+\ga-\eps} }{\a_2^{n+p+\de-\ga}} \; \Ga\pab{n+p+\de-\ga\, , n-\eps \, , \eps -\de-\ga -n }{1+n+\de-\ga-\eps \, , 1-\ga-\de-p-n \, , 1+p+n}
\label{definition S2}
\end{multline}
\bem
S_3\pa{p,\a_1,\a_2,\eps;s}= \ex{i\pi\pa{\de+p}} \ex{-i\f{\pi}{2}\eps}  \Ga\pab{-\ga+\de}{1-\ga-\de}
\f{\pi \a_1^{\ga+\de-p} \a_2^{\ga-\de} }{\sin \pa{\pi \eps}s^p}   \\
\times \sul{n\geq0}{} \paf{s}{i\a_1}^{n} \left[
\paf{i}{s}^{\eps} \Ga\pab{n-\ga-\de}{1+\de-\ga-n\, ,1+n}
\; _2\Phi_1\pa{\ba{c} \de-\ga \, ,  \ga+\de  \\ 1+n-p-\eps \ea; -i\f{s}{\a_2} } \right.
 \\
\left. - \a_1^{-\eps}\Ga\pab{n-\ga-\de+\eps}{1+\de-\ga-n-\eps\, ,1+n+\eps}
\; _2\Phi_1\pa{\ba{c} \de-\ga \, ,  \ga+\de  \\ 1+n-p \ea; -i\f{s}{\a_2} }   \right]
\label{definition S3}
\end{multline}
\end{lemme}
\Proof
We consider the Appell function of the second kind as in \eqref{fonction de Appell} and assume
 the following dependence between the parameters
\beq
a=p+\eps+ b+ c \;\;\; \e{and} \;\;\; f=d=1+b+c=1+a-p-\eps
\enq
with $\eps$ small and complex.
The series expansion of the second Appell function can be re-summed into a Mellin-Barnes type integral representation
\cite{SlaterGeneralizedHypergeormetricFunctions}:
\bem
\Ga\pab{a \, , b \, }{d \, } F_2\pa{a  ; \ba{ccc} b &c\\
      d & f\ea ; x, y}=  \Int{\mathscr{C}}{} \f{\dd s}{2i\pi} \Ga\pab{a+s, b+s, -s}{f+s} \\
       \times \, _2 F_1\pa{ \ba{c} a+s, c \\ f \ea ; y } \pa{-x}^{s} \;\; .
\label{Barnes type representation for Appel function}
\end{multline}
One possible choice of the contour $\msc{C}$  is depicted in
fig. \ref{contour pour la fonction d'Appel}
 and the integral is convergent provided $\abs{\arg\pa{-x}}+\abs{\Im\pa{y}}<\pi$.  This can be seen using
the asymptotics of $\Ga\pa{z}$ in the region $\abs{\arg\pa{s}}<\pi$. One has
\beq
\Ga\pab{a+s \, , b+s \, ,-s }{f+s} = 2i\pi \e{sgn}\pa{\Im\pa{ s} } s^{a+b-f-1} \ex{-\pi \abs{\Im \pa{s}}} \pa{1+ \e{O}\pa{s^{-1}}} \, , \,
s\tend + i \infty \; .
\enq
\noindent We remind the asymptotics behavior of a hypergeometric function of a large argument \cite{BatemanHigherTranscendentalFunctions}:
\beq
_2F_1\pa{ \ba{c c}a+s \, , c \\ f \ea ; y  } =  \paa{\ex{i\pi c} \Ga\pab{f}{f-c} \pa{s y}^{-c} + \ex{s y +a y }\pa{s y}^{c-f}}
\pac{1+ \e{O}\pa{s^{-1}}} \; .
\enq
\noindent So that, all together
\bem
\abs{\Ga\pab{a+s \, , b+s \, ,-s }{f+s} \; _2F_1\pa{ \ba{c c}a+s \, , c \\ f \ea ; y  } \pa{-x}^{-s} }  \\
\leq C \abs{s^{a+b+2f-1+c}} \exp \paa{- \Im \pa{s}  \e{arg}\pa{-x} + \abs{\Im \pa{y} \, \Im \pa{s}} -\pi \abs{\Im\pa{ s}}} \; ,
\end{multline}
\noindent where C is some computable constant. The integrand of \eqref{Barnes type representation for Appel function}
is thus absolutely integrable provided $\abs{\arg\pa{-x}} + \abs{\Im\pa{ y} }<\pi$.

We now split \eqref{Barnes type representation for Appel function} in two by using the analytic continuation of
Gauss function for for large $y$ \eqref{Fonction de Gauss et continuation Analytique}:
\bem
\Ga\pab{a \, , b \, , c }{d \, , f } F_2\pa{a  ; \ba{ccc} b &c\\
      d & f\ea ; x, y}=  L_1+L_2 \\
L_1=\! \Int{\mathscr{C}}{} \!\!\f{\dd s}{2i\pi} \Ga\! \pab{a+s, b+s, -s,c-a-s}{d+s,f-a-s}
      \!\!  \, _2 F_1\pa{ \ba{c} a+s, a+s+1-f \\ a+1+s-c \ea ; \f{1}{y}  } \f{\pa{-x}^{s}}{\pa{-y}^{a+s}} \\
L_2=\Int{\mathscr{C}}{} \f{\dd s}{2i\pi} \Ga\pab{a-c+s, b+s, -s,c}{d+s,f-c}
       \, _2 F_1\pa{ \ba{c} c+1-f, c \\ c+1-a-s \ea ; y^{-1} } \f{\pa{-x}^{s}}{\pa{-y}^{c}}
\end{multline}

We were able to split the integral into two parts as each of them converges separately. The separate convergence
of the integrals is readily seen from the asymptotic of the Gauss hypergeomertic function \cite{BatemanHigherTranscendentalFunctions}:
\bem
_2F_1\pa{ \ba{c} \a+s, \beta+s \\ \ga+s \ea ; z }=
\pa{1-z}^{\ga-\a-\beta-s} \;\; _2F_1\pa{ \ba{c} \ga-\a, \ga-\beta \\ \ga+s \ea ; z } \\
\underset{s \tend +\infty}{\sim}
\pa{1-z}^{\ga-\a-\beta-s}\pa{1+\e{O}\pa{s^{-1}}} \quad, \;\;  \abs{z}<1.
\end{multline}
\noindent Hence putting
\beq
g\pa{x,y;s} = \Ga\pab{a+s,b+s,c-a-s,-s}{d+s,f-a-s}\!\!\! \, _2F_1\pa{ \ba{c} a+s, a-f+s+1 \\ a-c+s+1 \ea ; \f{1}{y}} \! \f{\pa{-x}^{s}}{\pa{-y}^s}
 \; ,
\enq
\noindent and using the asymptotic behavior of the Gamma function we get
\beq
\abs{g\pa{x,y;s}} \leq C \abs{s}^{\Re\pa{a+b+c-d-f-1}} \ex{-\pi\abs{\Im \pa{s}} -\Im\pa{s}\pa{\e{arg}\pa{-x}-\e{arg}\pa{1-y^{-1}}-\e{arg}\pa{-y}}} \; .
\enq
$L_1$ is thus convergent in the region defined by the equation
\beq
\abs{ \e{arg}\pa{-x}-\e{arg}\pa{1-y^{-1}}-\e{arg}\pa{-y} }<\pi \; ,\quad  \e{and} \quad  \abs{y}>1\; .
\enq
Similar calculation lead to the conclusion that $L_2$ is convergent  in the region \newline $\abs{\e{arg}\pa{-x}}< \pi$.

$L_1$ can be computed as a sum over the  poles located at the right of $\msc{C}$.
These are $s=n \; , \; n \in \mathbb{N}$ and $s=-b-\eps+n \; , \; n\in \mathbb{N}^*$. One eventually arrives to
\bem
L_1=-\f{\sin \pi \eps}{ \sin \pi\pa{b+\eps}} \pa{-y}^{-a} \sul{n\geq 0}{} \paf{-x}{y}^n
\Ga\pab{a+n,b+n,p+\eps+n}{d+n,1+n+\eps+p+b,1+n} \\
 \times \; _2F_1\pa{ \ba{c} a+n, p+n+\eps \\ 1+b+\eps+n+p \ea ; y^{-1} }
-\f{\sin \pi b}{ \sin \pi\pa{b+\eps}} \pa{-x}^{-b-\eps}\pa{-y}^{-p-c}\sul{n\geq 1}{}  \paf{-x}{y}^n \\
\times
\Ga\pab{p+c+n,n-\eps,p+n-b}{1+c-\eps+n,1+n-\eps-b,1+p+n} \; _2F_1\pa{ \ba{c} p+c+n, p+n-b \\ 1+n+p \ea ; y^{-1} }
\end{multline}

Similarly $L_2$ can be computed
as a sum over the poles located to the left of $\msc{C}$; these are $\paa{-b-n, \,-b-n-\eps}$
 where $n \in \mathbb{N} $. This becomes apparent when one normalizes the Gauss function:
\beq
_2F_1\pa{ \ba{c} a, b \\ c \ea ; z } =\Ga\pa{c}  \, _2\Phi_1\pa{ \ba{c} a, b \\ c \ea ; z } \; ,
\enq
\noindent so that $_2\Phi_1$ is an entire function of the parameters $a, b$ and $c$. The result reads:
\bem
\f{\pa{-1}^{p}\pa{-y}^{-c}}{\pa{-x}^{b+\eps}}\f{\sin\pi \pa{d-b}}{\sin\pi \eps} \sul{n\geq 0}{}
 \f{\pa{-x}^{\eps}}{x^n}\Ga\pab{b+n,n-c}{1+n}\,  _2\Phi_1\pa{ \ba{c} -b, c \\1+n-p-\eps \ea ; y^{-1} } \\
-\f{\sin\pi \pa{d-b-\eps}}{x^n \sin\pi \pa{d-b}} \Ga\pab{b+n+\eps,\eps+n-c}{1+\eps+n} \; _2\Phi_1\pa{ \ba{c} -b, c \\1+n-p \ea ; y^{-1} }
\end{multline}

The joint condition for the convergence of $L_1$ and $L_2$ defines an open subset O of $\Cx^2$.
Hence we can continue the series representation for the second Appell function to the largest open subset in $\Cx^2$ containing O where the series is convergent. In particular, the series representation is well defined for the range of parameters that we use.
Specifying the values of $a,b,c,d,f$ to the ones of the Lemma we obtain the claimed result. $\Box$

\subsection{Useful integrals}
We will use the series expansion  for $F_2$ in order to compute some integrals of products of $CHF$.
We remind the definitions of the functions $\tau\pa{\ga,\de;t}$ and $\varphi\pa{\ga,\de;t}$
\bem
\Ga\pab{1-2\ga,1-2\ga}{1+\de-\ga,1-\de-\ga}\tau\pa{\ga,\de;t}=-\Phi\pa{-\ga-\de,1-2\ga;-it}\Phi\pa{\de-\ga,1-2\ga;it}
\nonumber \\
+
\pa{\Dp{z}\Phi}\pa{-\ga-\de,1-2\ga;-it}\Phi\pa{\de-\ga,1-2\ga;it} \nonumber \\
 + \Phi\pa{-\ga-\de,1-2\ga;-it}\pa{\Dp{z}\Phi}\pa{\de-\ga,1-2\ga;it}
\end{multline}
\noindent and
\beqa
\Ga\pab{1-2\ga,1-2\ga}{1+\de-\ga,1-\de-\ga}\varphi\pa{\ga,\de;t}=\Phi\pa{-\ga-\de,1-2\ga;-it}\Phi\pa{\de-\ga,1-2\ga;it} \;\; . \nonumber
\eeqa
\begin{cor}
\label{Corolaire integrales CHF}
Let $\abs{\Re\pa{\de}}<\tf{1}{2}$, $\Re\pa{\ga}<\tf{1}{2}$ and $\tau\pa{\ga,\de;t}$, $\varphi\pa{\ga,\de;t}$ be as above, then
\bem
\Int{-\infty}{0} \dd t  \pa{\f{\ex{i\pi\de}}{\abs{t}^{2\ga}} \varphi\pa{\ga,\de;t}-1} =
2i\de  \quad , \quad
\Int{\R}{} \dd t \pa{\f{\ex{-i\pi\de \e{sgn}\pa{t}}}{\abs{t}^{2\ga}} \varphi\pa{\ga,\de;t} -1}  =0
\nonumber \\
 \Int{-\infty}{0} \dd t \pa{ \f{\ex{i\pi\de}}{\abs{t}^{2\ga}} \tau\pa{\ga,\de;t} +1-\f{2i\de}{t-1} } =
-2i\de   -\pi\ga  +i\pa{\ga+\de}\psi\pa{-\ga-\de} +i\pa{\de-\ga} \psi\pa{\de-\ga}  \\
\Int{\R}{}\dd t  \pa{  \f{ \ex{-i\pi \e{sgn}\pa{t} \de} }{\abs{t}^{2\ga}}  \tau\pa{\ga,\de;t}
+1-\f{2i\de }{t+\e{sgn}\pa{t}}} =-2 \pi\ga \\
\Int{\R}{}\dd t  \ln\pa{\abs{t}} \pa{\abs{t}^{-2\ga} \varphi\pa{\ga,0,t}-1}=- 2\pi\ga   \\
\Int{\R}{} \dd t \ln\pa{\abs{t}} \pa{\abs{t}^{-2\ga}\tau\pa{\ga,0,t}+1}=
 2 \pi \ga  \pa{ \psi\pa{1-\ga}-2\psi\pa{1-2\ga}+1}  \label{integrales de CHF}
\end{multline}
\noindent The Riemann integrability of the integrands is part of the conclusion.
\end{cor}
\Proof
Using the asymptotic behavior of Tricomi's CHF \eqref{asy-Psi} as well as \eqref{Phi s'ecrit comme Psi} one can readily convince oneself that for $\abs{\Re\pa{\de}}<\tf{1}{2}$
\bem
\Phi\pa{-\de-\ga,1-2\ga;-it}\Phi\pa{\de-\ga,1-2\ga;it}\underset{t\tend+\infty}{\sim} \\
\ex{i\pi\de}t^{2\ga}
\Ga\pab{1-2\ga,1-2\ga}{1+\de-\ga, 1-\de-\ga} \paa{1+\e{O}\paf{\ex{\pm it}}{t^{1\mp2\de}} }
\end{multline}
\bem
\pa{\Dp{z}\Phi}\pa{-\de-\ga,1-2\ga;-it}\Phi\pa{\de-\ga,1-2\ga;it} \\
+\Phi\pa{-\de-\ga,1-2\ga;-it}\pa{\Dp{z}\Phi}\pa{\de-\ga,1-2\ga;it}
\underset{t\tend+\infty}{\sim}\\
\ex{i\pi\de} t^{2\ga} \Ga\pab{1-2\ga,1-2\ga}{1+\de-\ga, 1-\de-\ga} \paa{\f{2i\de}{t}+\e{O}\paf{\ex{\pm it}}{t^{1\mp2\de}}}\; .
\end{multline}
\noindent Where the terms that are sub-leading to $\ex{\pm it}t^{\pm 2\de-1}$ are already absolutely integrable.
For $\abs{\Re\pa{ \de}}<\tf{1}{2}$, $\ex{\pm it}t^{\pm2\de-1}$ is only Riemann integrable, thus all integrals should be understood in this sense. Moreover the asymptotics in the region $t\tend -\infty$ can be inferred from those at $t \tend +\infty$ if one starts with the complex conjugated parameters $\de^*$ and $\ga^*$ and then takes the complex conjugate of the asymptotic series. This settles the question about the Riemann integrability of the different integrands.

The proof of this Corollary is straightforward although quite long. The principle of the proof is to
extract the divergent and constant terms in the $s \tend 0^+$ limit from the expansion of $f\pa{p,\a_1,\a_2,\eps;s}$ around $s\tend 0^+$.

We shall explain in detail how to obtain  the first four integrals. The remaining two are obtained in a similar fashion, although computations become more and more involved. At the end of the proof we list all the summation
identities that are necessary to compute the $s\tend 0^+$ limit in the other cases.
Recall the notation introduced in  \eqref{integrals of CHF and Appell functions}. Then one has
\bem
\Int{0}{+\infty}\f{\dd t}{t^{2\ga}}
\paa{ \varphi\pa{\ga,\de,t} - t^{2\ga} \ex{i\pi\de}}=
\lim_{s \tend 0^+}\Int{0}{+\infty}\f{\dd t}{t^{2\ga}} \ex{-s t}
\paa{ \varphi\pa{\ga,\de,t} - t^{2\ga} \ex{i\pi\de}}\\
= \Ga\pab{1+\de-\ga,1-\de-\ga}{1-2\ga,1-2\ga} \pa{ f\pa{1,1,1,0;s}-\f{\ex{i\pi\de}}{s}} \; .
\label{integral varphi Reimm int evaluation}
\end{multline}
\noindent  Since we compute the $s\tend 0^+$ limit, it is enough to determine $f\pa{1,1,1,0;s} $ up to $\e{o}\pa{1}$
with respect to $s\tend 0^+$. One easily sees that $S_1\pa{1,1,1,0;s}=0$, \textit{cf} \eqref{definition S1}.
Since we compute $S_2\pa{1,1,1,0;s}$ \eqref{definition S2} up to $\e{O}\pa{s}$ terms, we can already replace Gauss' function by 1; the latter only contributes to higher orders terms in s. Then, for this particular choice of the constants $p,\a_1,\a_2$ and $\eps$, the second sum boils down to
\beq
S_2\pa{1,1,1,0;s}=-i\ex{i\pi\de} \sul{n\geq 1}{} \Ga\pab{n}{n+2}+ \e{O}\pa{s} =-i\ex{i\pi\de} + \e{O}\pa{s} \; .
\enq
Finally, we estimate $S_3\pa{1,1,1,0;s}$ \eqref{definition S3}. This term is the most complicated one. Indeed, due to the presence of the factor $\tf{1}{\pa{s \sin\pi \eps}}$ in
front of the sum,  one must compute the linear in $\eps$ terms of the sum (the zeroth order vanishes as it should be). Also, it is enough to expand
the sum up to $\e{O}\pa{s\ln s}$ as such terms won't contribute to the result after the  $s\tend 0^+$ limit is performed.  After some computations
one gets
%
%
%
%\beq
%
%S_3= \pa{\eps\Ga\pab{-\ga-\de}{1-\ga+\de}\pa{-1+2i \de s } -is \eps \Ga\pab{1-\ga-\de}{\de-\ga} +\e{O}\pa{\eps,\ln s \, s^2} }
%
%\enq
%
%
%
%So that the the last line in \eqref{integrals of CHF and Appell functions} for this specific choice of parameters %gives a contribution
%
%
%
\beq
S_3\pa{1,1,1,0;s}=\Ga\pab{\de-\ga, -\de-\ga}{1+\de-\ga, 1-\de-\ga} \ex{i\pi\de} \pa{s^{-1}-2i\de} +i\ex{i\pi\de}\; + \e{O}\pa{s \ln s}.
\enq
Adding up all the three contributions, we see that the $s^{-1}$ part cancels with the one coming from the regularization term in
\eqref{integral varphi Reimm int evaluation}. The remaining terms combine to give
\beq
\Int{0}{+\infty}\f{\dd t}{t^{2\ga}}
\paa{ \varphi\pa{\ga,\de,t} - t^{2\ga} \ex{i\pi\de}}=-2i\de \ex{i\pi\de}\; .
\label{integral CHF la plus simple}
\enq
Now, the first integral in the list of \eqref{integrales de CHF} is obtained by considering \eqref{integral CHF la plus simple} in the case of
parameters $\de^*$ and $\ga^*$,  changing variables $t\tend -t$, and then taking the complex conjugate of the whole expression.

We now explain how to evaluate the integrals involving $\tau\pa{\de,\ga;t}$. Since we have already established the value of integrals involving
 $\varphi\pa{\de,\ga;t}$, we only need to compute
\bem
\ex{-i\pi\de}\Int{0}{+\infty} \f{\dd t}{ t^{2\ga} }
\left\{ \pa{\Dp{z}\Phi}\pa{-\ga-\de,1-2\ga;-it}\, \Phi\pa{\de-\ga,1-2\ga;it}  \right. \\
\hspace{2cm}+ \Phi\pa{-\ga-\de,1-2\ga;-it}\, \pa{\Dp{z}\Phi}\!\pa{\de-\ga,1-2\ga;it}\\
\hspace{4.5cm}\left. -\f{2i\de t^{2\ga}\ex{i \pi \de}}{t+1} \Ga\pab{1-2\ga,1-2\ga}{1+\de-\ga,1-\de-\ga}  \!\!  \right\} \\
= -i  \lim_{s\tend 0^+} \left\{ \pac{\Dp{\a_2}f\pa{0,1,1,0;s} -\Dp{\a_1}f\pa{0,1,1,0;s}} \ex{-i\pi\de} \right. \hspace{2cm}\\
\left. +2\de\, \Ga\pab{1-2\ga,1-2\ga}{1+\de-\ga,1-\de-\ga}  \Psi\pa{1,1;s} \right\} \;\; .
\label{deuxieme integrale CHF}
\end{multline}
In this case, only the first term appearing in $S_1$ contributes:
\beq
\pac{\pa{\Dp{\a_2}-\Dp{\a_1}}.S_1}\pa{0,1,1,0;s}=i\ex{-i \pi \ga} \Ga\pab{1-2\ga,1-2\ga,\de+\ga}{\de-\ga}+ \e{O}\pa{s}
\; .
\enq
Already for $S_2$ one has to compute some less trivial sums
\bem
\pac{\pa{\Dp{\a_2}-\Dp{\a_1}}.S_2}\pa{0,1,1,0;s}\\
=i\ex{i\pi \de}  \Ga\pab{1-2\ga,1-2\ga}{-\de-\ga , \de-\ga}
\paa{\sul{n\geq 1}{} \f{1}{n\pa{n+\de-\ga}}+\f{1}{n\pa{n+\de+\ga}} } \pa{1+\e{O}\pa{s}}\hspace{1.2cm} \\
=i\ex{i\pi \de}  \Ga\pab{1-2\ga,1-2\ga}{-\de-\ga , \de-\ga}
\! \paa{\! \f{\psi\pa{1}-\psi\pa{1+\de-\ga}}{\ga-\de}+ \f{\psi\pa{1+\de+\ga}-\psi\pa{1}}{\ga+\de}\!}\!\pa{1+\e{O}\pa{s}} \; .
\end{multline}
And we have used $\sul{n\geq 1}{} \f{1}{n\pa{n-a}}=\f{\psi\pa{1}-\psi\pa{1-a}}{a}$. Finally,
\bem
\pac{\pa{\Dp{\a_2}-\Dp{\a_1}}.S_3}\pa{0,1,1,0;s}=-\ex{i\pi\de}\Ga\pab{\de-\ga,-\de-\ga}{1+\de-\ga,1-\de-\ga} \\
\paa{2\de\pac{\ln\pa{\tf{i}{s}}+2\psi\pa{1}-\psi\pa{-\de-\ga}-\psi\pa{1+\de-\ga}  }+1  }\pa{1+\e{O}\pa{s\ln s}}\; .
\end{multline}
Adding together the three contributions we get
\bem
\ex{-i\pi\de} \Ga\pab{1+\de-\ga,1-\de-\ga}{1-2\ga,1-2\ga} \pac{\pa{\Dp{\a_2}-\Dp{\a_1}}.f}\pa{0,1,1,0;s}= \\
-2\de\pac{ \ln\pa{\tf{i}{s}}+\psi\pa{1}-\psi\pa{-\de-\ga}-\psi\pa{1+\de-\ga}   } -\pa{\ga+\de}\psi\pa{\de-\ga} \\
+\f{2\de}{\ga-\de}+\pa{\ga-\de}\psi\pa{1+\de+\ga} + \f{ \pi\pa{\de-\ga}\ex{-i\pi\pa{\ga+\de}} }{\sin \pi \pa{\ga+\de}  }  + \e{O}\pa{s\ln s} \\
=2\de\pa{\ln s-\psi\pa{1}} -i\pi\ga +\pa{\ga+\de}\psi\pa{-\de-\ga} + \pa{\de-\ga}\psi\pa{\de-\ga} +\e{O}\pa{s\ln s} \; .
\end{multline}
We used the addition formulae for the $\psi$ function
\beq
\psi\pa{1+z}-\psi\pa{z}=\f{1}{z} \quad ,\quad \psi\pa{1+z}-\psi\pa{-z}=-\pi\cot\pi z \; ,
\enq
in order to obtain the last line.
The leading asymptotics of Tricomi's CHF around zero $\Psi\pa{1,1;s}=-\ln s+\psi\pa{1}+\e{O}\pa{s \ln s}$ allows to take the $s\tend 0^+$ limit
 in \eqref{deuxieme integrale CHF}. We get
\bem
\ex{-i\pi\de}\Int{0}{+\infty} \f{\dd t}{t^{2\ga}} \paa{\varphi\pa{\ga,\de;t}+\tau\pa{\ga,\de;t}-\f{2i\de \ex{i\pi\de}}{t+1} } \\
=-\pi\ga -i\pa{\ga+\de} \psi\pa{-\de-\ga}-i\pa{\de-\ga}\psi\pa{\de-\ga} \; .
\end{multline}
The other integrals involving $\tau$ are then obtained from the latter results by the standard manipulations that we have already described.

The value of the last two integrals appearing in the Corollary is obtained by a similar procedure. Namely,
\bem
\Int{0}{+\infty} \f{\dd t\, \ln t}{t^{2\ga}} \paa{\varphi\pa{\ga,0;t} - t^{2\ga} } \\
=\lim_{s\tend 0^+} \paa{ \Ga\pab{1-\ga,1-\ga}{1-2\ga,1-2\ga} \pa{\Dp{\eps}f}\pa{1,1,1,0;s}-\f{\psi\pa{1}-\ln s}{s} } \; .
\end{multline}
A long but straightforward computation yields
\beq
\Ga\pab{1-\ga,1-\ga}{1-2\ga,1-2\ga}\pa{\Dp{\eps}f}\pa{1,1,1,0;s}=-\pi\ga+\f{\psi\pa{1}-\ln s}{s} +\e{O}\pa{s\ln s} \; .
\enq
Similarly,
\bem
\Int{0}{+\infty} \f{\dd t \ln t}{t^{2\ga}} \paa{\varphi\pa{\ga,0;t}+\tau\pa{\ga,0;t}} \\
=\lim_{s\tend 0^+} -i\Ga\pab{1-\ga,1-\ga}{1-2\ga,1-2\ga} \pac{ \pa{\Dp{\a_2\, , \eps}^2-\Dp{\a_1 \, ,\eps}^2} f}\pa{0,1,1,0;s}\; .
\end{multline}
An even longer but as much straightforward computation leads to
\bem
\pac{\pa{ \Dp{\a_2 \eps}^2-\Dp{\a_1\eps}^2 }S_1}\pa{0,1,1,0;s}=
\f{\pi \ex{-i\pi\ga}}{\ga\sin\pi\ga} \paa{2\ga^{-1}\psi\pa{1-2\ga}-\psi\pa{\ga}+i\f{\pi}{2} }
\end{multline}
\bem
\pac{\pa{\Dp{\a_2  \eps}^2-\Dp{\a_1 \eps}}S_2}\pa{0,1,1,0;s}=
\f{i\pi}{2\ga} \pa{\pi\cot\pi\ga-\ga^{-1}} -\f{1}{2\ga^3}+\f{\pi^2}{2\ga} -\f{3\pi^2}{2\ga}\cot^2\pi\ga \\
+\ga^{-2}\pac{\psi\pa{1-\ga}+2\psi\pa{-\ga} -\psi\pa{1}} -2 \pi \ga^{-1} \cot\pi\ga \, \psi\pa{1+2\ga} \\
+\f{\pac{\psi\pa{1+\ga}+\psi\pa{1-\ga}}}{2\ga}\pac{\pi\cot\pi\ga-\ga^{-1} }
\end{multline}
\bem
\pac{\pa{\Dp{\a_2  \eps}^2-\Dp{\a_1 \eps}}S_3}\pa{0,1,1,0;s}=
\ga^{-2}\pac{\tf{i\pi}{2}+\psi\pa{1}-\psi\pa{1-\ga}-\psi\pa{-\ga}  } \; .
\end{multline}
Where we have dropped the $\e{O}\pa{s\ln s}$ symbol so as to lighten the formulae a little. In the intermediary
 computations  of the contribution issued from $S_2$ term we used the formulae below
\beqa
\sul{n\geq 1}{} \f{1}{n^2-\ga^2}&=& \f{\psi\pa{\ga}-\psi\pa{-\ga}}{2\ga} \; ; \\
 \sul{n \geq 1}{} \f{1}{\pa{n+\ga}\pa{n-\ga}^2}&=&\f{\psi\pa{-\ga}-\psi\pa{\ga}}{4\ga^2}+ \f{\psi'\pa{1-\ga}}{2\ga} \; \\
\sul{n \in \mathbb{Z}}{} \f{\psi\pa{n-\ga}}{n^2-\ga^2} &=& -\pi \f{\cot\pi\ga}{\ga} \paa{ \psi\pa{1+2\ga}+\pi\cot\pi\ga} \;. \hspace{1cm}
\eeqa
The last of these summation identities is maybe less standard. It follows from an $\eps$ differentiation at $\eps=0$ of Dougall's  formula for sums
 $\Ga$ functions \cite{BatemanHigherTranscendentalFunctions}:
\beq
\sul{n\in\mathbb{Z}}{} \Ga\pab{n+\ga,n-\ga+\eps}{1+n-\ga, 1+n+\ga}= \f{\pi^2}{\sin\pi\ga\, \sin\pi\pa{\eps-\ga} \, \Ga\pa{1-2\ga,1-2\ga-\eps}} \; .
\enq
Finally, using standard properties of the $\psi$ function we get that
\bem
-i\Ga\pab{1-\ga,1-\ga}{1-2\ga,1-2\ga} \pac{ \pa{\Dp{\a_2\, , \eps}^2-\Dp{\a_1 \, ,\eps}^2}f}\pa{0,1,1,0;s}\\
=\ga \pi \paa{\psi\pa{1-\ga}-2\psi\pa{1-2\ga}} +\e{O}\pa{s \ln s} \;\; .
\end{multline}


\begin{thebibliography}{10}

\bibitem{AchiezerKacFormulaforTruncatedWienerHopf}
N.I.~Achiezer, \emph{{"The continuous analogues of some Theorems on Toeplitz
  matrices."}}, Ukrainian Math. J. \textbf{16} (1964), 445--462.

\bibitem{BasorLocalizationThmForToeplitz}
E.L.~Basor, \emph{{"A localization theorem for Toeplitz determinants."}},
  Indiana Univ. Math. J. \textbf{\bf{28}} (1979), 975--983.

\bibitem{BasorTracyGeneralizedFischer-HartwigConjecture}
E.L.~Basor and C.~A. Tracy, \emph{{"The Fisher-Hartwig Conjecture and
  generalizations. Current problems in statistical mechanics"}}, Physica A
  \textbf{\bf{177}} (1991), 167--173.

\bibitem{BasorWidomWienreHopfwithOneFHSingularity}
E.L.~Basor and H.~Widom, \emph{{"Wiener-Hopf determinants with Fisher-Hartwig
  symbols."}}, Oper. Th. \textbf{\bf{147}} (2004), 131--149.

\bibitem{BatemanHigherTranscendentalFunctions}
H.~Bateman, \emph{{Higher transcendental functions}}, McGraw-Hill Books, Malabar, Florida, 1953.

\bibitem{BaxterToeplitzStrongSzego}
G.~Baxter, \emph{{"A norm inequality for a finite section Wiener-Hopf
  equation."}}, Illinois J.Math. \textbf{\bf 7} (1963), 97--103.

\bibitem{BottcherFHToeplitzPureJump}
A.~B\"{o}ttcher, \emph{{"Toeplitz determinants with piecewise continuous generating
  function."}}, Z. Anal. Anw. \textbf{\bf{1}} (1982), 23--39.

\bibitem{BottcherFHConjectureWienerHopf}
A.~B\"{o}ttcher, \emph{{"Wiener-Hopf determinants with rational symbols."}}, Math.
  Nachr. \textbf{\bf{144}} (1989), 39--64.

\bibitem{BottcherSilbWienerHopfAsympAllNuorBarNuZero}
A.~B\"{o}ttcher and B.~Silbermann, \emph{{"Wiener-Hopf determinants with symbols
  having zero of analytic type."}}, Seminar Analysis 1982/83, Inst. f. Math.,
  Akad. Wiss. DDR, Berlin, 1983.

\bibitem{BottcherSilbProofOfFHConjInSomePartCases}
A.~B\"{o}ttcher and B.~Silbermann, \emph{{"Toeplitz matrices and determinants with Fisher-Hartwig
  symbols."}}, J.Funct.Anal. \textbf{\bf{63}} (1985), 178--214.

\bibitem{BottcherSilAnalysisOfToeplitzOperators}
A.~B\"{o}ttcher and B.~Silbermann, \emph{{"Analysis of Toeplitz operators."}}, Springer Verlag, 1990.

\bibitem{BottcherSilberWidomProofFHConjWienerHopfDiscontinuity}
A.~B\"{o}ttcher, B.~Silbermann, and H.~Widom, \emph{{"A continuous analogue of the
  Fisher-Hartwig formula for piecewise continuous symbols."}}, J.Funct.Anal.
  \textbf{\bf{122}} (1994), 222--246.

\bibitem{DeiftOrthPlyAndRandomMatrixRHP}
P.A.~Deift, \emph{{"Orthogonal Polynomials and Random Matrices: A
  Riemann-Hilbert Approach."}}, Courant Lecture Notes 3, New-York University,
  1999.

\bibitem{DeiftItsKrasovskyAsymptoticsofToeplitsHankelWithFHSymbols}
P.A.~Deift, A.R.~Its, and I.~Krasovsky, \emph{{"Asymptotics of Toeplitz,
  Hankel, and Toeplitz+Hankel determinants with Fisher-Hartwig
  singularities."}}, Ann. Math., \textbf{174}, 1243-1299, (2011).

\bibitem{DeiftItsKrasovskyProofOfGeneralizedFHConjectureAndMoreAnnouncmentResu%
lts}
P.A.~Deift, A.R.~Its, and I.~Krasovsky, \emph{{"Toeplitz and Hankel determinants with singularities:
  announcement of results."}}, ArXiV math-FA/ 08092420.

\bibitem{DeiftItsZhouSineKernelOnUnionOfIntervals}
P.A.~Deift, A.R.~Its, and X.~Zhou, \emph{{"A Riemann-Hilbert approach to
  asymptotics problems arising in the theory of random matrix models and also
  in the theory of integrable statistical mechanics."}}, Ann. Math.
  \textbf{\bf{146}} (1997), 149--235.

\bibitem{DeiftKriechMcLaughVenakZhouOrthogonalPlyExponWeights}
P.A.~Deift, T.~Kriecherbauer, K.T.-R.~McLaughlin, S.~Venakides, and X.~Zhou,
  \emph{{"Strong asymptotics of orthogonal polynomials with respect to
  exponential weights."}}, Comm. Pure Appl. Math \textbf{\bf{52}} (1999),
  1491--1552.

\bibitem{DeiftKriechMcLaughVenakZhouOrthogonalPlyVaryingExponWeights}
P.A.~Deift, T.~Kriecherbauer, K.T.-R.~McLaughlin, S.~Venakides, and X.~Zhou, \emph{{"Uniform asymptotics for polynomials orthogonal with respect to
  varying exponential weights and application to universality questions in
  random matrix theory."}}, Comm. Pure Appl. Math \textbf{\bf{52}} (1999),
  1335--1425.

\bibitem{DeiftZhouSteepestDescentForOscillatoryRHP}
P.A.~Deift and X.~Zhou, \emph{{"A steepest descent method for oscillatory
  Riemann-Hilbert problems."}}, Intl. Math. Res. \textbf{\bf{6}} (1997),
  285--299.

\bibitem{EhrhardtAsymptoticBehaviorOfFischerHartwigToeplitzGeneralCase}
T.~Ehrhardt, \emph{{"Toeplitz determinants with several Fisher-Hartwig
  singularities."}}, Ph.D. thesis, {Fakult$\ddot{a}$t f$\ddot{u}$r Mathematik
  der Technischen Universit$\ddot{a}$t Chemnitz}, Chemnitz, Germany, 1997.

\bibitem{EhrhardtSilbermannOnePureFHSingularity}
T.~Ehrhardt and B.~Silbermann, \emph{{"Toeplitz determinants with one
  Fisher-Hartwig singularity."}}, J.Funct. Anal. \textbf{\bf{148}} (1997),
  229--256.

\bibitem{ErdelyiTablesOfIntTransf}
A.~Erdelyi, \emph{{"Tables of Integral Transforms."}}, McGraw-Hill Inc.,US,
  1954.

\bibitem{FischerHartwigTheConjecture}
M.E.~Fisher and R.E.~Hartwig, \emph{{"Toeplitz determinants: some applications,
  theorems and conjectures."}}, Adv. Chem. Phys. \textbf{\bf{15}} (1968),
  333--353.

\bibitem{FokasItsKitaevIsomonodromyPlusRHPforOrthPly}
A.S.~Fokas, A.R.~Its, and A.V.~Kitaev, \emph{{"The Isomonodromy approach to
  Matrix Models in 2D Quantum Gravity."}}, Comm. Math. Phys. \textbf{\bf{147}}
  (1992), 395--430.

\bibitem{GakhovBoundaryValueProblems}
F.D.~Gakhov, \emph{"Boundary value problems"}, Dover edition, General Publishing
  Company, Canada, Toronto, 1990.

\bibitem{HirschmanSzegolimittheoremwithweakerHypothesis}
I.I.~Hirschman, \emph{{"On a theorem of Szeg$\overset{�}{\e{o}}$, Kac and
  Baxter."}}, J. Analyse Math. \textbf{\bf{14}} (1965), 225--234.

\bibitem{IbragimovFinalFormulationOfStrongSzego}
I.A.~Ibragimov, \emph{{"A theorem of Szeg$\overset{�}{\e{o}}$"}}, Mat. Zametki
  \textbf{\bf{3}} (1968), 693--702.

\bibitem{ItsDifferentialMethodForParametrix}
A.R.~Its, \emph{{"Asymptotic behavior of the solutions to the nonlinear
              {S}chr\"odinger equation, and isomonodromic deformations of
              systems of linear differential equations."}}, Dokl. Akad. Nauk SSSR, \textbf{\bf{261}}, 1981, 14-18.


\bibitem{ItsKrasHankelDetAsymForJumpSing}
A.R.~Its and I.~Krasovsky, \emph{{"Hankel determinant and orthogonal
  polynomials for the Gaussian weight with a jump."}}, math.FA/07063192 (2007).

\bibitem{ItsIzerginKorepinSlavnovDifferentialeqnsforCorrelationfunctions}
A.R.~Its, A.G.~Izergin, V.E.~Korepin, and N.A.~Slavonv, \emph{{"Differential
  equations for quantum correlation functions."}}, Int. J. Mod. Physics
  \textbf{\bf{B4}} (1990), 1003--1037.

\bibitem{KacAcheizerTruncatedWienerHopf}
M.~Kac, \emph{{"Toeplitz matrices, translation kernels and related problem in
  probability."}}, Duke Math. J. \textbf{\bf 21} (1954), 501--510.

\bibitem{KarlueIntegralsOf2CHF}
E.~Karlue, \emph{{"Integrals for 2 Hypergeometric functions."}}, J. Phys. A:
  Math. Gen \textbf{\bf 23} (1990), 1969--1971.

\bibitem{KozKitMailSlaTerRHPapproachtoSuperSineKernel}
N.~Kitanine, K.K.~Kozlowski, J.-M.~Maillet, N.~A. Slavnov, and V.~Terras,
  \emph{{"The Riemann-Hilbert approach to a generalized sine kernel and
  applications."}}, Comm. Math. Phys.  \textbf{291}, 691-761, (2009).

\bibitem{KrasovskyToeplitzFHtypeOnArc}
I.~Krasovsky, \emph{{"Asymptotics for Toeplitz determinants on a circular
  arc."}}, math.ph/0401256 (2004).

\bibitem{KrasovskyHankelDetAsymForPowerLike}
I.~Krasovsky, \emph{{"Correlations of the characteristic polynomials in the Gaussian
  Unitary Ensemble or a singular Hankel determinant."}}, Duke Math. J.
  \textbf{\bf{139}} (2007), 581--619.

\bibitem{KuilajaarsMVVUniformAsymptoticsForModifiedJacobiOrthogonalPolynomials}
A.B.J.~Kuijlaars, K.T.-R.~McLaughlin, W.~Van Assche, and M.~Vanlessen,
  \emph{{"The Riemann-Hilbert approach to strong asymptotics for orthogonal
  polynomials on [-1,1]."}}, Advances in Math. \textbf{\bf 188} (2004),
  337--398.

\bibitem{SimonsInfiniteDimensionalDeterminants}
B.~Simons, \emph{{"Notes on infinite dimensional determinants."}}, Adv. Math.
  \textbf{\bf{24}} (1977), 244--273.

\bibitem{SlaterGeneralizedHypergeormetricFunctions}
L.J.~Slater, \emph{{"Generalized Hypergeometric Functions."}}, Cambridge
  University Press, 1966.

\bibitem{WidomSzegoLimitforBlockToeplitz}
H.~Widom, \emph{{"Asymptotic Behavior of Block Toeplitz Matrices and
  Determinants II."}}, Adv. Math \textbf{\bf{21}} (1976), 1--29.

\end{thebibliography}
\end{document}